\documentclass[12pt,twoside]{article}
\usepackage{amssymb}
\usepackage{amsmath}
\usepackage{amsfonts}

\begin{document}
\pagestyle{empty}

\pagestyle{myheadings}
\newcommand\testopari{\sc Barbu \ --- \ Colli \ --- \ Gilardi \ --- \ Marinoschi}
\newcommand\testodispari{\sc Feedback stabilization of a Cahn-Hilliard type system}
\markboth{\testodispari}{\testopari}

%\title{}
%\author{}
%\date{}
%\maketitle

\renewcommand{\theequation}{\arabic{section}.\arabic{equation}}

\thispagestyle{empty}
\begin{center}
{\Large Feedback stabilization of the Cahn-Hilliard type system }

{\Large for phase separation}

\bigskip

Viorel Barbu$^{1}$, Pierluigi Colli$^{2}$, Gianni Gilardi$^{2}$, Gabriela
Marinoschi$^{3,\ast }$

\footnotetext{$^{\ast }$Corresponding author.
\par
\textit{E-mail addresses}: vb41@uaic.ro (V. Barbu), pierluigi.colli@unipv.it
(P. Colli), gianni.gilardi@unipv.it (G. Gilardi),
gabriela.marinoschi@acad.ro (G. Marinoschi)
\par
{}}\bigskip

$^{1}$\textquotedblleft Al. I. Cuza\textquotedblright\ University, Ia\c{s}i,
and Octav Mayer Institute of Mathematics

Bdul. Carol I 8, Ia\c{s}i, Romania

$^{2}$Dipartimento di Matematica \textquotedblleft F.
Casorati\textquotedblright , Universit\`{a} di Pavia and IMATI-CNR, Pavia

via Ferrata 1, 27100 Pavia, Italy

$^{3}$\textquotedblleft Gheorghe Mihoc-Caius Iacob\textquotedblright\
Institute of Mathematical Statistics and Applied Mathematics of the Romanian
Academy Calea 13 Septembrie 13, 050711 Bucharest, Romania

\bigskip
\end{center}

\noindent Abstract. This article is concerned with the internal feedback
stabilization of the phase field system of Cahn-Hilliard type, modeling the
phase separation in a binary mixture. Under suitable assumptions on an
arbitrarily fixed stationary solution, we construct via spectral separation
arguments a feedback controller having its support in an arbitrary open
subset of the space domain, such that the closed loop nonlinear system
exponentially reach the prescribed stationary solution. This feedback
controller has a finite dimensional structure in the state space of
solutions. In particular, every constant stationary solution is admissible.

\medskip



\section{Introduction}

\setcounter{equation}{0}

We consider the celebrated Cahn-Hilliard system \textbf{(}see \cite{CH-1},
\cite{elliot}\textbf{), }which is coupled, following the phase field
approach introduced by Caginalp (see \cite{brokate}, \cite{caginalp}) with
the energy balance equation in order to describe the spontaneous separation
of the components in a binary mixture
\begin{equation}
(\theta +l_{0}\varphi )_{t}-\Delta \theta =0,\text{ in }(0,\infty )\times
\Omega ,  \label{1-0}
\end{equation}%
\begin{equation}
\varphi _{t}-\Delta \mu =0,\text{ in }(0,\infty )\times \Omega ,  \label{2-0}
\end{equation}%
\begin{equation}
\mu =-\nu \Delta \varphi +F^{\prime }(\varphi )-\gamma _{0}\theta ,\text{ in
}(0,\infty )\times \Omega .  \label{3-0}
\end{equation}%
We complete the system with standard homogeneous Neumann boundary conditions
\begin{equation}
\frac{\partial \theta }{\partial \nu }=\frac{\partial \varphi }{\partial \nu
}=\frac{\partial \mu }{\partial \nu }=0,\text{ on }(0,\infty )\times
\partial \Omega  \label{5-0}
\end{equation}%
and with the initial data
\begin{equation}
\theta (0)=\theta _{0},\text{ }\varphi (0)=\varphi _{0},\text{ in }\Omega .
\label{4-0}
\end{equation}%
The equations and conditions (\ref{1-0})-(\ref{4-0}) give rise to the
so-called conserved phase field system, the name being due also to the mass
conservation of $\varphi $, which is obtained by integrating (\ref{1-0}) in
space and time and using the boundary condition for $\mu $ in (\ref{5-0})
and the initial condition for $\varphi $ in (\ref{4-0}). Proper references
on conserved phase field system are \cite{caginalp-1988}, \cite%
{caginalp-1990} by Caginalp, and the recent contributions \cite%
{caginalp-2011} and \cite{miranville-2013}, where a review of models and
results is done as well. We also quote the contributions \cite{colli-1999},
\cite{rocca-2004} in which a conserved phase field model allowing further memory
effects is investigated.

In the system (\ref{1-0})-(\ref{4-0}) the variables $\theta ,$ $\varphi $
and $\mu $ represent the temperature, the order parameter and the chemical
potential, respectively, $\nu $ is the outward normal vector to the
boundary, $l_{0},$ $\gamma _{0}$ are positive constants with some physical
meaning, and $F^{\prime }$ is the derivative of the double-well potential%
\begin{equation}
F(\varphi )=\frac{(\varphi ^{2}-1)^{2}}{4}.  \label{6}
\end{equation}%
The space domain $\Omega $ is an open, bounded connected subset of $\mathbb{R%
}^{d},$ $d=1,2,3,$ with a sufficiently smooth boundary $\Gamma =\partial
\Omega ,$ and the time $t$ runs in $\mathbb{R}^{+}=(0,\infty ).$ This system
has been widely studied in the last decades from several points of view
including existence of attractors and optimal control. A list of recent
references can be found in \cite{cherfils} and \cite{colli}.

In this paper we shall treat the stabilization for the Cahn-Hilliard system
around a stationary solution by two controllers $(u,v)$ having their support
in an open subset $\omega $ of $\Omega ,$ and placed on the right-hand sides
of equations (\ref{1-0})-(\ref{2-0}). By introducing the expression of $\mu $
given by (\ref{3-0}) into (\ref{2-0}) the system to be stabilized reads
\begin{equation}
\varphi _{t}-\Delta \left( -\nu \Delta \varphi +F^{\prime }(\varphi )-\gamma
_{0}\theta \right) =1_{\omega }^{\ast }v,\text{ in }(0,\infty )\times \Omega
,  \label{1}
\end{equation}%
\begin{equation}
(\theta +l_{0}\varphi )_{t}-\Delta \theta =1_{\omega }^{\ast }u,\text{ in }%
(0,\infty )\times \Omega ,  \label{2}
\end{equation}%
\begin{equation}
\frac{\partial \varphi }{\partial \nu }=\frac{\partial (\Delta \varphi )}{%
\partial \nu }=\frac{\partial \theta }{\partial \nu }=0,\text{ on }(0,\infty
)\times \Gamma ,  \label{5}
\end{equation}%
\begin{equation}
\varphi (0)=\varphi _{0},\text{ }\theta (0)=\theta _{0},\text{ in }\Omega .
\label{4}
\end{equation}%
The second boundary condition in (\ref{5}) follows by (\ref{3-0}) and (\ref%
{5-0}).

We specify that the function denoted $1_{\omega }^{\ast }$ is chosen with
the following properties%
\begin{equation}
1_{\omega }^{\ast }\in C_{0}^{\infty }(\Omega ),\text{ \ supp }1_{\omega
}^{\ast }\subset \omega ,\text{ \ }1_{\omega }^{\ast }>0\text{ on }\omega
_{0},\text{ }  \label{1*}
\end{equation}%
where $\omega _{0}$ is an open subset of $\omega .$

The purpose is to stabilize exponentially the solution to (\ref{1})-(\ref{4}%
) around a stationary solution $(\varphi _{\infty },\theta _{\infty })$ of
the uncontrolled system, by means of the feedback control $(v,u)$ expressed
as a function $\mathcal{F}(\varphi ,\theta )$. This turns out to prove that%
\begin{equation}
\lim_{t\rightarrow \infty }(\varphi (t),\theta (t))=(\varphi _{\infty
},\theta _{\infty }),  \label{7}
\end{equation}%
with an exponential rate of convergence, provided that the initial datum $%
(\varphi _{0},\theta _{0})$ is in a suitable neighborhood of $(\varphi
_{\infty },\theta _{\infty }).$

At this point we observe that the set of stationary states of the
uncontrolled system (\ref{1})-(\ref{5}) (for $u=v=0)$ is not empty, because
this may have any constant solution $\theta _{\infty }$ with some constant
or not constant solution $\varphi _{\infty }.$ A discussion concerning the
solutions to the stationary system
\begin{eqnarray}
\nu \Delta ^{2}\varphi _{\infty }-\Delta F^{\prime }(\varphi _{\infty })
&=&0,\text{ in }\Omega ,  \notag \\
-\Delta \theta _{\infty } &=&0,\text{ in }\Omega ,  \label{17'} \\
\frac{\partial \varphi _{\infty }}{\partial \nu } &=&\frac{\partial \Delta
\varphi _{\infty }}{\partial \nu }=\frac{\partial \theta _{\infty }}{%
\partial \nu }=0,\text{ on }\Gamma  \notag
\end{eqnarray}%
is presented in Lemma A1 in Appendix. The result asserts that $\theta
_{\infty }$ is constant and $\varphi _{\infty }\in H^{4}(\Omega )\subset
C^{2}(\overline{\Omega }).$ Also, $\varphi _{\infty }$ may be constant or
not.

It should be mentioned that a simple analysis of the linearized system
around a stationary state $(\varphi _{\infty },\theta _{\infty })$ reveals
that, in general, not all solutions to the stationary system are
asymptotically stable and so, their stabilization via a feedback controller
with support in an arbitrary subset $\omega \subset \Omega $ is of crucial
importance.

The stabilization technique used first in \cite{RT} for parabolic equations
and then in \cite{vb-stab-NS}, \cite{vb-Las-Trig}- \cite{vb-gw-2002}, \cite%
{raymond-2010} for Navier-Stokes equations and nonlinear parabolic systems
is based on the design of the feedback controller as a linear combination of
the unstable modes of the corresponding linearized system.

\subsection{Main result}

All the proofs given in this work converge to the main result of
stabilization which is described below in a few words, for the reader's
convenience. To this end, we briefly introduce some notation and definitions
necessary to give the statement of the theorem. Some of them will be resumed
and explained later, at the appropriate places.

\subparagraph{Functional framework.}

Let us denote
\begin{equation*}
H=L^{2}(\Omega ),\text{ }V=H^{1}(\Omega ),
\end{equation*}%
with the standard scalar products, identify $H$ with its dual space and set $%
V^{\prime }=(H^{1}(\Omega ))^{\prime }.$ Let $A:D(A)\subset H\rightarrow H$
be the linear operator
\begin{equation}
A=-\Delta +I,\text{ \ }D(A)=\left\{ w\in H^{2}(\Omega );\text{ }\frac{%
\partial w}{\partial \nu }=0\text{ on }\Gamma \right\} .  \label{29}
\end{equation}%
The operator $A$ is $m$-accretive on $H$ and so we can define its fractional
powers $A^{\alpha }$, $\alpha \geq 0$ (see e.g., \cite{pazy}, p. 72).\ We
recall that $A^{\alpha }$ is a linear continuous positive and self-adjoint
operator on $H,$ with the domain
\begin{equation*}
D(A^{\alpha })=\{w\in H;\left\Vert A^{\alpha }w\right\Vert _{H}<\infty \}
\end{equation*}%
and the norm
\begin{equation}
\left\Vert w\right\Vert _{D(A^{\alpha })}=\left\Vert A^{\alpha }w\right\Vert
_{H}.  \label{29-0}
\end{equation}%
Moreover, $D(A^{\alpha })\subset H^{2\alpha }(\Omega ),$ with equality if
and only if $2\alpha <3/2$\textbf{.}

Let $F_{l}$ and $\gamma $ be positive constants that will be specified later
and let us denote by $I$ the identity operator$.$ We introduce the
self-adjoint operator $\mathcal{A}:D(\mathcal{A})\subset H\times
H\rightarrow H\times H,$
\begin{equation}
\mathcal{A}=\left[
\begin{array}{cc}
\nu \Delta ^{2}-F_{l}\Delta & \gamma \Delta \\
\gamma \Delta & -\Delta%
\end{array}%
\right] ,  \label{30}
\end{equation}%
having the domain
\begin{equation*}
D(\mathcal{A})=\left\{ w=(y,z)\in H^{2}(\Omega )\times H^{1}(\Omega );\text{
}\mathcal{A}w\in H\times H,\text{ }\frac{\partial y}{\partial \nu }=\frac{%
\partial \Delta y}{\partial \nu }=\frac{\partial z}{\partial \nu }=0\text{
on }\Gamma \right\} .
\end{equation*}%
We denote by $\lambda _{i}$ and $\{(\varphi _{i},\psi _{i})\}_{i\geq 1}$ the
eigenvalues and eigenvectors respectively, of $\mathcal{A}$. Since $\mathcal{%
A}$ is self-adjoint and its resolvent $(\lambda I+\mathcal{A})^{-1}$ is
compact (as seen in a later proof), the eigenvalues are real and there is a
finite number $N$ of nonpositive eigenvalues $\lambda _{i}\leq 0,$ $%
i=1,...,N.$ We introduce the operators $B$ and $B^{\ast }$ ($B^{\ast }$
being the adjoint of $B)$ as%
\begin{equation*}
B:\mathbb{R}^{N}\rightarrow H\times H,\text{ }B^{\ast }:H\times H\rightarrow
\mathbb{R}^{N},
\end{equation*}%
\begin{equation}
BW=\left[
\begin{tabular}{c}
$\sum\limits_{i=1}^{N}1_{\omega }^{\ast }\varphi _{i}w_{i}$ \\
$\sum\limits_{i=1}^{N}1_{\omega }^{\ast }\psi _{i}w_{i}$%
\end{tabular}%
\right] \text{ for all }W=\left[
\begin{tabular}{c}
$w_{1}$ \\
$...$ \\
$w_{N}$%
\end{tabular}%
\right] \in \mathbb{R}^{N},  \label{B}
\end{equation}%
and
\begin{equation}
B^{\ast }q=\left[
\begin{tabular}{c}
$\int_{\Omega }1_{\omega }^{\ast }(\varphi _{1}q_{1}+\psi _{1}q_{2})dx$ \\
$...$ \\
$\int_{\Omega }1_{\omega }^{\ast }(\varphi _{N}q_{1}+\psi _{N}q_{2})dx$%
\end{tabular}%
\right] \text{ for all }q=\left[
\begin{array}{c}
q_{1} \\
q_{2}%
\end{array}%
\right] \in H\times H.  \label{B*}
\end{equation}%
Moreover, let $R$ be a linear positive self-adjoint operator
\begin{equation*}
R:D(A^{1/2})\times D(A^{1/4})\rightarrow H\times H
\end{equation*}%
which is the solution of the algebraic Riccati equation
\begin{equation}
2R\mathcal{A}+RBB^{\ast }R=\left[
\begin{tabular}{ll}
$A^{3}$ & $0$ \\
$0$ & $A^{3/2}$%
\end{tabular}%
\right] .  \label{18'}
\end{equation}%
The existence of a solution $R$ to (\ref{18'}) will be proved later on. Now,
we are ready to present the stabilization result which is the main aim of
our paper. We call the \textit{closed loop system} the system (\ref{1})-(\ref%
{4}) in which the right-hand side $(1_{\omega }^{\ast }v,1_{\omega }^{\ast
}u)$ is replaced by a function depending on $(\varphi ,\theta )$ defined by
the means of $R,$ more exactly
\begin{equation}
(1_{\omega }^{\ast }v,1_{\omega }^{\ast }u)=-BB^{\ast }R(\varphi -\varphi
_{\infty },\alpha _{0}(\theta -\theta _{\infty }+l(\varphi -\varphi _{\infty
})),  \label{control}
\end{equation}%
where the parameter $\alpha _{0}$ is introduced below in (\ref{alpha0-1}).
The theorem below given for the Cahn-Hilliard system in terms of $(\theta
,\varphi )$ is a consequence of Theorem 3.1 in Section 3, which is the main
result of this paper.

\medskip

Let $(\theta _{\infty },\varphi _{\infty })$\ be a solution to the
stationary uncontrolled system\textit{\ }(\ref{1})-(\ref{5})\textit{\ }and
set
\begin{equation*}
\chi _{\infty }:=\left\Vert \nabla \varphi _{\infty }\right\Vert _{L^{\infty
}(\Omega )}+\left\Vert \Delta \varphi _{\infty }\right\Vert _{L^{\infty
}(\Omega )}.
\end{equation*}

\medskip

\noindent \textbf{Theorem 1.1}. \textit{\ There exists }$\chi _{0}>0$\textit{%
\ }(\textit{depending on the problem parameters, the domain and }$\left\Vert
\varphi _{\infty }\right\Vert _{L^{\infty }(\Omega )})$ \textit{such that
the following holds true. If }$\chi _{\infty }\leq \chi _{0},$\textit{\
there exists}{\huge \ }$\rho $\textit{\ such that for all pairs }$(\varphi
_{0},\theta _{0})\in D(A^{1/2})\times D(A^{1/4})$\textit{\ with }%
\begin{equation}
\left\Vert \varphi _{0}-\varphi _{\infty }\right\Vert
_{D(A^{1/2})}+\left\Vert \alpha _{0}(\theta _{0}-\theta _{\infty })+\alpha
_{0}l_{0}(\varphi _{0}-\varphi _{\infty })\right\Vert _{D(A^{1/4})}\leq \rho
,  \label{c3}
\end{equation}%
\textit{the closed loop system }(\ref{1})-(\ref{4}) \textit{with }$%
(1_{\omega }^{\ast }v,1_{\omega }^{\ast }u)$\textit{\ replaced by} (\ref%
{control})\textit{\textbf{\ }has a unique solution }%
\begin{eqnarray}
(\varphi ,\theta ) &\in &C([0,\infty );H\times H)\cap L^{2}(0,\infty
;D(A^{3/2})\times D(A^{3/4}))  \label{c4} \\
&&\cap W^{1,2}(0,\infty ;(D(A^{1/2})\times D(A^{1/4}))^{\prime }),  \notag
\end{eqnarray}%
\textit{which is exponentially stable, that is}%
\begin{eqnarray}
&&\left\Vert \varphi (t)-\varphi _{\infty }\right\Vert
_{D(A^{1/2})}+\left\Vert \alpha _{0}(\theta (t)-\theta _{\infty })+\alpha
_{0}l_{0}(\varphi (t)-\varphi _{\infty })\right\Vert _{D(A^{1/4})}
\label{c5} \\
&\leq &C_{P}e^{-kt}(\left\Vert \varphi _{0}\right\Vert
_{D(A^{1/2})}+\left\Vert \theta _{0}\right\Vert _{D(A^{1/4})}),  \notag
\end{eqnarray}%
\textit{for some positive constants }$k$\textit{\ and }$C_{P}.$

\medskip

In the previous relations the positive constants $k$\ and $C_{P}$\ depend on
$\Omega ,$ the problem parameters and $\left\Vert \varphi _{\infty
}\right\Vert _{L^{\infty }(\Omega )}.$ In addition, $C_{P}$ depends on the
full norm $\left\Vert \varphi _{\infty }\right\Vert _{W^{2,\infty }(\Omega
)}.$

We remark that hypothesis $\chi _{\infty }\leq \chi _{0}$ is trivially
satisfied if $\varphi _{\infty }$ is a constant. Thus, any constant
stationary solution can be stabilized. This is stressed in the following
corollary.

\medskip

\noindent \textbf{Corollary 1.2}. \textit{Assume }$\varphi _{\infty }$%
\textit{\ to be constant.} \textit{Then, there exists }$\rho $\textit{\ such
that for all pairs }$(\varphi _{0},\theta _{0})\in D(A^{1/2})\times
D(A^{1/4})$\textit{\ satisfying }(\ref{c3})\textit{\ the unique solution to
the closed loop system is exponentially stable.}

\subsection{\textbf{A few} \textbf{preliminaries and plan of the paper}}

We prefer to make a function transformation%
\begin{equation}
\sigma =\alpha _{0}(\theta +l_{0}\varphi ),\text{ }  \label{8}
\end{equation}%
with $\alpha _{0}>0$ chosen such that
\begin{equation}
\frac{\gamma _{0}}{\alpha _{0}}=\alpha _{0}l_{0}=:\gamma >0,  \label{alpha0}
\end{equation}%
that is
\begin{equation}
\alpha _{0}=\sqrt{\frac{\gamma _{0}}{l_{0}}}.  \label{alpha0-1}
\end{equation}%
This transformation will give the possibility to work later on with a
self-adjoint operator acting on the linear part of the system. We observe
that if $l_{0}=\gamma _{0}$ (which usually does not occur in the model) we
directly obtain the self-adjoint linear operator.

Writing the system (\ref{1})-(\ref{5}) in the variables $\varphi $ and $%
\sigma $ and using (\ref{alpha0}) and the notation%
\begin{equation}
l:=\gamma _{0}l_{0},  \label{l}
\end{equation}%
we get the equivalent nonlinear system
\begin{equation}
\varphi _{t}+\nu \Delta ^{2}\varphi -\Delta F^{\prime }(\varphi )-l\Delta
\varphi +\gamma \Delta \sigma =1_{\omega }^{\ast }v,\text{ in }(0,\infty
)\times \Omega ,  \label{13}
\end{equation}%
\begin{equation}
\sigma _{t}-\Delta \sigma +\gamma \Delta \varphi =1_{\omega }^{\ast }u,\text{
in }(0,\infty )\times \Omega ,  \label{14}
\end{equation}%
\begin{equation}
\frac{\partial \varphi }{\partial \nu }=\frac{\partial \Delta \varphi }{%
\partial \nu }=\frac{\partial \sigma }{\partial \nu }=0,\text{ in }(0,\infty
)\times \Gamma ,  \label{16}
\end{equation}%
\begin{equation}
\varphi (0)=\varphi _{0},\text{ }\sigma (0)=\sigma _{0}:=\alpha _{0}(\theta
_{0}+l_{0}\varphi _{0}),\text{ in }\Omega ,  \label{15}
\end{equation}%
with the new meaning of $u,$ namely, $\alpha _{0}$ times the old $u.$ We
shall study in fact the stabilization for this transformed system. It is
obvious that if the stabilization $\lim\limits_{t\rightarrow \infty
}(\varphi (t),\sigma (t))=(\varphi _{\infty },\sigma _{\infty })$ is proved
for system (\ref{13})-(\ref{15}), whenever the initial datum $(\varphi
_{0},\sigma _{0})$ is in a neighborhood of $(\varphi _{\infty },\sigma
_{\infty })$, then this implies the stabilization (\ref{7}) for the
corresponding system (\ref{1})-(\ref{4}). We shall discuss this at the
appropriate place. Here, $\sigma _{\infty }$ is defined as $\alpha
_{0}(\theta _{\infty }+l_{0}\varphi _{\infty })$ and in general it can be
constant or not, depending on the same property for $\varphi _{\infty }.$
The stationary system in terms of $\varphi _{\infty }$ and $\sigma _{\infty
} $ reads

\begin{eqnarray}
\nu \Delta ^{2}\varphi _{\infty }-\Delta F^{\prime }(\varphi _{\infty
})-l\Delta \varphi _{\infty }+\gamma \Delta \sigma _{\infty } &=&0,\text{ in
}\Omega ,  \notag \\
-\Delta \sigma _{\infty }+\gamma \Delta \varphi _{\infty } &=&0,\text{ in }%
\Omega ,  \label{17} \\
\frac{\partial \varphi _{\infty }}{\partial \nu } &=&\frac{\partial \Delta
\varphi _{\infty }}{\partial \nu }=\frac{\partial \sigma _{\infty }}{%
\partial \nu }=0,\text{ on }\Gamma .  \notag
\end{eqnarray}%
Next, we rewrite the difference between system (\ref{13})-(\ref{15}) and
system (\ref{17}) by denoting
\begin{equation}
y=\varphi -\varphi _{\infty },\text{ }z=\sigma -\sigma _{\infty },
\label{18}
\end{equation}%
\begin{equation}
y_{0}=\varphi _{0}-\varphi _{\infty }\text{, }z_{0}=\sigma _{0}-\sigma
_{\infty }.  \label{22+0}
\end{equation}%
We have
\begin{equation}
y_{t}+\nu \Delta ^{2}y-\Delta (F^{\prime }(y+\varphi _{\infty })-F^{\prime
}(\varphi _{\infty }))-l\Delta y+\gamma \Delta z=1_{\omega }^{\ast }v,\text{
in }(0,\infty )\times \Omega ,  \label{19}
\end{equation}%
\begin{equation}
z_{t}-\Delta z+\gamma \Delta y=1_{\omega }^{\ast }u,\text{ in }(0,\infty
)\times \Omega ,  \label{20}
\end{equation}%
\begin{equation}
y(0)=y_{0},\text{ }z(0)=z_{0},\text{ in }\Omega ,  \label{21}
\end{equation}%
\begin{equation}
\frac{\partial y}{\partial \nu }=\frac{\partial \Delta y}{\partial \nu }=%
\frac{\partial z}{\partial \nu }=0,\text{ on }(0,\infty )\times \Gamma ,
\label{22}
\end{equation}%
and we shall stabilize it around the state $(0,0)$ for the initial datum $%
(y_{0},z_{0})$ lying in a neighborhood of $(0,0).$

This result is formulated in Theorem 3.1. Since its proof is technical and
long, we shall split parts of it in several Propositions, according to a
strategy following the steps below. As we have already specified, a central
part is the stabilization of the linearized system. We mention that actually
this will be not exactly the linearized system corresponding to the
nonlinear one, but the modified linear system (\ref{23})-(\ref{26}) (given
in Section 2), which is more convenient to be used in this case. Here are
the steps:

(i) Proof of the stabilization of the linear system (\ref{23})-(\ref{26}) by
a finite dimensional control, in Proposition 2.2.

(ii) Introduction and representation of $R,$ calculation of the feedback
control involving the operators $B,$ $B^{\ast }$ and $R,$ and stabilization
of the linear system (\ref{23})-(\ref{26}) by this feedback control in
Section 2, Propositions 2.3, 2.4 and Remark 2.6.

(iii) Proof of the existence of a unique solution to the \textit{nonlinear
closed loop system} (\ref{19})-(\ref{22}) (with $(u,v)$ expressed in terms
of $(y,z)$ by means of $B,$ $B^{\ast },$ $R)$ and stabilization of this
solution, in Section 3, Theorem 3.1. As a matter of fact, this is the main
result of stabilization given for the system in $(y,z).$

(iv) Retrieval of the result presented in Theorem 1.1 for the solution $%
(\varphi ,\theta )$, as a consequence of Theorem 3.1.

\subparagraph{Notation.}

We denote by $C$ or $C_{i},$ $i=1,2,...$ several positive constants possibly
depending on the system structure ($\nu ,l,\gamma )$, domain, space
dimension, and possibly on the norms of $\varphi _{\infty }.$ However, we
shall locally specify the dependence of the constants on $\varphi _{\infty
}. $ A symbol like $C_{\delta }$ with Greek subscripts denotes (possibly
different) constants that depend on the respective parameter, in addition.
Also, we mark precise constants which can be involved in essential proof
arguments by certain (small or capital) letters and specify them in the
text. Whenever no confusion may arise we shall not indicate the arguments of
the functions, for simplicity.

We shall denote by $(\cdot ,\cdot )$ a pair in a product space and by $%
(\cdot ,\cdot )_{X}$ the scalar product in a space $X.$ The norms in $%
L^{\infty }(\Omega )$ and $W^{2,\infty }(\Omega )$ are indicated by $%
\left\Vert \cdot \right\Vert _{\infty }$ and $\left\Vert \cdot \right\Vert
_{2,\infty },$ respectively.

\subparagraph{Tools.}

We repeatedly use the Sobolev embedding theorems
\begin{equation}
\left\Vert w\right\Vert _{L^{2r}(\Omega )}\leq C\left\Vert w\right\Vert
_{H^{\alpha }(\Omega )},\text{ \ }\alpha \geq d\left( \frac{1}{2}-\frac{1}{2r%
}\right) ,\text{ }d\leq 3,  \label{93}
\end{equation}%
(see e.g., \cite{brezis}, p. 285), its consequence%
\begin{equation}
\left\Vert w\right\Vert _{\infty }\leq C\left\Vert w\right\Vert
_{H^{2}(\Omega )},  \label{93-1}
\end{equation}%
and the elementary Young inequality%
\begin{equation}
ab\leq \delta a^{p}+C_{\delta }b^{p^{\prime }},\text{ }a\geq 0,\text{ }b\geq
0,\text{ }\delta >0,\text{ }p\in (1,\infty ),\text{ }\frac{1}{p}+\frac{1}{%
p^{\prime }}=1,  \label{93-2}
\end{equation}%
with $C_{\delta }$ depending on $p,$ besides $\delta .$

Moreover, we shall account for the following inequalities involving the
powers of $A:$

\begin{equation}
\left\Vert A^{\alpha }w\right\Vert _{H}\leq C\left\Vert A^{\alpha
_{1}}w\right\Vert _{H}^{\lambda }\left\Vert A^{\alpha _{2}}w\right\Vert
_{H}^{1-\lambda },\text{ for }\alpha =\lambda \alpha _{1}+(1-\lambda )\alpha
_{2},\text{ }\lambda \in \lbrack 0,1],\text{ }  \label{54+0}
\end{equation}%
\begin{equation}
\left\Vert A^{\alpha }w\right\Vert _{H}\leq C\left\Vert A^{\beta
}w\right\Vert _{H},\text{ if }\alpha <\beta ,  \label{54+1}
\end{equation}

\begin{equation}
\left\Vert A^{\alpha }w\right\Vert _{H^{\beta }(\Omega )}^{2}\leq
C\left\Vert A^{\alpha +\beta /2}w\right\Vert _{H}^{2},  \label{100-0}
\end{equation}%
with $C$ depending on the domain and the exponents.

\section{Stabilization of the linear system}

\setcounter{equation}{0}

In this section we shall deal with the linear system extracted from (\ref{19}%
)-(\ref{22}).

Let $\varphi _{\infty }\in C^{2}(\overline{\Omega })$ be the first component
of a solution to the stationary problem (\ref{17})$.$

We recall that $F$ is defined in (\ref{6}) and we develop $F^{\prime
}(y+\varphi _{\infty })$ in Taylor expansion and rewrite (\ref{19}) as
\begin{equation}
y_{t}+\nu \Delta ^{2}y-\Delta (F^{\prime \prime }(\varphi _{\infty
})y)-l\Delta y+\gamma \Delta z=\Delta F_{r}(y)+1_{\omega }^{\ast }v,
\label{19'}
\end{equation}%
where $F_{r}(y)$ is the rest of second order. Then, we define
\begin{equation}
\overline{F_{\infty }^{\prime \prime }}=\frac{1}{m_{\Omega }}\int_{\Omega
}F^{\prime \prime }(\varphi _{\infty }(\xi ))d\xi =\frac{3}{m_{\Omega }}%
\int_{\Omega }\varphi _{\infty }^{2}(\xi )d\xi -1,  \label{F-1}
\end{equation}%
where $m_{\Omega }$ is the measure of $\Omega .$ Thus, we have
\begin{equation}
F^{\prime \prime }(\varphi _{\infty }(x))=\overline{F_{\infty }^{\prime
\prime }}+g(x),  \label{F-2}
\end{equation}%
where
\begin{equation}
g(x):=\frac{1}{m_{\Omega }}\int_{\Omega }(F^{\prime \prime }(\varphi
_{\infty }(x))-F^{\prime \prime }(\varphi _{\infty }(\xi )))d\xi =\frac{3}{%
m_{\Omega }}\int_{\Omega }(\varphi _{\infty }^{2}(x)-\varphi _{\infty
}^{2}(\xi ))d\xi .  \label{g}
\end{equation}%
Plugging (\ref{F-1}) in (\ref{19'}) we get the following equivalent form of
the nonlinear system (\ref{19})-(\ref{22})
\begin{equation}
y_{t}+\nu \Delta ^{2}y-F_{l}\Delta y+\gamma \Delta z=\Delta
(F_{r}(y)+g(x)y)+1_{\omega }^{\ast }v,\text{ in }(0,\infty )\times \Omega ,
\label{19''}
\end{equation}%
\begin{equation}
z_{t}-\Delta z+\gamma \Delta y=1_{\omega }^{\ast }u,\text{ in }(0,\infty
)\times \Omega ,  \label{20''}
\end{equation}%
\begin{equation}
y(0)=y_{0},\text{ }z(0)=z_{0},\text{ in }\Omega ,  \label{21''}
\end{equation}%
\begin{equation}
\frac{\partial y}{\partial \nu }=\frac{\partial \Delta y}{\partial \nu }=%
\frac{\partial z}{\partial \nu }=0,\text{ in }(0,\infty )\times \Gamma ,
\label{22''}
\end{equation}%
where
\begin{equation}
F_{l}=\overline{F_{\infty }^{\prime \prime }}+l.  \label{Fl}
\end{equation}%
We note that $F_{l}$ also depends on $\Omega $ and on $\left\Vert \varphi
_{\infty }\right\Vert _{L^{2}(\Omega )}.$

Now, we introduce the linear system
\begin{equation}
y_{t}+\nu \Delta ^{2}y-F_{l}\Delta y+\gamma \Delta z=1_{\omega }^{\ast }v,%
\text{ in }(0,\infty )\times \Omega ,  \label{23}
\end{equation}%
\begin{equation}
z_{t}-\Delta z+\gamma \Delta y=1_{\omega }^{\ast }u,\text{ in }(0,\infty
)\times \Omega ,  \label{24}
\end{equation}%
\begin{equation}
y(0)=y_{0},\text{ }z(0)=z_{0},\text{ in }\Omega ,  \label{25}
\end{equation}%
\begin{equation}
\frac{\partial y}{\partial \nu }=\frac{\partial \Delta y}{\partial \nu }=%
\frac{\partial z}{\partial \nu }=0,\text{ in }(0,\infty )\times \Gamma ,
\label{26}
\end{equation}%
which is going to be studied in this Section, while the nonlinear system (%
\ref{19''})-(\ref{22''}) will be the object of Section 3.

Recalling the definition (\ref{30}) of the operator $\mathcal{A}$ we can
write (\ref{23})-(\ref{26}) as

\begin{equation}
\frac{d}{dt}(y(t),z(t))+\mathcal{A}(y(t),z(t))=1_{\omega }^{\ast }U(t),\text{
a.e. }t\in (0,\infty ),  \label{31-0}
\end{equation}%
\begin{equation}
(y(0),z(0))=(y_{0},z_{0}),  \label{32-0}
\end{equation}%
where $U(t)=(v(t),u(t)).$

Since the domain $\Omega $ is regular enough it follows that $D(\mathcal{A}%
)\subset H^{4}(\Omega )\times H^{2}(\Omega ).$ Also, we note that the
operator $\mathcal{A}$ is self-adjoint.

\subsection{Stabilization of the linear system by a finite dimensional
controller}

We set
\begin{equation}
\mathcal{H}=H\times H,\text{ }\mathcal{V}=D(A)\times D(A^{1/2}),\text{ }%
\mathcal{V}^{\prime }=(D(A)\times D(A^{1/2}))^{\prime },  \label{Hrond}
\end{equation}%
and note that $\mathcal{V\subset H}\subset \mathcal{V}^{\prime }$
algebraically and topologically, with compact injections. The scalar
products on $\mathcal{H}$ and $\mathcal{V}$ are defined by
\begin{eqnarray*}
((y,z),(\psi _{1},\psi _{2}))_{\mathcal{H}} &=&\int_{\Omega }(y\psi
_{1}+z\psi _{2})dx,\text{ } \\
((y,z),(\psi _{1},\psi _{2}))_{\mathcal{V}} &=&(\Delta y,\Delta \psi
_{1})_{H}+(y,\psi _{1})_{H}+(\nabla z,\nabla \psi _{2})_{\mathcal{H}%
}+(z,\psi _{2})_{H}.
\end{eqnarray*}%
We note that the second scalar product is equivalent to the one induced on $%
\mathcal{V}$ by the standard scalar product of $H^{2}(\Omega )\times
H^{1}(\Omega ).$

\medskip

\noindent \textbf{Proposition 2.1. }\textit{The operator}\textbf{\ }$%
\mathcal{A}$ \textit{is quasi }$m$\textit{-accretive on }$\mathcal{H},$ that
is $\lambda I+\mathcal{A}$ is $m$-accretive for some $\lambda >0,$ \textit{%
and its resolvent is compact.}

\textit{Let }$(y_{0},z_{0})\in \mathcal{H}$ \textit{and }$(v,u)\in L^{2}(0,T;%
\mathcal{H}).$\textit{\ Then, problem} (\ref{31-0})-(\ref{32-0}) \textit{has
a unique solution }%
\begin{equation}
(y,z)\in C([0,T];\mathcal{H})\cap L^{2}(0,T;\mathcal{V})\cap W^{1,2}([0,T];%
\mathcal{V}^{\prime }),\text{ \textit{for all} }T>0.  \label{32+1}
\end{equation}%
\textit{Moreover, }$(y,z)\in C((0,T];\mathcal{V})$ \textit{and} \textit{we
have the estimate}
\begin{eqnarray}
&&\left\Vert (y(t),z(t)\right\Vert _{\mathcal{H}}^{2}+\left\Vert
(y,z)\right\Vert _{L^{2}(0,T;\mathcal{V)}}^{2}+t\left\Vert
(y(t),z(t)\right\Vert _{\mathcal{V}}^{2}  \label{32-3} \\
&\leq &C\left( \left\Vert (y_{0},z_{0})\right\Vert _{\mathcal{H}%
}^{2}+\int_{0}^{T}\left\Vert 1_{\omega }^{\ast }U(s)\right\Vert _{\mathcal{H}%
}^{2}ds\right) ,\text{ \textit{for all} }t\in (0,T].  \notag
\end{eqnarray}%
\textit{The constant }$C$\textit{\ depends on }$\Omega ,$\textit{\ }$T$%
\textit{, the problem parameters and }$\left\Vert \varphi _{\infty
}\right\Vert _{L^{2}(\Omega )}.$

\medskip

\noindent \textbf{Proof. }The several constants we introduce in the proof
possibly depend on $\left\Vert \varphi _{\infty }\right\Vert _{L^{2}(\Omega
)}.$ We still denote by $\mathcal{A}$ the operator from $\mathcal{V}$ to $%
\mathcal{V}^{\prime }$ given\textbf{\ }by%
\begin{eqnarray}
\left\langle \mathcal{A}(y,z),(\psi _{1},\psi _{2})\right\rangle _{\mathcal{V%
}^{\prime },\mathcal{V}} &=&\int_{\Omega }(\nu \Delta y\cdot \Delta \psi
_{1}+F_{l}\nabla y\cdot \nabla \psi _{1}-\gamma \nabla z\cdot \nabla \psi
_{1})dx  \notag \\
&&+\int_{\Omega }(\nabla z\cdot \nabla \psi _{2}-\gamma \nabla y\cdot \nabla
\psi _{2})dx,  \label{32+2}
\end{eqnarray}%
for any\textbf{\ }$(\psi _{1},\psi _{2})\in \mathcal{V}$. As a matter of
fact this is the extension of $\mathcal{A}$ defined by (\ref{30}). We easily
see that $\mathcal{A}$ is bounded from $\mathcal{V}$ to $\mathcal{V}^{\prime
}$, that is
\begin{equation}
\left\Vert \mathcal{A}(y,z)\right\Vert _{\mathcal{V}^{\prime }}=\sup_{(\psi
_{1},\psi _{2})\in \mathcal{V},\left\Vert (\psi _{1},\psi _{2})\right\Vert _{%
\mathcal{V}}\leq 1}\left\vert \left\langle \mathcal{A}(y,z),(\psi _{1},\psi
_{2})\right\rangle _{\mathcal{V}^{\prime },\mathcal{V}}\right\vert \leq
C\left\Vert (y,z)\right\Vert _{\mathcal{V}},  \label{32'}
\end{equation}%
and that
\begin{equation}
\left\langle \mathcal{A}(y,z),(y,z)\right\rangle _{\mathcal{V}^{\prime },%
\mathcal{V}}\geq C_{1}\left\Vert (y,z)\right\Vert _{\mathcal{V}%
}^{2}-C_{2}\left\Vert (y,z)\right\Vert _{\mathcal{H}}^{2},\text{ for all }%
(y,z)\in \mathcal{V},  \label{32+3}
\end{equation}%
because
\begin{eqnarray*}
&&\left\langle \mathcal{A}(y,z),(y,z)\right\rangle _{\mathcal{V}^{\prime },%
\mathcal{V}}=\int_{\Omega }\left( \nu \left\vert \Delta y\right\vert
^{2}+F_{l}\left\vert \nabla y\right\vert ^{2}-2\gamma \nabla y\cdot \nabla
z+\left\vert \nabla z\right\vert ^{2}\right) dx \\
&\geq &\nu \left\Vert \Delta y\right\Vert _{H}^{2}-(\left\vert
F_{l}\right\vert +2\gamma ^{2})\left\Vert \nabla y\right\Vert _{H}^{2}+\frac{%
1}{2}\left\Vert \nabla z\right\Vert _{H}^{2} \\
&=&\nu \left\Vert y\right\Vert _{D(A)}^{2}+\frac{1}{2}\left\Vert
z\right\Vert _{H^{1}(\Omega )}^{2}-\nu \left\Vert y\right\Vert _{H}^{2}-%
\frac{1}{2}\left\Vert z\right\Vert _{H}^{2}-a_{0}\left\Vert \nabla
y\right\Vert _{H}^{2},
\end{eqnarray*}%
with $a_{0}=\left\vert F_{l}\right\vert +2\gamma ^{2}.$ Since
\begin{equation*}
a_{0}\left\Vert \nabla y\right\Vert _{H}^{2}\leq C\left\Vert y\right\Vert
_{D(A)}\left\Vert y\right\Vert _{H}\leq \frac{\nu }{2}\left\Vert
y\right\Vert _{D(A)}^{2}+\frac{C^{2}}{2\nu }\left\Vert y\right\Vert _{H}^{2}
\end{equation*}%
we obtain (\ref{32+3}). Based on these properties, the restriction of $%
\mathcal{A}$ to $\mathcal{H},$ $\mathcal{A}:D(\mathcal{A})\subset \mathcal{%
H\rightarrow H}$, previously defined, is quasi $m$-accretive on $\mathcal{H}$%
. This means that $\mathcal{A}+C_{2}I:\mathcal{V}\rightarrow \mathcal{V}%
^{\prime }$ is coercive, thus surjective (see \cite{vb-springer-2010}, p.
36).

Let $(y_{0},z_{0})\in \mathcal{H}$ and set $Y=(y,z).$ Since $1_{\omega
}^{\ast }U\in L^{2}(0,T;\mathcal{H})$ and $\mathcal{A}$ is symmetric, by the
Lions existence theorem (see \cite{lions-61}, Thm. 1.1, p. 46) problem (\ref%
{31-0})-(\ref{32-0}) has a unique solution satisfying (\ref{32+1}) and
\begin{equation}
\left\Vert Y(t)\right\Vert _{\mathcal{H}}^{2}+\int_{0}^{T}\left\Vert
Y(t)\right\Vert _{\mathcal{V}}^{2}\leq C\left( \left\Vert Y_{0}\right\Vert _{%
\mathcal{H}}^{2}+\int_{0}^{T}\left\Vert 1_{\omega }^{\ast }U(t))\right\Vert
_{\mathcal{H}}^{2}\right) ,\text{ }\forall t\in \lbrack 0,T].  \label{32-1-0}
\end{equation}%
This is a part of (\ref{32-3}) and the constant $C$ also depends on $T.$
Now, let us multiply formally (\ref{31-0}) by $t\frac{dY}{dt}(t)$ scalarly
in $\mathcal{H}.$ Since $\mathcal{A}$ is symmetric we have%
\begin{equation*}
t\left\Vert \frac{dY}{dt}(t)\right\Vert _{\mathcal{H}}^{2}+\frac{1}{2}\frac{d%
}{dt}t\left\langle \mathcal{A}Y(t),Y(t)\right\rangle _{\mathcal{V}^{\prime },%
\mathcal{V}}=\frac{1}{2}\left\langle \mathcal{A}Y(t),Y(t)\right\rangle _{%
\mathcal{V}^{\prime },\mathcal{V}}+\left( 1_{\omega }^{\ast }U(t),t\frac{dY}{%
dt}(t)\right) _{\mathcal{H}}.
\end{equation*}%
By integrating in time and applying the Young inequality we easily get%
\begin{equation*}
\int_{0}^{t}s\left\Vert \frac{dY}{ds}(s)\right\Vert _{\mathcal{H}%
}^{2}ds+t\left\langle \mathcal{A}Y(t),Y(t)\right\rangle _{\mathcal{V}%
^{\prime },\mathcal{V}}\leq \int_{0}^{t}\left\langle \mathcal{A}%
Y(s),Y(s)\right\rangle _{\mathcal{V}^{\prime },\mathcal{V}%
}ds+\int_{0}^{t}\left\Vert 1_{\omega }^{\ast }U(s)\right\Vert _{\mathcal{H}%
}^{2}ds.
\end{equation*}%
By (\ref{32+3}) and (\ref{32'}) we obtain
\begin{equation*}
C_{1}t\left\Vert Y(t)\right\Vert _{\mathcal{V}}^{2}-C_{2}t\left\Vert
Y(t)\right\Vert _{\mathcal{H}}^{2}\leq C\int_{0}^{t}\left\Vert
Y(s)\right\Vert _{\mathcal{V}}^{2}ds+\int_{0}^{t}\left\Vert 1_{\omega
}^{\ast }U(s)\right\Vert _{\mathcal{H}}^{2}ds,
\end{equation*}%
whence, by (\ref{32-1-0}) we finally get the complete estimate (\ref{32-3})
as claimed.

The above argument shows that $(\lambda I+\mathcal{A})^{-1}$ is well defined
for $\lambda \geq C_{2}.$

Let $(f_{1},f_{2})\in \mathcal{H}$ and denote $(\lambda I+\mathcal{A}%
)^{-1}(f_{1},f_{2})=(y,z).$ It is readily seen that (\ref{32+3}) implies
\begin{equation*}
\left\Vert (y,z)\right\Vert _{\mathcal{V}}^{2}\leq C\left\Vert
(f_{1},f_{2})\right\Vert _{\mathcal{H}}^{2},\text{ for }\lambda \geq C_{2},
\end{equation*}%
and some $C>0,$ whence it follows that $(\lambda I+\mathcal{A})^{-1}(E)$ is
relatively compact whenever $E$ is bounded in $\mathcal{H}.$\hfill $\square $

\medskip

We recall now that $\lambda _{i}$ and $\{(\varphi _{i},\psi _{i})\}_{i\geq
1} $ are the eigenvalues and eigenvectors of $\mathcal{A}$,
\begin{equation}
\mathcal{A}(\varphi _{i},\psi _{i})=\lambda _{i}(\varphi _{i},\psi _{i}),%
\text{ }i=1,2,...  \label{32+8}
\end{equation}%
We notice that one of the coefficients of $\mathcal{A}$ (see (\ref{32+2}))
is $F_{l}.$ Hence the eigenvalues and the eigenfunctions depend also on $%
\left\Vert \varphi _{\infty }\right\Vert _{L^{2}(\Omega )}.$

Since $\mathcal{A}$ is self-adjoint, its eigenvalues are real and
semi-simple, that is, $\mathcal{A}$ is diagonalizable (see \cite{kato}, p.
59). The eigenvectors corresponding to distinct eigenvalues are orthogonal.
Then, orthogonalizing the system $\{(\varphi _{i},\psi _{i})\}_{i}$ in the
space $\mathcal{H}$ we may assume that it is orthonormal and complete in $%
\mathcal{H}.$ Moreover, since the resolvent of $\mathcal{A}$ is compact,
there exists a finite number of nonpositive eigenvalues (see \cite{kato}, p.
187). The sequence $\{\lambda _{i}\}_{i}$ is increasing and every eigenvalue
is repeated according to its order of multiplicity. Let $N$ be the number of
these nonpositive eigenvalues, that is $\lambda _{i}\leq 0,$ for $i=1,...,N.$

Next we show that the system (\ref{31-0})-(\ref{32-0}) can be stabilized by
a finite dimensional control $U$ of the form
\begin{equation}
U(t,x)=(v(t,x),u(t,x))=\sum_{j=1}^{N}w_{j}(t)(\varphi _{j}(x),\psi _{j}(x)).
\label{U-spectral}
\end{equation}%
We rewrite (\ref{31-0})-(\ref{32-0}) as an \textit{open loop linear system}%
\begin{equation}
\frac{d}{dt}(y(t),z(t))+\mathcal{A}(y(t),z(t))=\sum_{j=1}^{N}w_{j}(t)1_{%
\omega }^{\ast }(\varphi _{j},\psi _{j}),\text{ a.e. }t\in (0,\infty ),
\label{31}
\end{equation}%
and take an arbitrary initial condition in $\mathcal{H},$
\begin{equation}
(y(0),z(0))=(y^{0},z^{0}).  \label{32}
\end{equation}

\medskip

\noindent \textbf{Proposition 2.2.} \textit{There exist }$w_{j}\in L^{2}(%
\mathbb{R}^{+}),$\textit{\ }$j=1,...,N,$\textit{\ such that the controller }%
\begin{equation}
U(t,x)=(v(t,x),u(t,x))=\sum_{j=1}^{N}w_{j}(t)(\varphi _{j}(x),\psi _{j}(x)),%
\text{ }t\geq 0,\text{ }x\in \Omega ,  \label{32+10}
\end{equation}%
\textit{stabilizes exponentially system} (\ref{31})-(\ref{32}), \textit{that
is, its solution }$(y,z)$ \textit{satisfies }%
\begin{equation}
\left\Vert y(t)\right\Vert _{H}+\left\Vert z(t)\right\Vert _{H}\leq
C_{P}e^{-kt}\left( \left\Vert y^{0}\right\Vert _{H}+\left\Vert
z^{0}\right\Vert _{H}\right) ,\text{ for all }t\geq 0.  \label{34}
\end{equation}%
\textit{Moreover, we have}
\begin{equation}
\left( \sum_{j=1}^{N}\int_{0}^{\infty }\left\vert w_{j}(t)\right\vert
^{2}dt\right) ^{1/2}\leq C_{P}\left( \left\Vert y^{0}\right\Vert
_{H}+\left\Vert z^{0}\right\Vert _{H}\right) .  \label{34-0}
\end{equation}%
\textit{In both formulas }$C_{P}$\textit{\ and }$k$ \textit{depend on the
problem parameters }($\nu ,$ $\gamma ,$ $l,$ $\Omega )$ \textit{and} $%
\left\Vert \varphi _{\infty }\right\Vert _{L^{2}(\Omega )}.$

\medskip

\noindent \textbf{Proof. }Let\textbf{\ }$T_{0}>0$ be arbitrary but fixed. We
prove that the solution is represented by a sum of two pairs of functions
such that the functions in the first pair vanish at $t=T_{0}$ and the
functions in the second pair decrease exponentially to 0, as $t\rightarrow
\infty .$ We split the proof in two parts.

\subparagraph{Part 1.}

We have the representation
\begin{equation}
(y(t,x),z(t,x))=\sum_{j=1}^{\infty }\xi _{j}(t)(\varphi _{j}(x),\psi
_{j}(x)),\text{ }(t,x)\in (0,\infty )\times \Omega ,  \label{32+9}
\end{equation}%
with $\xi _{j}\in C(\mathbb{R}^{+})$ and
\begin{equation}
\xi _{i}(0)=\xi _{i0}:=\int_{\Omega }(y^{0}\varphi _{j}(x)+z^{0}\psi
_{j}(x))dx,\text{ }i\geq 1.  \label{32-in}
\end{equation}%
We plug the expressions (\ref{32+9}) into (\ref{31})-(\ref{32}), getting%
\begin{equation*}
\sum_{j=1}^{\infty }\left( \xi _{j}^{\prime }(t)(\varphi _{j},\psi _{j})+\xi
_{j}(t)\lambda _{j}(\varphi _{j},\psi _{j})\right) =\sum_{j=1}^{N}1_{\omega
}^{\ast }w_{j}(t)(\varphi _{j},\psi _{j}).
\end{equation*}%
Taking into account that $((\varphi _{i},\psi _{i}),(\varphi _{j},\psi
_{j}))_{H\times H}=\delta _{ij},$ and multiplying scalarly the previous
equation by $(\varphi _{i},\psi _{i})$ in $H\times H,$ we obtain%
\begin{equation}
\xi _{i}^{\prime }+\lambda _{i}\xi _{i}=\sum_{j=1}^{N}w_{j}d_{ij},\text{ \ }%
\xi _{i}(0)=\xi _{i0},\text{ for }i\geq 1,  \label{35}
\end{equation}%
with
\begin{equation}
d_{ij}=\int_{\Omega }1_{\omega }^{\ast }(\varphi _{i}\varphi _{j}+\psi
_{i}\psi _{j})dx,\text{ }j=1,...,N,\text{ }i=1,...  \label{35-1}
\end{equation}%
Notice that $\left\vert d_{ij}\right\vert \leq \sup \left\vert 1_{\omega
}^{\ast }\right\vert .$ First, we discuss the subsystem extracted from (\ref%
{35}) by taking $i=1,...,N.$ It can be written in the form
\begin{equation}
X^{\prime }+MX=DW\text{ and }X(0)=X_{0},  \label{35-0}
\end{equation}%
where
\begin{equation*}
M=\left[
\begin{tabular}{lll}
$\lambda _{1}$ & $...$ & $0$ \\
$...$ & $...$ & $...$ \\
$0$ & $...$ & $\lambda _{N}$%
\end{tabular}%
\right] ,\text{ \ }X=\left[
\begin{tabular}{l}
$\xi _{1}$ \\
$...$ \\
$\xi _{N}$%
\end{tabular}%
\right] ,\text{ \ }X_{0}=\left[
\begin{tabular}{l}
$\xi _{10}$ \\
$...$ \\
$\xi _{N0}$%
\end{tabular}%
\right] ,\text{ \ }
\end{equation*}%
\begin{equation*}
D=\left[
\begin{tabular}{lll}
$d_{11}$ & $...$ & $d_{1N}$ \\
$...$ & $...$ & $...$ \\
$d_{N1}$ & $...$ & $d_{NN}$%
\end{tabular}%
\right] ,\text{ \ }W=\left[
\begin{tabular}{l}
$w_{1}$ \\
$...$ \\
$w_{N}$%
\end{tabular}%
\right] .\text{ }
\end{equation*}%
In the matrix $M$ each $\lambda _{j}$ is repeated according to its order of
multiplicity.

Next, we prove that, for every $T_{0}>0,$ system (\ref{35}) for $i=1,...,N,$
is null controllable on $[0,T_{0}].$ To do that, we first show that the
system $\{\sqrt{1_{\omega }^{\ast }}\varphi _{j},\sqrt{1_{\omega }^{\ast }}%
\psi _{j}\}_{j=1}^{N}$ is linearly independent in $\omega $ (since supp $%
1_{\omega }^{\ast }\subset \omega ).$ We assume that $\sum\limits_{j=1}^{N}%
\alpha _{j}(\sqrt{1_{\omega }^{\ast }}\varphi _{j},\sqrt{1_{\omega }^{\ast }}%
\psi _{j})=0$ in $\omega $ and deduce that $\alpha _{j}=0$ for $j=1,...,N.$
Our assumption reads $\sqrt{1_{\omega }^{\ast }}S=0$ in $\omega $ where $%
S:=\sum\limits_{j=1}^{N}\alpha _{j}(\varphi _{j},\psi _{j}).$ Thus, $S=0$ in
the open set $\omega _{0}$ since here $1_{\omega }^{\ast }>0.$ Now, we
observe that the elliptic system
\begin{equation*}
\nu \Delta ^{2}\varphi -F_{l}\Delta \varphi +\gamma \Delta \psi -\lambda
_{j}\varphi =\gamma \Delta \varphi -\Delta \psi -\lambda _{j}\psi =0
\end{equation*}%
has constant (thus analytic) coefficients, and so any solution $(\varphi
,\psi )$ to it is analytic (see \cite{morrey}). Thus, $S$ is analytic too,
whence $S=0$ in $\Omega .$ This implies that $\alpha _{j}=0$ for $j=1,...,N,$
since the system $\left\{ (\varphi _{j},\psi _{j})\right\} _{j=1}^{N}$ is
linearly independent in $\Omega .$

In conclusion, the system $\{(\sqrt{1_{\omega }^{\ast }}\varphi _{i},\sqrt{%
1_{\omega }^{\ast }}\psi _{i})\}_{i}$ is linearly independent on $\omega $
and so, the determinant of $[d_{ij}]_{i,j}$ is not zero. This implies that
any solution to
\begin{equation}
\sum_{i=1}^{N}d_{ij}p_{i}(t)=0,\text{ }t\in \lbrack 0,T_{0}],\text{ }%
j=1,...,N,  \label{36+0}
\end{equation}%
must be zero, that is $p_{i}(t)=0$ for all $i=1,...,N.$ So, the assumptions
of Lemma A2 in Appendix are trivially satisfied, whence it follows that
there are $w_{i}$ such that $\xi _{i}(T_{0})=0$ for all $i=1,...,N,$ and
\begin{equation}
\left( \int_{0}^{T_{0}}\sum_{i=1}^{N}\left\vert w_{i}(t)\right\vert
^{2}dt\right) ^{1/2}\leq C\sum_{i=1}^{N}\left\vert \xi _{i0}\right\vert ,
\label{36-0}
\end{equation}%
where $\xi _{i},$ $i=1,...,N,$ denote the solution to system (\ref{35}). It
follows by (\ref{32+9}) that $(y(T_{0}),z(T_{0}))=(0,0).$ By (\ref{32+10})
and (\ref{32-in}) we have
\begin{eqnarray}
\left( \int_{0}^{T_{0}}(\left\Vert v(t)\right\Vert _{H}^{2}+\left\Vert
u(t)\right\Vert _{H}^{2})dt\right) ^{1/2} &=&\left(
\int_{0}^{T_{0}}\sum_{i=1}^{N}\left\vert w_{i}(t)\right\vert ^{2}dt\right)
^{1/2}  \label{36+3} \\
&\leq &C\sum_{i=1}^{N}\left\vert \xi _{i0}\right\vert \leq C\left(
\left\Vert y^{0}\right\Vert _{H}+\left\Vert z^{0}\right\Vert _{H}\right) .
\notag
\end{eqnarray}%
From (\ref{35}), by the formula of variation of constants, we have
\begin{equation*}
\xi _{i}(t)=e^{-\lambda _{i}t}\xi
_{i0}+\sum_{j=1}^{N}d_{ij}\int_{0}^{t}e^{-\lambda _{i}(t-s)}w_{j}(s)ds,\text{
}t\geq 0
\end{equation*}%
and recalling (\ref{36-0}), we easily deduce the estimate
\begin{equation}
\left\vert \xi _{i}(t)\right\vert \leq C\left( \left\Vert y^{0}\right\Vert
_{H}+\left\Vert z^{0}\right\Vert _{H}\right) ,\text{ for }t\in \lbrack
0,T_{0}]\text{ and }i=1,...,N.  \label{36+1}
\end{equation}%
The finite dimensional controller steers into the origin, at $t=T_{0},$ the
solution $\{\xi _{j}\}_{j=1}^{N}.$ We extend $w_{i}$ and $\xi _{i}$ by 0 at
the right of $t=T_{0}$, and take as a new controller
\begin{equation}
\widetilde{U}(t)=\left\{
\begin{array}{l}
(v(t),u(t))\text{ for }t<T_{0} \\
0,\text{ \ \ \ \ \ \ \ \ \ \ \ for }t\geq T_{0}%
\end{array}%
\right.  \label{50}
\end{equation}%
and $(y(t),z(t))=(0,0)$ for $t\geq T_{0}.$ For this controller, (\ref{36+3})
remains valid if we replace $T_{0}$ by $+\infty .$ What we have obtained is
exactly (\ref{34-0}).

\subparagraph{Part 2.}

We come back to (\ref{35}) and discuss it for $i\geq N+1.$ We show that it
is stabilized exponentially in origin by the finite dimensional controller (%
\ref{50}). Now, we assume $t>T_{0}$ and recall that $w_{j}(t)=0$ for $%
t>T_{0}.$ In (\ref{35}) we apply again by the formula of variation of
constants and compute an estimate for $\xi _{i}.$ Taking into account the
fact that $\lambda _{i}$ is positive for $i\geq N+1,...,$ and $\lambda
_{N+1}\leq \lambda _{i}$ for $i\geq N+2,$ we have%
\begin{eqnarray*}
\left\vert \xi _{i}(t)\right\vert &\leq &e^{-\lambda _{N+1}t}\left\vert \xi
_{i0}\right\vert +CN\int_{0}^{T_{0}}e^{-\lambda _{N+1}(t-s)}\left\vert
w_{j}(s)\right\vert ds \\
&\leq &e^{-\lambda _{N+1}t}\left\vert \xi _{i0}\right\vert +C_{1}e^{-\lambda
_{N+1}t}\left( \int_{0}^{T_{0}}e^{2\lambda _{N+1}s}ds\right) ^{1/2}\left(
\int_{0}^{T_{0}}\left\vert w_{j}(s)\right\vert ^{2}ds\right) ^{1/2} \\
&\leq &C_{2}e^{-\lambda _{N+1}t}\left( \left\vert \xi _{i0}\right\vert
+\left( \int_{0}^{T_{0}}\left\vert w_{j}(s)\right\vert ^{2}ds\right)
^{1/2}\right) ,\text{ for }t>T_{0},\text{ }i\geq N+1,
\end{eqnarray*}%
where $C,$ $C_{1},C_{2}$ are constants independent of $(y^{0},z^{0}).$
Hence, by (\ref{32+9}), (\ref{36-0}) and by the Bessel inequality we obtain
\begin{equation*}
\left\Vert y(t)\right\Vert _{H}+\left\Vert z(t)\right\Vert _{H}\leq C\left(
\left\Vert y^{0}\right\Vert _{H}+\left\Vert z^{0}\right\Vert _{H}\right) ,%
\text{ for }t>T_{0}.
\end{equation*}%
For $t\leq T_{0}$ we have a similar estimate as in (\ref{36+1}) and all
these together lead to (\ref{34})$,$ as claimed. \hfill $\square $

\subsection{Feedback stabilization of the linear system}

This subsection is devoted to the determination of a feedback controller
(depending on the solution $(y,z))$ which stabilizes exponentially the
solution to (\ref{31})-(\ref{32}). We begin with the study of a minimization
problem which is the key for this purpose.

We recall the operator $A$ defined by (\ref{29}) and consider the quadratic
minimization problem%
\begin{equation}
\Phi (y^{0},z^{0})=\text{ }\underset{W\in L^{2}(0,\infty ;\mathbb{R}^{N})}{%
\text{Min}}\left\{ J(W)=\frac{1}{2}\int_{0}^{\infty }\left( \left\Vert
A^{3/2}y(t)\right\Vert _{H}^{2}+\left\Vert A^{3/4}z(t)\right\Vert
_{H}^{2}+\left\Vert W(t)\right\Vert _{\mathbb{R}^{N}}^{2}\right) dt\right\}
\label{51}
\end{equation}%
subject to (\ref{31})-(\ref{32}). Here $W$ is the function $%
(w_{1},...,w_{N})\in L^{2}(0,\infty ;\mathbb{R}^{N})$ occurring in (\ref{31}%
). We note that $D(\Phi )=\{(y^{0},z^{0})\in H\times H;$ $\Phi
(y^{0},z^{0})<\infty \}.$

In the next proofs we may also refer to $U=(v,u)\in L^{2}(0,\infty ;H\times
H)$ given by (\ref{32+10}), where $\{(\varphi _{j},\psi _{j})\}_{j=1}^{N}$
are the eigenvectors of the operator $\mathcal{A}$ corresponding to the
unstable eigenvalues $\lambda _{j}\leq 0.$

The constants we shall introduce can depend on $\left\Vert \varphi _{\infty
}\right\Vert _{L^{2}(\Omega )}.$

\medskip

\noindent \textbf{Proposition 2.3. }\textit{For each pair }$(y^{0},z^{0})\in
D(A^{1/2})\times D(A^{1/4})$\textit{, problem} (\ref{51}) \textit{has a
unique\ optimal solution}
\begin{equation}
(\{w_{j}^{\ast }\}_{j=1}^{N},y^{\ast },z^{\ast })\in L^{2}(\mathbb{R}^{+};%
\mathbb{R}^{N}))\times L^{2}(\mathbb{R}^{+};D(A^{3/2}))\times L^{2}(\mathbb{R%
}^{+};D(A^{3/4}))  \label{52}
\end{equation}%
\textit{and }%
\begin{equation}
c_{1}\left( \left\Vert A^{1/2}y^{0}\right\Vert _{H}^{2}+\left\Vert
A^{1/4}z^{0}\right\Vert _{H}^{2}\right) \leq \Phi (y^{0},z^{0})\leq
c_{2}\left( \left\Vert A^{1/2}y^{0}\right\Vert _{H}^{2}+\left\Vert
A^{1/4}z^{0}\right\Vert _{H}^{2}\right) .  \label{52+0}
\end{equation}%
\textit{If }$(y^{0},z^{0})\in D(A)\times D(A^{1/2}),$\textit{\ we have}
\begin{eqnarray}
&&\left( \left\Vert Ay^{\ast }(t)\right\Vert _{H}^{2}+\left\Vert
A^{1/2}z^{\ast }(t)\right\Vert _{H}^{2}\right) +\int_{0}^{t}\left(
\left\Vert A^{2}y^{\ast }(s)\right\Vert _{H}^{2}+\left\Vert Az^{\ast
}(s)\right\Vert _{H}^{2}\right) ds  \label{52+1} \\
&\leq &c_{3}\left( \left\Vert Ay^{0}\right\Vert _{H}^{2}+\left\Vert
A^{1/2}z^{0}\right\Vert _{H}^{2}\right) ,\text{ \textit{for all} }t\geq 0,
\notag
\end{eqnarray}%
\textit{where} $c_{1},$ $c_{2}$\textit{, }$c_{3}$ \textit{are positive
constants} (\textit{depending on }$\Omega ,$ \textit{the problem parameters
and the quantity} $\left\Vert \varphi _{\infty }\right\Vert _{L^{2}(\Omega )}).$

\medskip

\noindent \textbf{Proof. }For all $(y^{0},z^{0})\in H\times H,$ it follows
by Proposition 2.2 that there exist $w_{j}\in L^{2}(\mathbb{R}^{+})$ such
that (\ref{31})-(\ref{32}) has a solution with the properties (\ref{34})-(%
\ref{34-0}).

First, we rewrite (\ref{31}) by calculating the operator $\mathcal{A}$ given
by (\ref{30}) in terms of the operator $A=-\Delta +I.$ We get
\begin{equation}
\mathcal{A(}y,z)=\left[
\begin{array}{c}
\nu A^{2}y+(F_{l}-2\nu )Ay-(F_{l}-\nu )y-\gamma Az+\gamma z \\
-\gamma Ay+\gamma y+Az-z%
\end{array}%
\right] .  \label{53}
\end{equation}%
Also, we recall the interpolation relations (\ref{54+0})-(\ref{54+1}).

Now, let $(y^{0},z^{0})\in D(A^{1/2})\times D(A^{1/4}).$ We multiply (\ref%
{31}), where for the moment the right-hand side is written for simplicity $%
(1_{\omega }^{\ast }v,1_{\omega }^{\ast }u),$ by $(Ay(t),A^{1/2}z(t))$
scalarly in $H\times H$ and obtain
\begin{eqnarray}
&&\frac{1}{2}\frac{d}{dt}\left( \left\Vert A^{1/2}y(t)\right\Vert
_{H}^{2}+\left\Vert A^{1/4}z(t)\right\Vert _{H}^{2}\right) +\nu \left\Vert
A^{3/2}y(t)\right\Vert _{H}^{2}+\left\Vert A^{3/4}z(t)\right\Vert _{H}^{2}
\notag \\
&=&-(F_{l}-2\nu )(Ay(t),Ay(t))_{H}+(F_{l}-\nu )(y(t),Ay(t))_{H}+\gamma
(Az(t),Ay(t))_{H}  \notag \\
&&-\gamma (z(t),Ay(t))_{H}+\gamma (Ay(t),A^{1/2}z(t))_{H}-\gamma
(y(t),A^{1/2}z(t))_{H}+(z(t),A^{1/2}z(t))_{H}  \notag \\
&&+\int_{\Omega }1_{\omega }^{\ast }v(t)Ay(t)dx+\int_{\Omega }1_{\omega
}^{\ast }u(t)A^{1/2}z(t)dx,\text{ a.e. }t>0.  \label{54}
\end{eqnarray}%
Next, we use the interpolation properties (\ref{54+0}), (\ref{54+1}) and the
Young formula for the following terms:
\begin{eqnarray*}
-(F_{l}-2\nu )(Ay(t),Ay(t))_{H} &\leq &\left\vert F_{l}-2\nu \right\vert
\left\Vert Ay(t)\right\Vert _{H}^{2}\leq C\left( \left\Vert
A^{3/2}y(t)\right\Vert _{H}^{2/3}\left\Vert A^{0}y(t)\right\Vert
_{H}^{1/3}\right) ^{2} \\
&\leq &\delta \left\Vert A^{3/2}y(t)\right\Vert _{H}^{2}+C_{\delta
}\left\Vert y(t)\right\Vert _{H}^{2},
\end{eqnarray*}%
\begin{eqnarray*}
(F_{l}-\nu )(y(t),Ay(t))_{H} &\leq &\left\vert F_{l}-\nu \right\vert
\left\Vert y(t)\right\Vert _{H}\left\Vert Ay(t)\right\Vert _{H} \\
&\leq &C\left\Vert y(t)\right\Vert _{H}\left( \left\Vert
A^{3/2}y(t)\right\Vert _{H}^{2/3}\left\Vert y(t)\right\Vert _{H}^{1/3}\right)
\\
&=&C\left\Vert A^{3/2}y(t)\right\Vert _{H}^{2/3}\left\Vert y(t)\right\Vert
_{H}^{4/3}\leq \delta \left\Vert A^{3/2}y(t)\right\Vert _{H}^{2}+C_{\delta
}\left\Vert y(t)\right\Vert _{H}^{2},
\end{eqnarray*}%
\begin{eqnarray*}
&&\gamma (Az(t),Ay(t))_{H}=\gamma (A^{1/2}z(t),A^{3/2}y(t))_{H}\leq
C\left\Vert A^{1/2}z(t)\right\Vert _{H}\left\Vert A^{3/2}y(t)\right\Vert _{H}
\\
&\leq &\frac{C}{\delta }\left\Vert A^{1/2}z(t)\right\Vert _{H}^{2}+\delta
\left\Vert A^{3/2}y(t)\right\Vert _{H}^{2}\leq \frac{C}{\delta }\left\Vert
A^{3/4}z(t)\right\Vert _{H}^{4/3}\left\Vert z(t)\right\Vert
_{H}^{2/3}+\delta \left\Vert A^{3/2}y(t)\right\Vert _{H}^{2} \\
&\leq &\delta \left\Vert A^{3/4}z(t)\right\Vert _{H}^{2}+C_{\delta
}\left\Vert z(t)\right\Vert _{H}^{2}+\delta \left\Vert
A^{3/2}y(t)\right\Vert _{H}^{2},
\end{eqnarray*}%
\begin{eqnarray*}
\gamma (z(t),Ay(t))_{H} &\leq &\gamma \left\Vert z(t)\right\Vert
_{H}\left\Vert Ay(t)\right\Vert _{H}\leq \left\Vert Ay(t)\right\Vert
_{H}^{2}+\gamma ^{2}\left\Vert z(t)\right\Vert _{H}^{2} \\
&\leq &\delta \left\Vert A^{3/2}y(t)\right\Vert _{H}^{2}+C_{\delta
}\left\Vert y(t)\right\Vert _{H}^{2}+C\left\Vert z(t)\right\Vert _{H}^{2},
\end{eqnarray*}%
\begin{eqnarray*}
\gamma (Ay(t),A^{1/2}z(t))_{H} &\leq &\gamma \left\Vert Ay(t)\right\Vert
_{H}\left\Vert A^{1/2}z(t)\right\Vert _{H}\leq \frac{C}{\delta }\left\Vert
Ay(t)\right\Vert _{H}^{2}+\delta \left\Vert A^{1/2}z(t)\right\Vert _{H}^{2}
\\
&\leq &\delta \left\Vert A^{3/2}y(t)\right\Vert _{H}^{2}+C_{\delta
}\left\Vert y(t)\right\Vert _{H}^{2}+C_{1}\delta \left\Vert
A^{3/4}z(t)\right\Vert _{H}^{2},
\end{eqnarray*}%
\begin{equation*}
\gamma (y(t),A^{1/2}z(t))_{H}\leq \delta \left\Vert A^{1/2}z(t)\right\Vert
_{H}^{2}+C_{\delta }\left\Vert y(t)\right\Vert _{H}^{2}\leq C_{1}\delta
\left\Vert A^{3/4}z(t)\right\Vert _{H}^{2}+C_{\delta }\left\Vert
y(t)\right\Vert _{H}^{2},
\end{equation*}%
\begin{equation*}
(z(t),A^{1/2}z(t))_{H}\leq \delta \left\Vert A^{1/2}z(t)\right\Vert
_{H}^{2}+C_{\delta }\left\Vert z(t)\right\Vert _{H}^{2}\leq C_{1}\delta
\left\Vert A^{3/4}z(t)\right\Vert _{H}^{2}+C_{\delta }\left\Vert
z(t)\right\Vert _{H}^{2},
\end{equation*}%
\begin{equation*}
\int_{\Omega }1_{\omega }^{\ast }v(t)Ay(t)dx\leq C\left\Vert v(t)\right\Vert
_{H}\left\Vert Ay(t)\right\Vert _{H}\leq C_{\delta }\left\Vert
v(t)\right\Vert _{H}^{2}+\delta \left\Vert A^{3/2}y(t)\right\Vert _{H}^{2},
\end{equation*}%
\begin{equation*}
\int_{\Omega }1_{\omega }^{\ast }u(t)A^{1/2}z(t)dx\leq C\left\Vert
u(t)\right\Vert _{H}\left\Vert A^{1/2}z(t)\right\Vert _{H}\leq C_{\delta
}\left\Vert u(t)\right\Vert _{H}^{2}+C_{1}\delta \left\Vert
A^{3/4}z(t)\right\Vert _{H}^{2}.
\end{equation*}%
Plugging all these in (\ref{54}), choosing $\delta $ small enough and
recalling (\ref{34}), we obtain
\begin{eqnarray*}
&&\frac{1}{2}\frac{d}{dt}\left( \left\Vert A^{1/2}y(t)\right\Vert
_{H}^{2}+\left\Vert A^{1/4}z(t)\right\Vert _{H}^{2}\right) +\frac{\nu }{2}%
\left\Vert A^{3/2}y(t)\right\Vert _{H}^{2}+\frac{1}{2}\left\Vert
A^{3/4}z(t)\right\Vert _{H}^{2} \\
&\leq &C(\left\Vert y(t)\right\Vert _{H}^{2}+\left\Vert z(t)\right\Vert
_{H}^{2}+\left\Vert u(t)\right\Vert _{H}^{2}+\left\Vert v(t)\right\Vert
_{H}^{2}) \\
&\leq &C\left\{ e^{-kt}(\left\Vert y^{0}\right\Vert _{H}^{2}+\left\Vert
z^{0}\right\Vert _{H}^{2})+\left\Vert u(t)\right\Vert _{H}^{2}+\left\Vert
v(t)\right\Vert _{H}^{2}\right\} .
\end{eqnarray*}%
Integrating in time and using (\ref{34-0}) we get
\begin{eqnarray}
&&\left\Vert A^{1/2}y(t)\right\Vert _{H}^{2}+\left\Vert
A^{1/4}z(t)\right\Vert _{H}^{2}  \label{61} \\
&&+C_{\nu }\int_{0}^{t}\left( \left\Vert A^{3/2}y(s)\right\Vert
_{H}^{2}+\left\Vert A^{3/4}z(s)\right\Vert _{H}^{2}+\left\Vert
u(s)\right\Vert _{H}^{2}+\left\Vert v(s)\right\Vert _{H}^{2}\right) ds
\notag \\
&\leq &\left\Vert A^{1/2}y^{0}\right\Vert _{H}^{2}+\left\Vert
A^{1/4}z^{0}\right\Vert _{H}^{2}+\frac{C}{k}(1-e^{-kt})(\left\Vert
y^{0}\right\Vert _{H}^{2}+\left\Vert z^{0}\right\Vert _{H}^{2})+C(\left\Vert
y^{0}\right\Vert _{H}^{2}+\left\Vert z^{0}\right\Vert _{H}^{2})  \notag \\
&\leq &C\left( \left\Vert A^{1/2}y^{0}\right\Vert _{H}^{2}+\left\Vert
A^{1/4}z^{0}\right\Vert _{H}^{2}\right) ,\text{ for all }t>0,  \notag
\end{eqnarray}%
where $C_{\nu }=\min \left\{ \nu ,1\right\} $ and $C$ denotes constants
depending on the problem parameters $\nu ,\gamma ,F_{l}.$ From here we
deduce that
\begin{eqnarray}
&&\int_{0}^{\infty }\left( \left\Vert A^{3/2}y(t)\right\Vert
_{H}^{2}+\left\Vert A^{3/4}z(t)\right\Vert _{H}^{2}+\left\Vert
u(t)\right\Vert _{H}^{2}+\left\Vert v(t)\right\Vert _{H}^{2}\right) dt
\label{62} \\
&\leq &C(\left\Vert A^{1/2}y^{0}\right\Vert _{H}^{2}+\left\Vert
A^{1/4}z^{0}\right\Vert _{H}^{2})\leq c_{2}(\left\Vert
A^{1/2}y^{0}\right\Vert _{H}^{2}+\left\Vert A^{1/4}z^{0}\right\Vert
_{H}^{2}),  \notag
\end{eqnarray}%
where $c_{2}$ depends on the problem parameters. This is the right
inequality in (\ref{52+0}).

Now, we take in (\ref{51}) a minimizing sequence $\{W^{n}\}_{n\geq 1},$ $%
W^{n}=(w_{1}^{n},...,w_{N}^{n})$ such that
\begin{equation*}
(u_{n}(t),v_{n}(t))=\sum\limits_{j=1}^{N}w_{j}^{n}(t)(\varphi _{j},\psi
_{j}).
\end{equation*}%
We can assume that
\begin{equation}
d\leq \frac{1}{2}\int_{0}^{\infty }\left( \left\Vert
A^{3/2}y_{n}(t)\right\Vert _{H}^{2}+\left\Vert A^{3/4}z_{n}(t)\right\Vert
_{H}^{2}+\left\Vert W^{n}(t)\right\Vert _{\mathbb{R}^{N}}^{2}\right) dt\leq
d+\frac{1}{n},  \label{63}
\end{equation}%
where $d$ is the positive infimum of $J(W)$ in (\ref{51}) and $(y_{n},z_{n})$
is the solution to (\ref{31})-(\ref{32}) corresponding to $W^{n}.$ By (\ref%
{63}) we have a subsequence $\{n\rightarrow \infty \}$ such that
\begin{equation*}
w_{j}^{n}\rightarrow w_{j}^{\ast }\text{ weakly in }L^{2}(\mathbb{R}^{+}),%
\text{ }j=1,...,N,
\end{equation*}%
\begin{equation*}
y_{n}\rightarrow y^{\ast }\text{ weakly in }L^{2}(\mathbb{R}^{+};D(A^{3/2})),
\end{equation*}%
\begin{equation*}
z_{n}\rightarrow z^{\ast }\text{ weakly in }L^{2}(\mathbb{R}^{+};D(A^{3/4})).
\end{equation*}%
Also, by (\ref{31}) we have
\begin{equation*}
\frac{d}{dt}(y_{n},z_{n})\rightarrow \frac{d}{dt}(y,z)\text{ weakly in }%
L^{2}(\mathbb{R}^{+};\mathcal{V}^{\prime }),
\end{equation*}%
where $\mathcal{V}^{\prime }$ is defined in (\ref{Hrond}). Since $%
(u_{n}(t),v_{n}(t))=\sum\limits_{j=1}^{N}w_{j}^{n}(t)(\varphi _{j},\psi
_{j}) $ we have
\begin{equation*}
(u_{n},v_{n})\rightarrow (u^{\ast },v^{\ast
})=\sum\limits_{j=1}^{N}w_{j}^{\ast }(t)(\varphi _{j},\psi _{j})\text{
weakly in }L^{2}(\mathbb{R}^{+};H\times H).
\end{equation*}%
Thus, $(y^{\ast },z^{\ast })$ solves the system (\ref{31})-(\ref{32})
corresponding to $W^{\ast }:=(w_{1}^{\ast },...,w_{N}^{\ast }).$ Moreover,
passing to the limit in (\ref{63}) we get on the basis of the weakly lower
semicontinuity of $J$ that $J(W^{\ast })=d.$

The uniqueness follows by the fact that $J$ is strictly convex and the state
system is linear.

Moreover, by (\ref{54}) we can write%
\begin{eqnarray*}
&&\int_{0}^{t}\left( \nu \left\Vert A^{3/2}y(s)\right\Vert
_{H}^{2}+\left\Vert A^{3/4}z(s)\right\Vert _{H}^{2}\right) ds \\
&=&\frac{1}{2}\left( \left\Vert A^{1/2}y^{0}\right\Vert _{H}^{2}+\left\Vert
A^{1/4}z^{0}\right\Vert _{H}^{2}\right) -\frac{1}{2}\left( \left\Vert
A^{1/2}y(t)\right\Vert _{H}^{2}+\left\Vert A^{1/4}z(t)\right\Vert
_{H}^{2}\right) \\
&&-(F_{l}-2\nu )\int_{0}^{t}(Ay(s),Ay(s))_{H}ds+(F_{l}-\nu
)\int_{0}^{t}(y(s),Ay(s))_{H}ds \\
&&+\gamma \int_{0}^{t}(Az(s),Ay(s))_{H}ds \\
&&-\gamma \int_{0}^{t}(z(s),Ay(s))_{H}ds+\gamma
\int_{0}^{t}(Ay(s),A^{1/2}z(s))_{H}ds-\gamma
\int_{0}^{t}(y(s),A^{1/2}z(s))_{H}ds \\
&&+\int_{0}^{t}(z(s),A^{1/2}z(s))_{H}ds+\int_{0}^{t}\int_{\Omega }1_{\omega
}^{\ast }v(s)Ay(s)dxdt+\int_{0}^{t}\int_{\Omega }1_{\omega }^{\ast
}u(s)A^{1/2}z(s)dxds.
\end{eqnarray*}%
We are going to derive a basic inequality by arguing as we did for all terms
on the right-hand side in (\ref{54}) in order to get (\ref{61}), but
suitably changing $\delta $ and $C_{\delta }$ in the use of the Young
inequality. For instance, we have
\begin{eqnarray*}
\left\vert (F_{l}-2\nu )(Ay(t),Ay(t))_{H}\right\vert &\leq &\left\vert
F_{l}-2\nu \right\vert \left\Vert Ay(t)\right\Vert _{H}^{2}\leq C\left(
\left\Vert A^{3/2}y(t)\right\Vert _{H}^{2/3}\left\Vert A^{0}y(t)\right\Vert
_{H}^{1/3}\right) ^{2} \\
&\leq &C_{\delta }\left\Vert A^{3/2}y(t)\right\Vert _{H}^{2}+\delta
\left\Vert y(t)\right\Vert _{H}^{2},
\end{eqnarray*}%
which implies
\begin{equation*}
-(F_{l}-2\nu )(Ay(t),Ay(t))_{H}\geq -C_{\delta }\left\Vert
A^{3/2}y(t)\right\Vert _{H}^{2}-\delta \left\Vert y(t)\right\Vert _{H}^{2}.
\end{equation*}%
By treating all the terms in the same way we arrive at
\begin{eqnarray*}
&&\int_{0}^{t}\left( \nu \left\Vert A^{3/2}y(s)\right\Vert
_{H}^{2}+\left\Vert A^{3/4}z(s)\right\Vert _{H}^{2}\right) ds \\
&\geq &\frac{1}{2}\left( \left\Vert A^{1/2}y^{0}\right\Vert
_{H}^{2}+\left\Vert A^{1/4}z^{0}\right\Vert _{H}^{2}\right) -\frac{1}{2}%
\left( \left\Vert A^{1/2}y(t)\right\Vert _{H}^{2}+\left\Vert
A^{1/4}z(t)\right\Vert _{H}^{2}\right) \\
&&-C_{\delta }\int_{0}^{t}\left( \left\Vert A^{3/2}y(s)\right\Vert
_{H}^{2}+\left\Vert A^{3/4}z(s)\right\Vert _{H}^{2}\right) ds-C_{1}\delta
\int_{0}^{t}(\left\Vert y(s)\right\Vert _{H}^{2}+\left\Vert z(s)\right\Vert
_{H}^{2})ds \\
&&-\delta \int_{0}^{t}(\left\Vert u(s)\right\Vert _{H}^{2}+\left\Vert
v(s)\right\Vert _{H}^{2})ds,
\end{eqnarray*}%
where $C_{1}>0$ could be computed. By relying on (\ref{34}) and (\ref{34-0})
we obtain
\begin{eqnarray*}
&&(C_{\delta }+\max \{\nu ,1\})\int_{0}^{t}\left( \left\Vert
A^{3/2}y(s)\right\Vert _{H}^{2}+\left\Vert A^{3/4}z(s)\right\Vert
_{H}^{2}\right) ds \\
&\geq &\frac{1}{2}\left( \left\Vert A^{1/2}y^{0}\right\Vert
_{H}^{2}+\left\Vert A^{1/4}z^{0}\right\Vert _{H}^{2}\right) -\frac{1}{2}%
\left( \left\Vert A^{1/2}y(t)\right\Vert _{H}^{2}+\left\Vert
A^{1/4}z(t)\right\Vert _{H}^{2}\right) \\
&&-2C_{1}\delta C_{P}^{2}\int_{0}^{t}e^{-2ks}ds(\left\Vert y^{0}\right\Vert
_{H}^{2}+\left\Vert z^{0}\right\Vert _{H}^{2})-2\delta C_{P}^{2}(\left\Vert
y^{0}\right\Vert _{H}^{2}+\left\Vert z^{0}\right\Vert _{H}^{2}).
\end{eqnarray*}%
Computing the integral, using (\ref{54+1}) and choosing $\delta $ small
enough we get%
\begin{eqnarray}
&&(C+\max \{\nu ,1\})\int_{0}^{t}\left( \left\Vert A^{3/2}y(s)\right\Vert
_{H}^{2}+\left\Vert A^{3/4}z(s)\right\Vert _{H}^{2}\right) ds  \label{250} \\
&\geq &\frac{1}{4}\left( \left\Vert A^{1/2}y^{0}\right\Vert
_{H}^{2}+\left\Vert A^{1/4}z^{0}\right\Vert _{H}^{2}\right) -\frac{1}{2}%
\left( \left\Vert A^{1/2}y(t)\right\Vert _{H}^{2}+\left\Vert
A^{1/4}z(t)\right\Vert _{H}^{2}\right) .  \notag
\end{eqnarray}%
Since the last term on the right-hand side is a continuous $L^{1}$ function,
one can take a sequence $t_{j}\nearrow \infty $ such that
\begin{equation*}
\left\Vert A^{1/2}y(t_{j})\right\Vert _{H}^{2}+\left\Vert
A^{1/4}z(t_{j})\right\Vert _{H}^{2}\rightarrow 0.
\end{equation*}%
Passing to the limit in (\ref{250}) along such a sequence we obtain
\begin{equation*}
\int_{0}^{\infty }\left( \left\Vert A^{3/2}y(s)\right\Vert
_{H}^{2}+\left\Vert A^{3/4}z(s)\right\Vert _{H}^{2}\right) ds\geq
c_{1}\left( \left\Vert A^{1/2}y^{0}\right\Vert _{H}^{2}+\left\Vert
A^{1/4}z^{0}\right\Vert _{H}^{2}\right) ,
\end{equation*}%
where $c_{1}>0$ depends only on the problem parameters and $\left\Vert
\varphi _{\infty }\right\Vert _{L^{2}(\Omega )}.$

This relation written for the optimal pair $(W^{\ast },(y^{\ast },z^{\ast
})) $ implies that
\begin{eqnarray*}
\Phi (y^{0},z^{0}) &=&\frac{1}{2}\int_{0}^{\infty }\left( \left\Vert
A^{3/2}y^{\ast }(t)\right\Vert _{H}^{2}+\left\Vert A^{3/4}z^{\ast
}(t)\right\Vert _{H}^{2}+\left\Vert W^{\ast }(t)\right\Vert _{\mathbb{R}%
^{N}}^{2}\right) dt \\
&\geq &c_{1}\left( \left\Vert A^{1/2}y^{0}\right\Vert _{H}^{2}+\left\Vert
A^{1/4}z^{0}\right\Vert _{H}^{2}\right) ,
\end{eqnarray*}%
that is the left inequality in (\ref{52+0}).

Relation (\ref{62}), valid also for the optimal pair, leads to the
right-hand side of (\ref{52+0}).

The next calculation will be done in view of proving (\ref{52+1}).

We recall (\ref{53}) and multiply (\ref{31}) by $(A^{2}y,\alpha Az)$
scalarly in $H\times H,$ with $\alpha $ a positive number that will
specified later. We obtain
\begin{eqnarray}
&&\frac{1}{2}\frac{d}{dt}\left( \left\Vert Ay(t)\right\Vert _{H}^{2}+\alpha
\left\Vert A^{1/2}z(t)\right\Vert _{H}^{2}\right) +\nu \left\Vert
A^{2}y(t)\right\Vert _{H}^{2}+\alpha \left\Vert Az(t)\right\Vert _{H}^{2}
\label{63+1} \\
&=&-(F_{l}-2\nu )(Ay(t),A^{2}y(t))_{H}+(F_{l}-\nu
)(y(t),A^{2}y(t))_{H}+\gamma (Az(t),A^{2}y(t))_{H}  \notag \\
&&-\gamma (z(t),A^{2}y(t))_{H}+\alpha (z(t),Az(t))_{H}+\alpha \gamma
(Ay(t),Az(t))_{H}-\alpha \gamma (y(t),Az(t))_{H}  \notag \\
&&+\int_{\Omega }1_{\omega }^{\ast }v(t)A^{2}y(t)dx+\alpha \int_{\Omega
}1_{\omega }^{\ast }u(t)Az(t)dx.  \notag
\end{eqnarray}%
As previously, we have
\begin{eqnarray*}
-(F_{l}-2\nu )(Ay(t),A^{2}y(t))_{H} &\leq &\delta \left\Vert
A^{2}y(t)\right\Vert _{H}^{2}+C_{\delta }\left\Vert Ay(t)\right\Vert _{H}^{2}
\\
&\leq &\delta \left\Vert A^{2}y(t)\right\Vert _{H}^{2}+\delta \left\Vert
A^{3/2}y(t)\right\Vert _{H}^{2}+C_{\delta }\left\Vert y(t)\right\Vert
_{H}^{2},
\end{eqnarray*}%
\begin{equation*}
(F_{l}-\nu )(y(t),A^{2}y(t))_{H}\leq \delta \left\Vert A^{2}y(t)\right\Vert
_{H}^{2}+C_{\delta }\left\Vert y(t)\right\Vert _{H}^{2},
\end{equation*}%
\begin{equation*}
\gamma (Az(t),A^{2}y(t))_{H}\leq \delta \left\Vert A^{2}y(t)\right\Vert
_{H}^{2}+C_{\delta }\left\Vert Az(t)\right\Vert _{H}^{2},
\end{equation*}%
\begin{equation*}
-\gamma (z(t),A^{2}y(t))_{H}\leq \delta \left\Vert A^{2}y(t)\right\Vert
_{H}^{2}+C_{\delta }\left\Vert z(t)\right\Vert _{H}^{2},
\end{equation*}%
\begin{equation*}
\alpha (z(t),Az(t))_{H}\leq \frac{\alpha }{8}\left\Vert Az(t)\right\Vert
_{H}^{2}+C_{\alpha }\left\Vert z(t)\right\Vert _{H}^{2},
\end{equation*}%
\begin{eqnarray*}
\alpha \gamma (Ay(t),Az(t))_{H} &\leq &8\gamma ^{2}\alpha \left\Vert
Ay(t)\right\Vert _{H}^{2}+\frac{\alpha }{8}\left\Vert Az(t)\right\Vert
_{H}^{2} \\
&\leq &\delta \left\Vert A^{3/2}y(t)\right\Vert _{H}^{2}+C_{\delta }\alpha
^{2}\left\Vert y(t)\right\Vert _{H}^{2}+\frac{\alpha }{8}\left\Vert
Az(t)\right\Vert _{H}^{2},
\end{eqnarray*}%
\begin{equation*}
-\alpha \gamma (y(t),Az(t))_{H}\leq \frac{\alpha }{8}\left\Vert
Az(t)\right\Vert _{H}^{2}+C_{\alpha }\left\Vert y(t)\right\Vert _{H}^{2},
\end{equation*}%
\begin{equation*}
\int_{\Omega }1_{\omega }^{\ast }v(t)A^{2}y(t)dx\leq \delta \left\Vert
A^{2}y(t)\right\Vert _{H}^{2}+C_{\delta }\left\Vert v(t)\right\Vert _{H}^{2},
\end{equation*}%
\begin{equation*}
\int_{\Omega }1_{\omega }^{\ast }u(t)Az(t)dx\leq \frac{\alpha }{8}\left\Vert
Az(t)\right\Vert _{H}^{2}+C\left\Vert u(t)\right\Vert _{H}^{2}.
\end{equation*}%
Plugging all these relations in (\ref{63+1}), using $\left\Vert
A^{3/2}y(t)\right\Vert _{H}^{2}\leq C\left\Vert A^{2}y(t)\right\Vert
_{H}^{2} $ and choosing $\delta $ small enough we get
\begin{eqnarray*}
&&\frac{1}{2}\frac{d}{dt}\left( \left\Vert Ay(t)\right\Vert _{H}^{2}+\alpha
\left\Vert A^{1/2}z(t)\right\Vert _{H}^{2}\right) +\frac{\nu }{2}\left\Vert
A^{2}y(t)\right\Vert _{H}^{2}+\left( \frac{\alpha }{2}-C_{2}\right)
\left\Vert Az(t)\right\Vert _{H}^{2} \\
&\leq &C_{\alpha }\left( \left\Vert y(t)\right\Vert _{H}^{2}+\left\Vert
z(t)\right\Vert _{H}^{2}+\left\Vert v(t)\right\Vert _{H}^{2}+\left\Vert
u(t)\right\Vert _{H}^{2}\right) ,
\end{eqnarray*}%
where $C_{2}$ can be computed and depends only the system parameters, and $%
C_{\alpha }$ depends on $\alpha ,$ in addition. We choose, for instance, $%
\alpha =4C_{2},$ integrate from 0 to $t,$ use (\ref{34})-(\ref{34-0}) in
order to find
\begin{eqnarray*}
&&\left\Vert Ay(t)\right\Vert _{H}^{2}+\left\Vert A^{1/2}z(t)\right\Vert
_{H}^{2}+\int_{0}^{t}\left( \left\Vert A^{2}y(s)\right\Vert
_{H}^{2}+\left\Vert Az(s)\right\Vert _{H}^{2}\right) ds \\
&\leq &C\left( \left\Vert Ay^{0}\right\Vert _{H}^{2}+\left\Vert
A^{1/2}z^{0}\right\Vert _{H}^{2}\right) +CC_{p}^{2}(\left\Vert
y^{0}\right\Vert _{H}^{2}+\left\Vert z^{0}\right\Vert _{H}^{2}).
\end{eqnarray*}%
Writing this relation for the optimal pair we obtain (\ref{52+1}), as
claimed. \hfill $\square $

\medskip

Let us point out a first consequence of Proposition 2.3. It is not difficult
to check that the functional
\begin{equation*}
\left\Vert \cdot \right\Vert _{\Phi }=\sqrt{\Phi }:D(A^{1/2})\times
D(A^{1/4})\rightarrow \mathbb{R}
\end{equation*}%
is a norm satisfying the parallelogram law. Then, (\ref{52+0}) implies that $%
\left\Vert \cdot \right\Vert _{\Phi }$ is a Hilbert norm on $%
D(A^{1/2})\times D(A^{1/4})$ equivalent to the natural one. In addition, $%
\Phi $ is a quadratic functional. Moreover, if we denote $\left( \cdot
,\cdot \right) _{\Phi }$ the corresponding scalar product we can introduce%
\begin{equation}
R:\Xi :=D(A^{1/2})\times D(A^{1/4})\rightarrow \Xi ^{\prime
}=(D(A^{1/2})\times (D(A^{1/4}))^{\prime }  \label{64-0}
\end{equation}%
such that
\begin{equation*}
\left\langle R(y^{0},z^{0}),(Y,Z)\right\rangle _{\Xi ^{\prime },\Xi }=\frac{1%
}{2}\left( (y^{0},z^{0}),(Y,Z)\right) _{\Phi }\text{ for all }%
(y^{0},z^{0}),(Y,Z)\in \Xi .
\end{equation*}%
In fact $R$ coincides with $2R_{\Xi },$ where $R_{\Xi }$ is the Riesz
operator associated to $\left\Vert \cdot \right\Vert _{\Phi }.$ In
particular, we have
\begin{equation}
\Phi (y^{0},z^{0})=\frac{1}{2}\left( R(y^{0},z^{0}),(y^{0},z^{0})\right)
\text{ for all }(y^{0},z^{0})\in D(A^{1/2})\times D(A^{1/4}).  \label{64}
\end{equation}%
Moreover, $R(y^{0},z^{0})$ is the G\^{a}teaux derivative of the function $%
\Phi $ at $(y^{0},z^{0}).$ Indeed, for any $(Y,Z)\in D(A^{1/2})\times
D(A^{1/4})=\Xi $ we have%
\begin{eqnarray*}
&&\Phi ^{\prime }(y^{0},z^{0})(Y,Z)=\lim_{\lambda \rightarrow 0}\frac{\Phi
(y^{0}+\lambda Y,z^{0}+\lambda Z)-\Phi (y^{0},z^{0})}{\lambda } \\
&=&\frac{1}{2}\lim_{\lambda \rightarrow 0}\frac{\left\langle R(y^{0}+\lambda
Y,z^{0}+\lambda Z),(y^{0}+\lambda Y,z^{0}+\lambda Z)\right\rangle _{\Xi
^{\prime },\Xi }-\left\langle R(y^{0},z^{0}),(y^{0},z^{0})\right\rangle
_{\Xi ^{\prime },\Xi }}{\lambda } \\
&=&\left\langle R(y^{0},z^{0}),(Y,Z)\right\rangle _{\Xi ^{\prime },\Xi },
\end{eqnarray*}%
hence
\begin{equation}
\Phi ^{\prime }(y^{0},z^{0})=R(y^{0},z^{0}),\text{ for all }(y^{0},z^{0})\in
D(A^{1/2})\times D(A^{1/4}).  \label{65}
\end{equation}

Since $\Phi $ is coercive by (\ref{52+0}) we can define the restriction of $%
R $ to $H\times H$ (denoted still by $R$) having the domain
\begin{equation*}
D(R)=\{(y^{0},z^{0})\in \Xi ;\text{ }R(y^{0},z^{0})\in H\times H\}.
\end{equation*}%
It also turns out that $R$ is self-adjoint. Moreover, $R$ can be written of
the form
\begin{equation}
R=\left[
\begin{array}{cc}
R_{11} & R_{12} \\
R_{21} & R_{22}%
\end{array}%
\right] .  \label{R}
\end{equation}%
We shall give more details about this in the next Proposition which also
provides a representation result for the optimal solution to (\ref{51}).

Let us recall the operators $B:\mathbb{R}^{N}\rightarrow H\times H,$ $%
B^{\ast }:H\times H\rightarrow \mathbb{R}^{N},$ defined in (\ref{B}) and (%
\ref{B*}),
\begin{equation*}
BW=\left[
\begin{tabular}{c}
$\sum\limits_{i=1}^{N}1_{\omega }^{\ast }\varphi _{i}w_{i}$ \\
$\sum\limits_{i=1}^{N}1_{\omega }^{\ast }\psi _{i}w_{i}$%
\end{tabular}%
\right] \text{ for all }W=\left[
\begin{tabular}{c}
$w_{1}$ \\
$...$ \\
$w_{N}$%
\end{tabular}%
\right] \in \mathbb{R}^{N},
\end{equation*}%
and
\begin{equation*}
B^{\ast }q=\left[
\begin{tabular}{c}
$\int_{\Omega }1_{\omega }^{\ast }(\varphi _{1}q_{1}+\psi _{1}q_{2})dx$ \\
$...$ \\
$\int_{\Omega }1_{\omega }^{\ast }(\varphi _{N}q_{1}+\psi _{N}q_{2})dx$%
\end{tabular}%
\right] \text{, for all }q=\left[
\begin{array}{c}
q_{1} \\
q_{2}%
\end{array}%
\right] \in H\times H.
\end{equation*}%
Then, (\ref{31})-(\ref{32}) can be rewritten as
\begin{eqnarray}
\frac{d}{dt}(y(t),z(t))+\mathcal{A}(y(t),z(t)) &=&BW(t),\text{ a.e. }t>0,
\label{31'} \\
(y(0),z(0)) &=&(y^{0},z^{0}).  \notag
\end{eqnarray}

\medskip

\noindent \textbf{Proposition 2.4. }\textit{Let }$W^{\ast }=\{w_{i}^{\ast
}\}_{i=1}^{N}$ \textit{and} $(y^{\ast },z^{\ast })$ \textit{be optimal for
problem} (\ref{51}), \textit{corresponding to} $(y^{0},z^{0})\in
D(A^{1/2})\times D(A^{1/4}).$ \textit{Then, }$W^{\ast }$ \textit{is
expressed as }%
\begin{equation}
W^{\ast }(t)=-B^{\ast }R(y^{\ast }(t),z^{\ast }(t)),\text{ \textit{for all} }%
t>0.  \label{66}
\end{equation}%
\textit{\ Moreover, }$R$ \textit{has the following properties}%
\begin{eqnarray}
2c_{1}\left\Vert (y^{0},z^{0})\right\Vert _{D(A^{1/2})\times D(A^{1/4})}^{2}
&\leq &(R(y^{0},z^{0}),(y^{0},z^{0}))_{H\times H}\leq 2c_{2}\left\Vert
(y^{0},z^{0})\right\Vert _{D(A^{1/2})\times D(A^{1/4})}^{2},  \notag \\
\text{ \textit{for all} }(y^{0},z^{0}) &\in &D(A^{1/2})\times D(A^{1/4}),
\label{69}
\end{eqnarray}%
\begin{equation}
\left\Vert R(y^{0},z^{0})\right\Vert _{H\times H}\leq C_{R}\left\Vert
(y^{0},z^{0})\right\Vert _{D(A)\times D(A^{1/2})},\text{ \textit{for}
\textit{all} }(y^{0},z^{0})\in D(A)\times D(A^{1/2}),  \label{68}
\end{equation}%
\textit{\ and satisfies the Riccati algebraic equation }(\ref{18'})\textit{,
that is}%
\begin{eqnarray}
&&2\left( R(\overline{y},\overline{z}),\mathcal{A}(\overline{y},\overline{z}%
)\right) _{H\times H}+\left\Vert B^{\ast }R(\overline{y},\overline{z}%
)\right\Vert _{\mathbb{R}^{N}}^{2}  \label{67} \\
&=&\left\Vert A^{3/2}\overline{y}\right\Vert _{H}^{2}+\left\Vert A^{3/4}%
\overline{z}\right\Vert _{H}^{2},\text{ \textit{for} \textit{all} }(%
\overline{y},\overline{z})\in D(A^{3/2})\times D(A^{3/4}).  \notag
\end{eqnarray}%
\textit{Here,} $c_{1},$ $c_{2}$\textit{, }$C_{R}$ \textit{are constants} ($%
c_{1}$, $c_{2}$ \textit{are the same as in} (\ref{52+0}) \textit{and depend
on the problem parameters, }$\Omega $\textit{\ and} $\left\Vert \varphi
_{\infty }\right\Vert _{L^{2}(\Omega )},$ \textit{and }$C_{R}$\textit{\
depends only on} $\Omega ).$

\medskip

\noindent \textbf{Proof. }We organize the proof in two steps.

\subparagraph{Step 1.}

Inequalities (\ref{69}) immediately follow from (\ref{64}) and (\ref{52+0}).

Next we prove (\ref{68}) and (\ref{66}).

Let $T$ be positive and arbitrary. We recall that by the dynamic programming
principle (see e.g., \cite{vb-optim-94}, p. 104), the minimization problem (%
\ref{51}) is equivalent to the following problem
\begin{equation}
\underset{W\in L^{2}(0,T;\mathbb{R}^{N})}{\text{Min}}\left\{ \frac{1}{2}%
\int_{0}^{T}\left( \left\Vert A^{3/2}y(t)\right\Vert _{H}^{2}+\left\Vert
A^{3/4}z(t)\right\Vert _{H}^{2}+\left\Vert W(t)\right\Vert _{\mathbb{R}%
^{N}}^{2}\right) dt+\Phi (y(T),z(T))\right\}  \label{70}
\end{equation}%
subject to (\ref{31})-(\ref{32}). Thus, a solution to (\ref{51}) is a
solution to (\ref{70}) on $(0,T)$ and conversely.

We introduce the adjoint system
\begin{eqnarray}
\frac{d}{dt}(p^{T},q^{T})(t)-\mathcal{A}(p^{T}(t),q^{T}(t)) &=&(A^{3}y^{\ast
}(t),A^{3/2}z^{\ast }(t)),\text{ in }(0,T)\times \Omega ,  \label{71} \\
(p^{T}(T),q^{T}(T)) &=&-R(y^{\ast }(T),z^{\ast }(T)),\text{ in }\Omega ,
\notag
\end{eqnarray}%
by recalling that $\mathcal{A}$ is self-adjoint. We have used (\ref{65}) for
writing the final condition at $t=T.$ For the moment we indicate this
solution by $(p^{T},q^{T}),$ to show its dependence on $T.$ Later, we shall
prove that actually it is independent of $T.$ By the maximum principle in (%
\ref{70}), we have that
\begin{equation}
W^{\ast }(t)=B^{\ast }(p^{T}(t),q^{T}(t)),\text{ a.e. }t\in (0,T)  \label{72}
\end{equation}%
\textbf{(}see \cite{lions-control}, p. 114; see also \cite{vb-optim-94}, p.
190\textbf{).}

For proving (\ref{68}), let $(y^{0},z^{0})\in D(A)\times D(A^{1/2}).$

Since
\begin{equation*}
R(y^{\ast }(T),z^{\ast }(T))\in (D(A^{1/2}))^{\prime }\times
(D(A^{1/4}))^{\prime }\subset V^{\prime }\times V^{\prime }
\end{equation*}%
and
\begin{equation*}
(A^{3}y^{\ast },A^{3/2}z^{\ast })\in L^{2}(0,T;V^{\prime }\times V^{\prime
}),
\end{equation*}%
it follows that (\ref{71}) has a unique solution
\begin{equation}
(p^{T},q^{T})\in L^{2}(0,T;H\times H)\cap C([0,T];V^{\prime }\times
V^{\prime })  \label{72-1}
\end{equation}%
(see \cite{baiocchi}, Thm. 7.1, p. 291). We shall prove that $(p^{T},q^{T})$
is in $C([0,T);H\times H)$. For the reader's convenience we give the
argument, adapting some ideas from the proof in \cite{vb-gw-2002}. We define
\begin{equation}
(\widetilde{p},\widetilde{q})=\widetilde{A}(p^{T},q^{T})\text{ }  \label{74}
\end{equation}%
where $\widetilde{A}$ is the operator
\begin{equation*}
\widetilde{A}=\left[
\begin{array}{cc}
A^{-1} & 0 \\
0 & A^{-1/2}%
\end{array}%
\right] .
\end{equation*}%
By recalling (\ref{53}) we see that $\mathcal{A}$ and $\widetilde{A}$
commute. Thus, we replace (\ref{74}) in (\ref{71}), obtaining the system
\begin{eqnarray}
\frac{d}{dt}(\widetilde{p},\widetilde{q})(t)-\mathcal{A}(\widetilde{p}(t),%
\widetilde{q}(t)) &=&(A^{2}y^{\ast }(t),Az^{\ast }(t)),\text{ in }%
(0,T)\times \Omega ,  \label{75} \\
(\widetilde{p}(T),\widetilde{q}(T)) &=&-\widetilde{A}R(y^{\ast }(T),z^{\ast
}(T))\text{, in }\Omega .  \notag
\end{eqnarray}%
According to (\ref{52+1}), we have $(A^{2}y^{\ast },Az^{\ast })\in
L^{2}(0,T;H\times H)$ and by (\ref{72-1}) we obtain $\widetilde{A}R(y^{\ast
}(T),z^{\ast }(T))\in V\times H$. By applying a backward version of
Proposition 2.1, formula (\ref{32-3}) we see that system (\ref{75}) has a
unique solution
\begin{equation*}
(\widetilde{p},\widetilde{q})\in C([0,T);D(A)\times D(A^{1/2}))\text{ }
\end{equation*}%
and so $(p^{T},q^{T})\in C([0,T);H\times H).$ Next, we prove the relation
\begin{equation}
R(y^{0},z^{0})=-(p^{T}(0),q^{T}(0)).  \label{76}
\end{equation}%
To this end, let us consider two solutions to (\ref{70}), $(W^{\ast
},y^{\ast },z^{\ast })$ and $(W_{1}^{\ast },y_{1}^{\ast },z_{1}^{\ast })$,
corresponding to $(y^{0},z^{0})$ and $(y^{1},z^{1}),$ respectively, both in $%
D(A)\times D(A^{1/2}).$ Using the subdifferential inequality $\left\Vert
v\right\Vert _{X}-\left\Vert v_{1}\right\Vert _{X}\leq 2(v,v-v_{1})_{X},$
which holds in any Hilbert space $X,$ and the relation $\Phi ^{\prime }=R$
we compute
\begin{eqnarray}
&&\Phi (y^{0},z^{0})-\Phi (y_{1},z_{1})  \label{77} \\
&\leq &\int_{0}^{T}\left\{ (A^{3/2}y^{\ast }(t),A^{3/2}(y^{\ast
}(t)-y_{1}^{\ast }(t)))_{H}+(A^{3/4}z^{\ast }(t),A^{3/4}(z^{\ast
}(t)-z_{1}^{\ast }(t)))_{H}\right\} dt  \notag \\
&&+\int_{0}^{T}(W^{\ast }(t),W^{\ast }(t)-W_{1}^{\ast }(t))_{\mathbb{R}%
^{N}}dt  \notag \\
&&+(R(y^{\ast }(T),z^{\ast }(T)),(y^{\ast }(T)-y_{1}^{\ast }(T),z^{\ast
}(T)-z_{1}^{\ast }(T)))_{H\times H}.  \notag
\end{eqnarray}%
It is clear that $W^{\ast }=(w_{1}^{\ast },...,w_{N}^{\ast })$ and $%
W_{1}^{\ast }=(w_{11}^{\ast },...,w_{1N}^{\ast }),$ respectively. By
multiplying (\ref{71}) by $(y^{\ast }(t)-y_{1}^{\ast }(t),z^{\ast
}(t)-z_{1}^{\ast }(t)),$ integrating by parts and using the difference of
the state equations (\ref{31'}), written for both solutions, we obtain that
\begin{eqnarray}
&&\frac{d}{dt}((p^{T}(t),q^{T}(t)),(y^{\ast }(t)-y_{1}^{\ast }(t),z^{\ast
}(t)-z_{1}^{\ast }(t)))_{H\times H}  \label{78} \\
&=&(A^{3/2}y^{\ast }(t),A^{3/2}(y^{\ast }(t)-y_{1}^{\ast
}(t)))_{H}+(A^{3/4}z^{\ast }(t),A^{3/4}(z^{\ast }(t)-z_{1}^{\ast }(t)))_{H}
\notag \\
&&+((p^{T}(t),q^{T}(t)),BW^{\ast }(t)-BW_{1}^{\ast }(t))_{H\times H},\text{
a.e. }t>0.  \notag
\end{eqnarray}%
Now, we integrate (\ref{78}) over $(0,T)$ and use the final condition in (%
\ref{71}) and (\ref{72}), to get
\begin{eqnarray*}
&&-(R(y^{\ast }(T),z^{\ast }(T)),(y^{\ast }(T)-y_{1}^{\ast }(T),z^{\ast
}(T)-z_{1}^{\ast }(T)))_{H\times H} \\
&&-((p^{T}(0),q^{T}(0)),(y^{\ast }(0)-y_{1}^{\ast }(0),z^{\ast
}(0)-z_{1}^{\ast }(0)))_{H\times H} \\
&=&\int_{0}^{T}\left\{ (A^{3/2}y^{\ast }(t),A^{3/2}(y^{\ast }(t)-y_{1}^{\ast
}(t)))_{H}+(A^{3/4}z^{\ast }(t),A^{3/4}(z^{\ast }(t)-z_{1}^{\ast
}(t)))_{H}\right\} dt \\
&&+\int_{0}^{T}(W^{\ast }(t),W^{\ast }(t)-W_{1}^{\ast }(t))_{\mathbb{R}%
^{N}}dt,
\end{eqnarray*}%
whence by (\ref{77}) we finally obtain
\begin{equation}
\Phi (y^{0},z^{0})-\Phi (y^{1},z^{1})\leq
-((p^{T}(0),q^{T}(0)),(y^{0}-y^{1},z^{0}-z^{1}))_{H\times H}.  \label{79}
\end{equation}%
This implies that
\begin{equation*}
-(p^{T}(0),q^{T}(0))\in \partial \Phi (y^{0},z^{0}).
\end{equation*}%
Since, as seen earlier, $\Phi $ is differentiable on $D(A^{1/2})\times
D(A^{1/4})$ it follows that
\begin{equation*}
-(p^{T}(0),q^{T}(0))=\Phi ^{\prime }(y^{0},z^{0})=R(y^{0},z^{0}),
\end{equation*}%
as claimed in (\ref{76}). This implies, since we have proved that $%
(p^{T},q^{T})\in C([0,\infty );H\times H),$ that $(p^{T}(0),q^{T}(0))\in
H\times H$ and so
\begin{equation}
R(y^{0},z^{0})\in H\times H\text{ for all }(y^{0},z^{0})\in D(A)\times
D(A^{1/2}).  \label{79-1}
\end{equation}

On the other hand, one can easily see that $R$ is a linear closed operator
from $D(A)\times D(A^{1/2})$ to $H\times H$, and so by the closed graph
theorem we conclude that it is continuous (see e.g, \cite{brezis}, Thm. 2.9,
p. 37), that is $R\in \mathcal{L}(D(A)\times D(A^{1/2});H\times H)$, as
claimed by (\ref{68}).

We define the restriction of $R$ to $H\times H,$ still denoted by $R.$ Thus,
its domain contains $D(A)\times D(A^{1/2})$.

Now, we resume (\ref{72}) which extends by the continuity (\ref{72-1}) at $%
t=T,$ in $V^{\prime }.$%
\begin{equation}
W^{\ast }(T)=B^{\ast }(p^{T}(T),q^{T}(T)).  \label{72'}
\end{equation}%
Moreover, since $(y^{\ast }(t),z^{\ast }(t))\in D(A)\times D(A^{1/2})$ for
all $t\geq 0,$ by (\ref{52+1}), we have by (\ref{79-1}) that $R(y^{\ast
}(t),z^{\ast }(t))\in H\times H$ for all $t\geq 0.$ In particular, this is
true for $t=T$ and so using the final condition in (\ref{71}) we get
\begin{equation}
(p^{T}(T),q^{T}(T))=-R(y^{\ast }(T),z^{\ast }(T))\in H\times H.  \label{71'}
\end{equation}%
This relation combined with (\ref{72'}) implies
\begin{equation*}
W^{\ast }(T)=-B^{\ast }R(y^{\ast }(T),z^{\ast }(T))
\end{equation*}%
where $T$ is arbitrary. Therefore, it can be written for any $t,$ as in (\ref%
{66}), as claimed.

By (\ref{66}) and by the definition (\ref{B*}) and (\ref{32+10}) we can write%
\begin{equation}
w_{j}=-(B^{\ast }R(y^{\ast }(t),z^{\ast }(t)))_{j}=-\int_{\Omega }1_{\omega
}^{\ast }(\varphi _{j}R_{1}(y^{\ast }(t),z^{\ast }(t))+\psi
_{j}R_{2}(y^{\ast }(t),z^{\ast }(t)))dx  \label{wj-fin}
\end{equation}%
and by (\ref{R}) we get
\begin{equation}
R_{1}(y^{\ast }(t),z^{\ast }(t))=R_{11}y^{\ast }(t)+R_{12}z^{\ast }(t),\text{
\ }R_{2}(y^{\ast }(t),z^{\ast }(t))=R_{21}y^{\ast }(t)+R_{22}z^{\ast }(t).
\label{R1-R2}
\end{equation}%
In particular,
\begin{eqnarray*}
(v^{\ast },u^{\ast }) &=&\left( \sum\limits_{j=1}^{N}\varphi
_{j}w_{j},\sum\limits_{j=1}^{N}\psi _{j}w_{j}\right) \\
&=&\left( -\sum_{j=1}^{N}\varphi _{j}(B^{\ast }R(y^{\ast }(t),z^{\ast
}(t)))_{j},-\sum_{j=1}^{N}\psi _{j}(B^{\ast }R(y^{\ast }(t),z^{\ast
}(t)))_{j}\right)
\end{eqnarray*}%
which implies by (\ref{B*}) the representation%
\begin{eqnarray}
v^{\ast }(t,x) &=&-\sum_{j=1}^{N}\varphi _{j}(x)(B^{\ast }R(y^{\ast
}(t),z^{\ast }(t)))_{j}  \label{v} \\
&=&-\sum_{j=1}^{N}\varphi _{j}(x)\int_{\Omega }1_{\omega }^{\ast }(\varphi
_{j}R_{1}(y^{\ast }(t),z^{\ast }(t))+\psi _{j}R_{2}(y^{\ast }(t),z^{\ast
}(t)))(\xi )d\xi ,  \notag
\end{eqnarray}%
\begin{eqnarray}
u^{\ast }(t,x) &=&-\sum_{j=1}^{N}\psi _{j}(x)(B^{\ast }R(y^{\ast
}(t),z^{\ast }(t)))_{j}  \label{u} \\
&=&-\sum_{j=1}^{N}\psi _{j}(x)\int_{\Omega }1_{\omega }^{\ast }(\varphi
_{j}R_{1}(y^{\ast }(t),z^{\ast }(t))+\psi _{j}R_{2}(y^{\ast }(t),z^{\ast
}(t)))(\xi )d\xi .  \notag
\end{eqnarray}%
Finally, it follows by (\ref{66})\ that%
\begin{equation}
1_{\omega }^{\ast }U(t)=1_{\omega }^{\ast }(v^{\ast }(t),u^{\ast
}(t))=-BB^{\ast }R(y^{\ast }(t),z^{\ast }(t)).  \label{control-fin}
\end{equation}

\subparagraph{Step 2.}

We pass now to the proof of (\ref{67}). It is enough to consider $%
(y^{0},z^{0})\in D(A)\times D(A^{1/2}).$ Since $(W^{\ast },y^{\ast },z^{\ast
})$ is the solution to both (\ref{51}) and (\ref{70}) written with $T=t$
where $t\geq 0$ is arbitrary and the minimum is $\Phi (y^{0},z^{0})$ we can
write
\begin{eqnarray*}
&&\Phi (y^{0},z^{0})=\frac{1}{2}\int_{0}^{t}\left( \left\Vert A^{3/2}y^{\ast
}(s)\right\Vert _{H}^{2}+\left\Vert A^{3/4}z^{\ast }(s)\right\Vert
_{H}^{2}+\left\Vert W^{\ast }(s)\right\Vert _{\mathbb{R}^{N}}^{2}\right) ds
\\
&&+\frac{1}{2}\int_{t}^{\infty }\left( \left\Vert A^{3/2}y^{\ast
}(s)\right\Vert _{H}^{2}+\left\Vert A^{3/4}z^{\ast }(s)\right\Vert
_{H}^{2}+\left\Vert W^{\ast }(s)\right\Vert _{\mathbb{R}^{N}}^{2}\right) ds
\\
&=&\Phi (y^{0},z^{0})-\Phi (y^{\ast }(t),z^{\ast }(t))+\frac{1}{2}%
\int_{t}^{\infty }\left( \left\Vert A^{3/2}y^{\ast }(s)\right\Vert
_{H}^{2}+\left\Vert Az^{\ast }(s)\right\Vert _{H}^{2}+\left\Vert W^{\ast
}(s)\right\Vert _{\mathbb{R}^{N}}^{2}\right) ds
\end{eqnarray*}%
and so
\begin{equation}
\Phi (y^{\ast }(t),z^{\ast }(t))=\frac{1}{2}\int_{t}^{\infty }\left(
\left\Vert A^{3/2}y^{\ast }(s)\right\Vert _{H}^{2}+\left\Vert A^{3/4}z^{\ast
}(s)\right\Vert _{H}^{2}+\left\Vert W^{\ast }(s)\right\Vert _{\mathbb{R}%
^{N}}^{2}\right) ds,\text{ }  \label{81}
\end{equation}%
for any $t\geq 0.$ Now, we want to differentiate (\ref{81}) with respect to $%
t.$ To this aim, we recall (\ref{64}) and that $R$ is symmetric. Thus,
\begin{eqnarray*}
\frac{d}{dt}\Phi (y^{\ast }(t),z^{\ast }(t)) &=&\frac{1}{2}\frac{d}{dt}%
(R(y^{\ast }(t),z^{\ast }(t)),(y^{\ast }(t),z^{\ast }(t)))_{H\times H} \\
&=&\left( R(y^{\ast }(t),z^{\ast }(t)),\frac{d}{dt}(y^{\ast }(t),z^{\ast
}(t))\right) _{H\times H}.
\end{eqnarray*}%
Hence, taking into account (\ref{66}) we obtain for a.e. $t>0$ (since $%
A^{3/2}y^{\ast }(t)$ is defined only for a.e. $t)$ that
\begin{eqnarray}
&&\left( R(y^{\ast }(t),z^{\ast }(t)),\frac{d}{dt}(y^{\ast }(t),z^{\ast
}(t))\right) _{H\times H}+\frac{1}{2}\left( \left\Vert A^{3/2}y^{\ast
}(t)\right\Vert _{H}^{2}+\left\Vert A^{3/4}z^{\ast }(t)\right\Vert
_{H}^{2}\right)  \label{83} \\
&&+\frac{1}{2}\left\Vert B^{\ast }R(y^{\ast }(t),z^{\ast }(t))\right\Vert _{%
\mathbb{R}^{N}}^{2}=0.  \notag
\end{eqnarray}%
Now, we come back to the system (\ref{31'}) in which the right-hand side is
replaced by (\ref{control-fin}). This becomes a \textit{closed loop} system
with the right-hand side $-BB^{\ast }R(y^{\ast }(t),z^{\ast }(t)).$

We take into account that by (\ref{68})
\begin{equation*}
\left\Vert BB^{\ast }R(y^{\ast }(t),z^{\ast }(t))\right\Vert _{H\times
H}\leq C_{1}\left\Vert R(y^{\ast }(t),z^{\ast }(t))\right\Vert _{H\times
H}\leq C_{2}\left\Vert (y^{\ast }(t),z^{\ast }(t))\right\Vert _{D(A)\times
D(A^{1/2})}
\end{equation*}%
a.e. $t>0,$ for all $(y^{\ast }(t),z^{\ast }(t))\in D(A)\times D(A^{1/2}).$

We show that $-(\mathcal{A}+BB^{\ast }R)$ generates a $C_{0}$-semigroup in $%
H\times H$, using Lemma A3 in Appendix.

In our case, we particularize $E=D(A^{2})\times D(A),$ $F=H\times H,$ $L=%
\mathcal{A}$ and $M=BB^{\ast }R.$ Then, the operator $\mathcal{A}+BB^{\ast
}R $ is quasi $m$-accretive and we have the result. Thus, for $%
(y^{0},z^{0})\in H\times H$ we have%
\begin{equation*}
\mathcal{A}(y^{\ast }(t),z^{\ast }(t)),\text{ }BB^{\ast }R(y^{\ast
}(t),z^{\ast }(t)),\text{ }\frac{d}{dt}(y^{\ast }(t),z^{\ast }(t))\in
C((0,\infty );H\times H)
\end{equation*}%
(see \cite{brezis-73}, p. 72, for a basic result). Then, we can replace $%
\frac{d}{dt}(y^{\ast }(t),z^{\ast }(t))$ from (\ref{31'}) and plug it in (%
\ref{83}). On account of (\ref{66}) we have
\begin{eqnarray*}
&&\left( R(y^{\ast }(t),z^{\ast }(t)),-\mathcal{A}(y^{\ast }(t),z^{\ast
}(t))\right) _{H\times H}+\frac{1}{2}\left( \left\Vert A^{3/2}y^{\ast
}(t)\right\Vert _{H}^{2}+\left\Vert A^{3/4}z^{\ast }(t)\right\Vert
_{H}^{2}\right) \\
&&+\frac{1}{2}\left\Vert B^{\ast }R(y^{\ast }(t),z^{\ast }(t))\right\Vert _{%
\mathbb{R}^{N}}^{2}=(R(y^{\ast }(t),z^{\ast }(t)),BB^{\ast }R(y^{\ast
}(t),z^{\ast }(t))_{H\times H},\text{ }t\geq 0
\end{eqnarray*}%
which implies (\ref{67}) (written with a generic notation $(\overline{y},%
\overline{z})\in D(A^{3/2})\times D(A^{3/4}))$, as claimed.

\hfill $\square $

\medskip

\noindent \textbf{Remark 2.5. }We note that the previous equation can be
still written%
\begin{eqnarray*}
&&2\left( R(y^{\ast }(t),z^{\ast }(t)),\mathcal{A}(y^{\ast }(t),z^{\ast
}(t))\right) _{H\times H}+(B^{\ast }R(y^{\ast }(t),z^{\ast }(t)),B^{\ast
}R(y^{\ast }(t),z^{\ast }(t)))_{\mathbb{R}^{N}\times \mathbb{R}^{N}} \\
&=&(\widehat{\mathcal{A}}(y^{\ast }(t),z^{\ast }(t)),(y^{\ast }(t),z^{\ast
}(t)))_{H\times H},
\end{eqnarray*}%
where $\widehat{\mathcal{A}}:D(A^{3/2})\times D(A^{3/4})\rightarrow
(D(A^{3/2})\times D(A^{3/4}))^{\prime }$ is defined by%
\begin{align*}
&\quad\ \left\langle \widehat{\mathcal{A}}(y,z),(\psi _{1},\psi _{2})\right\rangle
_{(D(A^{3/2})\times D(A^{3/4}))^{\prime },D(A^{3/2})\times
D(A^{3/4})}\\
&=(A^{3/2}y,A^{3/2}\psi _{1})_{H\times H}+(A^{3/4}z,A^{3/4}\psi
_{2})_{H\times H}
\end{align*}%
for all $(\psi _{1},\psi _{2})\in D(A^{3/2})\times D(A^{3/4}).$ By $\widehat{%
\mathcal{A}}$ we still denote its restriction to $H\times H,$
\begin{equation*}
\widehat{\mathcal{A}}=\left[
\begin{tabular}{ll}
$A^{3}$ & $0$ \\
$0$ & $A^{3/2}$%
\end{tabular}%
\right]
\end{equation*}%
with the domain $D(\widehat{\mathcal{A}})=\{(y,z)\in D(A^{3/2})\times
D(A^{3/4});$ $\widehat{\mathcal{A}}(y,z)\in H\times H\}.$

Since%
\begin{eqnarray*}
((R(y,z),\mathcal{A}(y,z))_{H\times H} &=&(\mathcal{A}R(y,z),(y,z))_{H\times
H}, \\
((\mathcal{A}(y,z),R(y,z))_{H\times H} &=&((y,z),\mathcal{A}R(y,z))_{H\times
H}
\end{eqnarray*}%
we get $\mathcal{A}R=R\mathcal{A},$ thus
\begin{equation*}
2R\mathcal{A}(y^{\ast }(t),z^{\ast }(t))+RBB^{\ast }R(y^{\ast }(t),z^{\ast
}(t))=\widehat{\mathcal{A}}(y^{\ast }(t),z^{\ast }(t)),
\end{equation*}%
for all $t>0.$ Letting $t\rightarrow 0$ we obtain
\begin{equation*}
2R\mathcal{A}(y^{0},z^{0})+RBB^{\ast }R(y^{0},z^{0})=\widehat{\mathcal{A}}%
(y^{0},z^{0}),
\end{equation*}%
for all $(y^{0},z^{0})\in D(A)\times D(A^{1/2})$ and so the Riccati equation
takes the form%
\begin{equation}
2R\mathcal{A+}RBB^{\ast }R=\widehat{\mathcal{A}}.  \label{Riccati}
\end{equation}

\medskip

\noindent \textbf{Remark 2.6. }Just as a remark, we observe that the linear
system is exponentially stabilized to $(0,0)$ by the feedback controller
just constructed. To sustain this assertion we recall a generalization of
Datko's result (see Lemma A4 in Appendix, see also \cite{pazy}, p.116). In
our case, the operator $\mathcal{D}=-(\mathcal{A}+BB^{\ast }R)$ generates a $%
C_{0}$-semigroup in $H\times H$ and, as seen earlier, equation (\ref{31})
\begin{equation*}
\frac{d}{dt}(y(t),z(t))+(\mathcal{A+}BB^{\ast }R)(y(t),z(t))=0,\text{ }t\geq
0
\end{equation*}%
has the solution $(y(t),z(t))$ with the property%
\begin{eqnarray*}
\int_{0}^{\infty }\left( \left\Vert y(t)\right\Vert _{H}^{2}+\left\Vert
z(t)\right\Vert _{H}^{2}\right) dt &\leq &C\int_{0}^{\infty }\left(
\left\Vert A^{3/2}y(t)\right\Vert _{H}^{2}+\left\Vert A^{3/4}z(t)\right\Vert
_{H}^{2}\right) dt \\
&\leq &C\left( \left\Vert y^{0}\right\Vert _{D(A^{1/2})}^{2}+\left\Vert
z^{0}\right\Vert _{D(A^{1/4})}^{2}\right) <\infty .
\end{eqnarray*}%
Hence
\begin{equation*}
\left\Vert y(t)\right\Vert _{H}^{2}+\left\Vert z(t)\right\Vert _{H}^{2}\leq
Ce^{-\kappa t}\left( \left\Vert y^{0}\right\Vert
_{D(A^{1/2})}^{2}+\left\Vert z^{0}\right\Vert _{D(A^{1/4})}^{2}\right) ,%
\text{ for all }t\geq 0,
\end{equation*}%
as claimed.

\section{Feedback stabilization of the nonlinear system}

\setcounter{equation}{0}

We recall that $B$ and $B^{\ast }$ are defined by (\ref{B}) and (\ref{B*}),
and that $R_{1}=(R_{11},R_{12}),$ $R_{2}=(R_{21},R_{22})$ are given by (\ref%
{R1-R2}).

In this section we shall deal with the nonlinear system (\ref{19''})-(\ref%
{22''}) in which the right-hand side $(1_{\omega }^{\ast }v,1_{\omega
}^{\ast }u)$ is replaced by the feedback controller determined in the
previous section, that is, we replace $1_{\omega }^{\ast }U(t)=(1_{\omega
}^{\ast }v(t),1_{\omega }^{\ast }u(t))$ by $-BB^{\ast }R(y(t),z(t)).$ As (%
\ref{31-0})-(\ref{32-0}) is the abstract form of (\ref{23})-(\ref{26}), the
abstract form of the nonlinear system (\ref{19''})-(\ref{22''}) with
replaced right-hand side reads
\begin{eqnarray}
\frac{d}{dt}(y(t),z(t))+\mathcal{A}(y(t),z(t)) &=&\mathcal{G}(y(t))-BB^{\ast
}R(y(t),z(t)),\text{ a.e. }t>0,  \label{82+1} \\
(y(0),z(0)) &=&(y_{0},z_{0}),  \notag
\end{eqnarray}%
where $(y_{0},z_{0})$ is fixed now by (\ref{22+0}), $\mathcal{G}%
(y(t))=(G(y(t)),0)$ and
\begin{equation}
G(y)=\Delta F_{r}(y)+\Delta (g(x)y).  \label{82+2}
\end{equation}%
We recall that $F_{r}$ is the rest of second order of the Taylor expansion
of $F^{\prime }(y+\varphi _{\infty })$ and $g$ is defined by (\ref{g}).
Using the rest in integral form we have
\begin{equation}
F_{r}(y)=y^{2}\int_{0}^{1}(1-s)F^{\prime \prime \prime }(\varphi _{\infty
}+sy)ds=y^{3}+3\varphi _{\infty }y^{2},  \label{88}
\end{equation}%
and assuming that all operations make sense (this will be checked later) we
get
\begin{eqnarray}
G(y) &=&\sum_{j=1}^{7}I_{j},  \label{82+3} \\
I_{1}(y) &=&3y^{2}\Delta y,\text{ }I_{2}(y)=6y\left\vert \nabla y\right\vert
^{2},\text{ }I_{3}(y)=12y\nabla y\cdot \nabla \varphi _{\infty },\text{ }%
I_{4}(y)=3y^{2}\Delta \varphi _{\infty },  \notag \\
I_{5}(y) &=&6\varphi _{\infty }y\Delta y,\text{ }I_{6}(y)=6\varphi _{\infty
}\left\vert \nabla y\right\vert ^{2},\text{ }  \notag \\
I_{7}(y) &=&\Delta (gy)=g\Delta y+y\Delta g+2\nabla y\cdot \nabla g.  \notag
\end{eqnarray}

As usually, in the sequel, $\varphi _{\infty }$\ is the first component of a
stationary solution of the uncontrolled system\textit{\ }(\ref{13})-(\ref{16}%
).\textit{\ }We set
\begin{equation}
\chi _{\infty }:=\left\Vert \nabla \varphi _{\infty }\right\Vert _{\infty
}+\left\Vert \Delta \varphi _{\infty }\right\Vert _{\infty }.  \label{hi}
\end{equation}

\medskip

\noindent \textbf{Theorem 3.1.}\textit{\ There exists }$\chi _{0}>0$\textit{%
\ }(\textit{depending on the problem parameters, the domain and }$\left\Vert
\varphi _{\infty }\right\Vert _{\infty })$ \textit{such that the following
holds true. If }$\chi _{\infty }\leq \chi _{0},$\textit{\ there exists}%
{\huge \ }$\rho $\textit{\ such that for all pairs }$(y_{0},z_{0})\in
D(A^{1/2})\times D(A^{1/4})$\ \textit{with }%
\begin{equation}
\left\Vert y_{0}\right\Vert _{D(A^{1/2})}+\left\Vert z_{0}\right\Vert
_{D(A^{1/4})}\leq \rho ,  \label{ro}
\end{equation}%
\textit{the closed loop system }(\ref{82+1})\textit{\textbf{\ }has a unique
solution }%
\begin{eqnarray}
(y,z) &\in &C([0,\infty );H\times H)\cap L^{2}(0,\infty ;D(A^{3/2})\times
D(A^{3/4}))  \label{y-z} \\
&&\cap W^{1,2}(0,\infty ;(D(A^{1/2})\times D(A^{1/4}))^{\prime }),  \notag
\end{eqnarray}%
\textit{which is exponentially stable, that is}%
\begin{equation}
\left\Vert y(t)\right\Vert _{D(A^{1/2})}+\left\Vert z(t)\right\Vert
_{D(A^{1/4})}\leq C_{P}e^{-kt}(\left\Vert y_{0}\right\Vert
_{D(A^{1/2})}+\left\Vert z_{0}\right\Vert _{D(A^{1/4})}),  \label{83+5}
\end{equation}%
\textit{for some positive constants }$k$\textit{\ and }$C_{P}.$

\medskip

In the previous relations the positive constants $k$\ and $C_{P}$\ depend on
$\Omega ,$ the problem parameters and $\left\Vert \varphi _{\infty
}\right\Vert _{\infty }.$ In addition, $C_{P}$ depends on the full norm $%
\left\Vert \varphi _{\infty }\right\Vert _{2,\infty }.$

\medskip

\noindent \textbf{Proof. }The proof of this theorem will be done in three
steps regarding the existence, uniqueness and stabilization. First,
existence and uniqueness are proved on every interval $[0,T]$ and then they
will be extended to the whole $[0,\infty ).$

\subparagraph{Step 1.}

Existence of the solution is proved on every interval $[0,T]$ by the
Schauder fixed point theorem. Let $(y_{0},z_{0})\in D(A^{1/2})\times
D(A^{1/4}).$

Let $r$ be positive and bounded by a constant which will be specified later.
For $T>0$ arbitrary, but fixed, we introduce the set
\begin{eqnarray}
S_{T} &=&\left\{ (y,z)\in L^{2}(0,T;H\times H);\text{ }\sup_{t\in
(0,T)}\left( \left\Vert y(t)\right\Vert _{D(A^{1/2})}^{2}+\left\Vert
z(t)\right\Vert _{D(A^{1/4})}^{2}\right) \right.  \label{ST} \\
&&\left. +\int_{0}^{T}\left( \left\Vert A^{3/2}y(t)\right\Vert
_{H}^{2}+\left\Vert A^{3/4}z(t)\right\Vert _{H}^{2}\right) dt\leq
r^{2}\right\} .  \notag
\end{eqnarray}%
Let $0<\varepsilon <1/4.$ Clearly, $S_{T}$ is a convex closed subset of $%
L^{2}(0,T;D(A^{3/2-\varepsilon })\times H).$

We fix $(\overline{y},\overline{z})\in S_{T}$ and consider the Cauchy problem%
\begin{eqnarray}
\frac{d}{dt}(y(t),z(t))+\mathcal{A}(y(t),z(t))+BB^{\ast }R(y(t),z(t)) &=&%
\mathcal{G}(\overline{y}(t)),\text{ a.e. }t\in (0,T),  \label{100} \\
(y(0),z(0)) &=&(y_{0},z_{0}).  \notag
\end{eqnarray}%
We prove that such a problem is well-posed and define $\Psi
_{T}:S_{T}\rightarrow L^{2}(0,T;D(A^{3/2-\varepsilon })\times H)$ by $\Psi
_{T}(\overline{y},\overline{z})=(y,z)$ the solution to (\ref{100}). We shall
prove that:

\ i) $\Psi _{T}(S_{T})\subset S_{T}$ provided that $r$ is well chosen;

ii) $\Psi _{T}(S_{T})$ is relatively compact in $L^{2}(0,T;D(A^{3/2-%
\varepsilon })\times H);$

iii) $\Psi _{T}$ is continuous in the $L^{2}(0,T;D(A^{3/2-\varepsilon
})\times H)$ norm.

i) We assert that $G(\overline{y})\in L^{2}(0,T;H\times H),$ relying on the
calculation which shall be made a little later, concluded by (\ref{est-G}).
Then, we recall that $\mathcal{A}+BB^{\ast }R$ is $m$-accretive in $H\times
H,$ as proved in Proposition 2.4, second part, and so, for $(y_{0},z_{0})\in
D(A^{1/2})\times D(A^{1/4})\subset H\times H$ it follows that the Cauchy
problem (\ref{100}) has a unique solution
\begin{equation}
(y,z)\in L^{2}(\delta ,T;D(\mathcal{A}))  \label{exist}
\end{equation}%
with $\delta >0$ arbitrary (see the last part of Proposition 2.4), which
implies that $\mathcal{A}(y(t),z(t))\in H\times H$ a.e. $t>0$ (see \cite%
{brezis-73}, p. 72). Therefore, also (recall (\ref{53})) $y(t)\in D(A^{2})$
and $z(t)\in D(A)$ a.e. Moreover, by Proposition 2.1 we have
\begin{equation}
(y,z)\in C([0,T];H\times H)\cap L^{2}(0,T;D(A)\times D(A^{1/2}))\cap
W^{1,2}([0,T];(D(A)\times D(A^{1/2}))^{\prime }).  \label{sol-Lions}
\end{equation}%
In particular, $(y(t),z(t))\in D(A)\times D(A^{1/2})$ a.e. $t\in (0,T)$ and
so $R(y(t),z(t))\in H\times H$ a.e. $t\in (0,T).$

Next, we have to prove that $(y,z)\in S_{T}$ provided that $r$ is well
chosen. To this end we multiply (\ref{100}) by $R(y(t),z(t))\in H\times H$
scalarly in $H\times H$ and get
\begin{eqnarray*}
&&\frac{1}{2}\frac{d}{dt}(R(y(t),z(t)),(y(t),z(t)))_{H\times H}+(\mathcal{A}%
(y(t),z(t)),R(y(t),z(t)))_{H\times H} \\
&=&-\left\Vert B^{\ast }R(y(t),z(t))\right\Vert _{\mathbb{R}^{N}}^{2}+(%
\mathcal{G}(\overline{y}(t)),R(y(t),z(t)))_{H\times H},\text{ a.e. }t>0.
\end{eqnarray*}%
Therefore, using the Riccati equation (\ref{67}), in the form%
\begin{equation*}
2(\mathcal{A}(y(t),z(t)),R(y(t),z(t)))_{H\times H}+\left\Vert B^{\ast
}R(y(t),z(t))\right\Vert _{\mathbb{R}^{N}}^{2}=\left\Vert
A^{3/2}y(t)\right\Vert _{H}^{2}+\left\Vert A^{3/4}z(t)\right\Vert _{H}^{2}
\end{equation*}%
and recalling (\ref{68}) and (\ref{54+1}) we obtain
\begin{eqnarray*}
&&\frac{1}{2}\frac{d}{dt}(R(y(t),z(t)),(y(t),z(t)))_{H\times H} \\
&&+\frac{1}{2}\left( \left\Vert A^{3/2}y(t)\right\Vert _{H}^{2}+\left\Vert
A^{3/4}z(t)\right\Vert _{H}^{2}+\left\Vert B^{\ast }R(y(t),z(t))\right\Vert
_{\mathbb{R}^{N}}^{2}\right) \\
&\leq &\left\Vert \mathcal{G}(\overline{y}(t))\right\Vert _{H\times
H}\left\Vert R(y(t),z(t))\right\Vert _{H\times H}\leq C_{R}\left\Vert G(%
\overline{y}(t))\right\Vert _{H}\left\Vert (y(t),z(t))\right\Vert
_{D(A)\times D(A^{1/2})} \\
&\leq &CC_{R}\left\Vert G(\overline{y}(t))\right\Vert _{H}\left( \left\Vert
A^{3/2}y(t)\right\Vert _{H}^{2}+\left\Vert A^{3/4}z(t)\right\Vert
_{H}^{2}\right) ^{1/2},\text{ a.e. }t\in (0,T),
\end{eqnarray*}%
with $C_{R}$ from (\ref{68}) and $C$ depending on $\Omega $. Integrating
over $(0,t)$ and using then Young's inequality and (\ref{69}) we
successively get%
\begin{eqnarray}
&&c_{1}^{\prime }(\left\Vert y(t)\right\Vert _{D(A^{1/2})}^{2}+\left\Vert
z(t)\right\Vert _{D(A^{1/4})}^{2})+\int_{0}^{t}\left( \left\Vert
A^{3/2}y(s)\right\Vert _{H}^{2}+\left\Vert A^{3/4}z(s)\right\Vert
_{H}^{2}\right) ds  \notag \\
&\leq &c_{2}^{\prime }(\left\Vert y_{0}\right\Vert
_{D(A^{1/2})}^{2}+\left\Vert z_{0}\right\Vert
_{D(A^{1/4})}^{2})+C_{R}^{\prime }\int_{0}^{t}\left\Vert G(\overline{y}%
(s))\right\Vert _{H}^{2}ds,\text{ }t\in (0,T),  \label{103}
\end{eqnarray}%
where $c_{1}^{\prime },$ $c_{2}^{\prime },$ $C_{R}^{\prime }$ are
proportional to $c_{1},$ $c_{2},$ $C_{R}^{2}.$ Now let us prove that $G(%
\overline{y})\in L^{2}(0,T;H)$ and for that we shall estimate each term $%
I_{j}$ in (\ref{82+3}) and show that $I_{j}(\overline{y})\in L^{2}(0,T;H),$
for $j=1,...,7.$ In the computations below we shall use the interpolation
inequalities (\ref{54+0})-(\ref{100-0}), and (\ref{93}). The constants we
shall introduce do not depend on $\varphi _{\infty }.$ By (\ref{82+3}) we
compute%
\begin{eqnarray*}
\left\Vert I_{1}(\overline{y})\right\Vert _{H}^{2} &=&C\left\Vert \overline{y%
}^{2}\Delta \overline{y}\right\Vert _{H}^{2}=C\int_{\Omega }\overline{y}%
^{4}(\Delta \overline{y})^{2}dx\leq C\left( \int_{\Omega }\overline{y}%
^{8}dx\right) ^{1/2}\left( \int_{\Omega }(\Delta \overline{y})^{4}dx\right)
^{1/2} \\
&=&C\left\Vert \overline{y}\right\Vert _{L^{8}(\Omega )}^{4}\left\Vert
\Delta \overline{y}\right\Vert _{L^{4}(\Omega )}^{2}\leq C\left\Vert
\overline{y}\right\Vert _{H^{\alpha _{1}}(\Omega )}^{4}\left\Vert \Delta
\overline{y}\right\Vert _{H^{\alpha _{2}}(\Omega )}^{2} \\
&\leq &C\left\Vert \overline{y}\right\Vert _{H^{\alpha _{1}}(\Omega
)}^{4}\left\Vert \overline{y}-A\overline{y}\right\Vert _{H^{\alpha
_{2}}(\Omega )}^{2}\leq C\left\Vert A^{0}\overline{y}\right\Vert _{H^{\alpha
_{1}}(\Omega )}^{4}\left\Vert A\overline{y}\right\Vert _{H^{\alpha
_{2}}(\Omega )}^{2} \\
&\leq &C\left\Vert A^{\alpha _{1}/2}\overline{y}\right\Vert
_{H}^{4}\left\Vert A^{1+\alpha _{2}/2}\overline{y}\right\Vert _{H}^{2},
\end{eqnarray*}%
where $\alpha _{1}\geq \frac{9}{8}$ and $\alpha _{2}\geq \frac{3}{4}.$
Furthermore, we have
\begin{eqnarray}
\left\Vert I_{1}(\overline{y})\right\Vert _{H}^{2} &\leq &C\left( \left\Vert
A^{3/2}\overline{y}\right\Vert _{H}^{(\alpha _{1}-1)/2}\left\Vert A^{1/2}%
\overline{y}\right\Vert _{H}^{(3-\alpha _{1})/2}\right) ^{4}\left(
\left\Vert A^{3/2}\overline{y}\right\Vert _{H}^{(\alpha _{2}+1)/2}\left\Vert
A^{1/2}\overline{y}\right\Vert _{H}^{(1-\alpha _{2})/2}\right) ^{2}  \notag
\\
&=&C\left\Vert A^{3/2}\overline{y}\right\Vert _{H}^{2\alpha _{1}+\alpha
_{2}-1}\left\Vert A^{1/2}\overline{y}\right\Vert _{H}^{7-(2\alpha
_{1}+\alpha _{2})},\text{ }\alpha _{1}\geq \frac{9}{8},\text{ }\alpha
_{2}\geq \frac{3}{4}.  \label{I1}
\end{eqnarray}%
Since $\left\Vert A^{1/2}\overline{y}\right\Vert _{H}^{7-(2\alpha
_{1}+\alpha _{2})}=\left\Vert A^{1/2}\overline{y}\right\Vert
_{H}^{3-(2\alpha _{1}+\alpha _{2})}\left\Vert A^{1/2}\overline{y}\right\Vert
_{H}^{4}\leq C\left\Vert A^{3/2}\overline{y}\right\Vert _{H}^{3-(2\alpha
_{1}+\alpha _{2})}\left\Vert A^{1/2}\overline{y}\right\Vert _{H}^{4},$ we get%
\begin{equation*}
\left\Vert I_{1}(\overline{y})\right\Vert _{H}^{2}\leq C\left\Vert A^{3/2}%
\overline{y}\right\Vert _{H}^{2}\left\Vert A^{1/2}\overline{y}\right\Vert
_{H}^{4}.
\end{equation*}%
But $(\overline{y},\overline{z})\in S_{T}$ and therefore
\begin{equation}
\int_{0}^{T}\left\Vert I_{1}(\overline{y})\right\Vert _{H}^{2}dt\leq
Cr^{4}\int_{0}^{T}\left\Vert A^{3/2}\overline{y}\right\Vert _{H}^{2}dt\leq
Cr^{6}.  \label{I1-1}
\end{equation}%
For the second term we infer that
\begin{eqnarray}
\left\Vert I_{2}(\overline{y})\right\Vert _{H}^{2} &=&C\left\Vert \overline{y%
}\left\vert \nabla \overline{y}\right\vert ^{2}\right\Vert
_{H}^{2}=C\int_{\Omega }\overline{y}^{2}\left\vert \nabla \overline{y}%
\right\vert ^{4}dx\leq C\left( \int_{\Omega }\overline{y}^{4}dx\right)
^{1/2}\left( \int_{\Omega }\left\vert \nabla \overline{y}\right\vert
^{8}dx\right) ^{1/2}  \notag \\
&=&C\left\Vert \overline{y}\right\Vert _{L^{4}(\Omega )}^{2}\left\Vert
\nabla \overline{y}\right\Vert _{L^{8}(\Omega )}^{4}\leq C\left\Vert
\overline{y}\right\Vert _{H^{\alpha _{1}}(\Omega )}^{2}\left\Vert \nabla
\overline{y}\right\Vert _{H^{\alpha _{2}}(\Omega )}^{4}  \label{I2} \\
&\leq &C\left\Vert A^{\alpha _{1}/2}\overline{y}\right\Vert
_{H}^{2}\left\Vert A^{1/2+\alpha _{2}/2}\overline{y}\right\Vert _{H}^{4}
\notag \\
&\leq &C\left( \left\Vert A^{3/2}\overline{y}\right\Vert _{H}^{(\alpha
_{1}-1)/2}\left\Vert A^{1/2}\overline{y}\right\Vert _{H}^{(3-\alpha
_{1})/2}\right) ^{2}\left( \left\Vert A^{3/2}\overline{y}\right\Vert
_{H}^{\alpha _{2}/2}\left\Vert A^{1/2}\overline{y}\right\Vert
_{H}^{(2-\alpha _{2})/2}\right) ^{4}  \notag \\
&\leq &C\left\Vert A^{3/2}\overline{y}\right\Vert _{H}^{\alpha _{1}+2\alpha
_{2}-1}\left\Vert A^{1/2}\overline{y}\right\Vert _{H}^{7-(\alpha
_{1}+2\alpha _{2})},\text{ \ }\alpha _{1}\geq \frac{3}{4},\text{ }\alpha
_{2}\geq \frac{9}{8}.  \notag
\end{eqnarray}%
Thus, by choosing $\alpha _{1}=3/4$ and $\alpha _{2}=9/8,$ we have
\begin{equation}
\left\Vert I_{2}(\overline{y})\right\Vert _{H}^{2}\leq C\left\Vert A^{3/2}%
\overline{y}\right\Vert _{H}^{2}\left\Vert A^{1/2}\overline{y}\right\Vert
_{H}^{4}.  \label{I2-1}
\end{equation}%
For $I_{3}$ to $I_{6}$ we have the following estimates:
\begin{eqnarray}
&&\left\Vert I_{3}(\overline{y})\right\Vert _{H}^{2}=C\left\Vert \overline{y}%
\nabla \overline{y}\cdot \nabla \varphi _{\infty }\right\Vert _{H}^{2}\leq
C\left\Vert \nabla \varphi _{\infty }\right\Vert _{\infty }^{2}\left\Vert
\overline{y}\right\Vert _{L^{4}(\Omega )}^{2}\left\Vert \nabla \overline{y}%
\right\Vert _{L^{4}(\Omega )}^{2}  \label{I3} \\
&\leq &C\left\Vert \nabla \varphi _{\infty }\right\Vert _{\infty
}^{2}\left\Vert \overline{y}\right\Vert _{H^{1}(\Omega )}^{2}\left\Vert
\overline{y}\right\Vert _{H^{2}(\Omega )}^{2}\leq C\left\Vert \nabla \varphi
_{\infty }\right\Vert _{\infty }^{2}\left\Vert A^{1/2}\overline{y}%
\right\Vert _{H}^{2}\left\Vert A\overline{y}\right\Vert _{H}^{2}  \notag \\
&\leq &C\left\Vert \nabla \varphi _{\infty }\right\Vert _{\infty
}^{2}\left\Vert A^{3/2}\overline{y}\right\Vert _{H}^{2}\left\Vert A^{1/2}%
\overline{y}\right\Vert _{H}^{2},\text{ }  \notag
\end{eqnarray}%
\begin{eqnarray}
&&\left\Vert I_{4}(\overline{y})\right\Vert _{H}^{2}=C\left\Vert \overline{y}%
^{2}\Delta \varphi _{\infty }\right\Vert _{H}^{2}\leq C\left\Vert \Delta
\varphi _{\infty }\right\Vert _{\infty }^{2}\left\Vert \overline{y}%
\right\Vert _{L^{4}(\Omega )}^{4}\leq C\left\Vert \Delta \varphi _{\infty
}\right\Vert _{\infty }^{2}\left\Vert \overline{y}\right\Vert _{H^{1}(\Omega
)}^{4}  \label{I4} \\
&\leq &C\left\Vert \Delta \varphi _{\infty }\right\Vert _{\infty
}^{2}\left\Vert A^{1/2}\overline{y}\right\Vert _{H}^{4}\leq C\left\Vert
\Delta \varphi _{\infty }\right\Vert _{\infty }^{2}\left\Vert A^{3/2}%
\overline{y}\right\Vert _{H}^{2}\left\Vert A^{1/2}\overline{y}\right\Vert
_{H}^{2},  \notag
\end{eqnarray}%
\begin{eqnarray}
&&\left\Vert I_{5}(\overline{y})\right\Vert _{H}^{2}=C\left\Vert \varphi
_{\infty }\overline{y}\Delta \overline{y}\right\Vert _{H}^{2}\leq
C\left\Vert \varphi _{\infty }\right\Vert _{\infty }^{2}\left\Vert \overline{%
y}\right\Vert _{L^{4}(\Omega )}^{2}\left\Vert \Delta \overline{y}\right\Vert
_{L^{4}(\Omega )}^{2}  \label{I5} \\
&\leq &C\left\Vert \varphi _{\infty }\right\Vert _{\infty }^{2}\left\Vert
\overline{y}\right\Vert _{H^{1}(\Omega )}^{2}\left\Vert \Delta \overline{y}%
\right\Vert _{H^{1}(\Omega )}^{2}\leq C\left\Vert \varphi _{\infty
}\right\Vert _{\infty }^{2}\left\Vert A^{1/2}\overline{y}\right\Vert
_{H}^{2}\left\Vert A^{3/2}\overline{y}\right\Vert _{H}^{2},  \notag
\end{eqnarray}%
\begin{eqnarray}
&&\left\Vert I_{6}(\overline{y})\right\Vert _{H}^{2}=C\left\Vert \varphi
_{\infty }\left\vert \nabla \overline{y}\right\vert ^{2}\right\Vert
_{H}^{2}\leq C\left\Vert \varphi _{\infty }\right\Vert _{\infty
}^{2}\left\Vert \nabla \overline{y}\right\Vert _{L^{4}(\Omega )}^{4}\leq
C\left\Vert \varphi _{\infty }\right\Vert _{\infty }^{2}\left\Vert \overline{%
y}\right\Vert _{H^{2}(\Omega )}^{4}  \label{I6} \\
&\leq &C\left\Vert \varphi _{\infty }\right\Vert _{\infty }^{2}\left\Vert A%
\overline{y}\right\Vert _{H}^{4}\leq C\left\Vert \varphi _{\infty
}\right\Vert _{\infty }^{2}\left\Vert A^{3/2}\overline{y}\right\Vert
_{H}^{2}\left\Vert A^{1/2}\overline{y}\right\Vert _{H}^{2},  \notag
\end{eqnarray}%
and finally%
\begin{eqnarray}
&&\left\Vert I_{7}(\overline{y})\right\Vert _{H}^{2}=\left\Vert g\Delta
\overline{y}+\overline{y}\Delta g+2\nabla \overline{y}\cdot \nabla
g\right\Vert _{H}^{2}  \label{I7} \\
&\leq &C\left\Vert g\right\Vert _{2,\infty }^{2}\left( \left\Vert A\overline{%
y}\right\Vert _{H}^{2}+\left\Vert \overline{y}\right\Vert
_{H}^{2}+\left\Vert A^{1/2}\overline{y}\right\Vert _{H}^{2}\right) \leq
C\left\Vert g\right\Vert _{2,\infty }^{2}\left\Vert A^{3/2}\overline{y}%
\right\Vert _{H}^{2}.  \notag
\end{eqnarray}%
Moreover, by (\ref{g}) we get
\begin{equation*}
\left\vert g(x)\right\vert \leq \frac{6\left\Vert \varphi _{\infty
}\right\Vert _{\infty }}{m_{\Omega }}\int_{\Omega }\left\vert \varphi
_{\infty }(x)-\varphi _{\infty }(\xi )\right\vert d\xi \leq 6\left\Vert
\varphi _{\infty }\right\Vert _{\infty }\left\Vert \nabla \varphi _{\infty
}\right\Vert _{\infty }d_{\Omega },
\end{equation*}%
\begin{equation}
\left\Vert g\right\Vert _{\infty }\leq C\left\Vert \varphi _{\infty
}\right\Vert _{\infty }\left\Vert \nabla \varphi _{\infty }\right\Vert
_{\infty },  \label{ginf}
\end{equation}%
with $d_{\Omega }$ the supremum of the geodesic distance of $\Omega .$ Next,%
\begin{equation}
\nabla g=6\varphi _{\infty }\nabla \varphi _{\infty },\text{ \ }\left\Vert
\nabla g\right\Vert _{\infty }\leq C\left\Vert \varphi _{\infty }\right\Vert
_{\infty }\left\Vert \nabla \varphi _{\infty }\right\Vert _{\infty }
\label{grad-g}
\end{equation}%
\begin{equation}
\Delta g=6\varphi _{\infty }\Delta \varphi _{\infty }+6\left\vert \nabla
\varphi _{\infty }\right\vert ^{2},\text{ \ }\left\Vert \Delta g\right\Vert
_{\infty }\leq C(\left\Vert \varphi _{\infty }\right\Vert _{\infty
}\left\Vert \Delta \varphi _{\infty }\right\Vert _{\infty }+\left\Vert
\nabla \varphi _{\infty }\right\Vert _{\infty }^{2}),  \label{lapl-g}
\end{equation}%
whence
\begin{equation}
\left\Vert g\right\Vert _{2,\infty }\leq C_{\Omega }\overline{g}_{\infty }
\label{g2inf}
\end{equation}%
where
\begin{equation}
\overline{g}_{\infty }=\left\Vert \varphi _{\infty }\right\Vert _{\infty
}\left\Vert \nabla \varphi _{\infty }\right\Vert _{\infty }+\left\Vert
\varphi _{\infty }\right\Vert _{\infty }\left\Vert \Delta \varphi _{\infty
}\right\Vert _{\infty }+\left\Vert \nabla \varphi _{\infty }\right\Vert
_{\infty }^{2},  \label{hi-inf}
\end{equation}%
with $C_{\Omega }$ a constant dependent on the domain $\Omega $. Finally,
collecting all the estimates above and recalling (\ref{82+3}) and (\ref{hi}%
), we obtain for $(\overline{y},\overline{z})\in S_{T}$ that
\begin{eqnarray*}
&&\int_{0}^{T}\left\Vert G(\overline{y}(t))\right\Vert
_{H}^{2}dt=\sum_{j=1}^{7}\int_{0}^{T}\left\Vert I_{j}(\overline{y}%
)\right\Vert _{H}^{2}dt=C\int_{0}^{T}2\left\Vert A^{3/2}\overline{y}%
\right\Vert _{H}^{2}\left\Vert A^{1/2}\overline{y}\right\Vert _{H}^{4}dt \\
&&+C\int_{0}^{T}\left( \left\Vert \nabla \varphi _{\infty }\right\Vert
_{\infty }^{2}+\left\Vert \Delta \varphi _{\infty }\right\Vert _{\infty
}^{2}+2\left\Vert \varphi _{\infty }\right\Vert _{\infty }^{2}\right)
\left\Vert A^{3/2}\overline{y}\right\Vert _{H}^{2}\left\Vert A^{1/2}%
\overline{y}\right\Vert _{H}^{2}dt \\
&&+C\left\Vert g\right\Vert _{2,\infty }^{2}\int_{0}^{T}\left\Vert A^{3/2}%
\overline{y}\right\Vert _{H}^{2}dt.
\end{eqnarray*}%
In view of (\ref{ST}), we conclude that
\begin{equation}
\int_{0}^{T}\left\Vert G(\overline{y}(t))\right\Vert _{H}^{2}dt\leq
C(r^{6}+\left\Vert \varphi _{\infty }\right\Vert _{2,\infty
}^{2}r^{4}+\left\Vert g\right\Vert _{2,\infty }^{2}r^{2}),  \label{est-G}
\end{equation}%
where we stress that $C$ is a constant independent of $\varphi _{\infty }.$
Going back to (\ref{103}), we can write it in the form
\begin{eqnarray}
&&\left\Vert y(t)\right\Vert _{D(A^{1/2})}^{2}+\left\Vert z(t)\right\Vert
_{D(A^{1/4})}^{2}+\int_{0}^{t}\left( \left\Vert A^{3/2}y(s)\right\Vert
_{H}^{2}+\left\Vert A^{3/4}z(s)\right\Vert _{H}^{2}\right) ds  \label{130} \\
&\leq &C_{P}\left( \left\Vert y_{0}\right\Vert _{D(A^{1/2})}^{2}+\left\Vert
z_{0}\right\Vert _{D(A^{1/4})}^{2}+\int_{0}^{t}\left\Vert G(\overline{y}%
(s))\right\Vert _{H}^{2}ds\right)  \notag
\end{eqnarray}%
where $C_{P}$ is a constant depending on $\nu ,\gamma ,l$ and $\left\Vert
\varphi _{\infty }\right\Vert _{L^{2}(\Omega )},$ and we would like to
impose that the left-hand side is $\leq r^{2}.$ To this aim it suffices that
the right-hand side is $\leq r^{2}.$ On account of (\ref{est-G}) we see that
the latter condition holds provided that
\begin{equation*}
\left\Vert y_{0}\right\Vert _{D(A^{1/2})}+\left\Vert z_{0}\right\Vert
_{D(A^{1/4})}\leq \rho \text{ }
\end{equation*}%
and
\begin{equation*}
C_{P}\rho ^{2}+CC_{P}(r^{6}+\left\Vert \varphi _{\infty }\right\Vert
_{2,\infty }^{2}r^{4}+\left\Vert g\right\Vert _{2,\infty }^{2}r^{2})\leq
r^{2},
\end{equation*}%
where $C_{P}$ and $C$ are the precise constants entering in (\ref{130}) and (%
\ref{est-G}). We notice that the first condition coincides with (\ref{ro}).
Then, if we assume that
\begin{equation}
C_{P}\rho ^{2}\leq \frac{1}{2}r^{2}\text{ \ that is, \ }\rho \leq
(2C_{P})^{-1/2}r  \label{130-2}
\end{equation}%
a sufficient condition for our bound is
\begin{equation*}
r^{4}+\left\Vert \varphi _{\infty }\right\Vert _{2,\infty
}^{2}r^{2}+\left\Vert g\right\Vert _{2,\infty }^{2}-\frac{1}{2CC_{P}}\leq 0.
\end{equation*}%
This is satisfied provided that
\begin{equation}
\left\Vert g\right\Vert _{2,\infty }^{2}\leq (2C_{P}C)^{-1/2}\text{ \ and }%
r\leq r_{1}  \label{130-0}
\end{equation}%
where $r_{1}>0$ is given by%
\begin{equation}
r_{1}^{2}=\frac{-\left\Vert \varphi _{\infty }\right\Vert _{2,\infty }^{2}+%
\sqrt{D_{\infty }}}{2}  \label{136}
\end{equation}%
with
\begin{equation*}
D_{\infty }=\left\Vert \varphi _{\infty }\right\Vert _{2,\infty
}^{4}-4(\left\Vert g\right\Vert _{2,\infty }^{2}-(2C_{P}C)^{-1}).
\end{equation*}%
Notice that the first condition in (\ref{130-0}) implies that $D_{\infty }$
is nonnegative and that $r_{1}$ is well-defined. Both $D_{\infty }$ and $%
r_{1}$ depend on the full norm $\left\Vert \varphi _{\infty }\right\Vert
_{2,\infty },$ i.e., on $\left\Vert g_{\infty }\right\Vert _{\infty }$ and $%
\chi _{\infty },$ but they are independent of $T.$

Now, we look for a sufficient condition for it. By (\ref{g2inf})-(\ref%
{hi-inf}) and (\ref{hi}) we have
\begin{equation*}
\left\Vert g\right\Vert _{2,\infty }\leq C_{\Omega }\overline{g}_{\infty
}\leq C_{\Omega }\left( \left\Vert \varphi _{\infty }\right\Vert _{\infty
}\chi _{\infty }+\chi _{\infty }^{2}\right) .
\end{equation*}%
Hence, the first inequality in (\ref{130-0}) holds if $C_{\Omega }\left(
\left\Vert \varphi _{\infty }\right\Vert _{\infty }\chi _{\infty }+\chi
_{\infty }^{2}\right) \leq (2C_{P}C)^{-1/2},$ that is
\begin{equation*}
\chi _{\infty }^{2}+\left\Vert \varphi _{\infty }\right\Vert _{\infty }\chi
_{\infty }-C_{\Omega }^{-1}(2C_{P}C)^{-1/2}\leq 0.
\end{equation*}%
But this is true whenever
\begin{equation}
\chi _{\infty }\leq \chi _{0}^{\prime }:=\frac{-\left\Vert \varphi _{\infty
}\right\Vert _{\infty }+\sqrt{\left\Vert \varphi _{\infty }\right\Vert
_{\infty }^{2}+4C_{\Omega }^{-1}(2C_{P}C)^{-1/2}}}{2}.  \label{hi-0}
\end{equation}%
We stress that $\chi _{0}^{\prime }$ depends on $\nu ,l,\gamma ,\Omega $ and
$\left\Vert \varphi _{\infty }\right\Vert _{\infty },$ but it is independent
of $T.$

So, if we assume (\ref{hi})-(\ref{ro}) with $\chi _{0}^{\prime }$ given by (%
\ref{hi-0}) and $\rho $ and $r$ satisfying (\ref{130-2}) and the second
constraint in (\ref{130-0}), coming back to (\ref{130}) we obtain
\begin{equation}
\left\Vert y(t)\right\Vert _{D(A^{1/2})}^{2}+\left\Vert z(t)\right\Vert
_{D(A^{1/4})}^{2}+\int_{0}^{t}\left( \left\Vert A^{3/2}y(s)\right\Vert
_{H}^{2}+\left\Vert A^{3/4}z(s)\right\Vert _{H}^{2}\right) ds\leq r^{2}
\label{131}
\end{equation}%
for all $t\in \lbrack 0,T].$

In conclusion, if we fix $r$ and $\rho $ satisfying (\ref{130-0}) and (\ref%
{130-2}), and if $\left\Vert \nabla \varphi _{\infty }\right\Vert _{\infty
}+\left\Vert \Delta \varphi _{\infty }\right\Vert _{\infty }$ is small as
specified above we have that $(y,z)\in S_{T},$ hence $\Psi _{T}$ maps $S_{T}$
into $S_{T}.$ It is understood that $r$ and $\rho $ are fixed in the sequel
(and do not depend on $T).$

ii) Let $(y,z)=\Psi _{T}(\overline{y},\overline{z}),$ with $(\overline{y},%
\overline{z})\in S_{T}.$ We observe that $(y,z)$ and $\frac{d}{dt}(y,z)$
remain bounded in $L^{2}(0,T;D(A^{3/2})\times D(A^{3/4}))$ and $%
W^{1,2}([0,T];(D(A)\times D(A^{1/2}))^{\prime }),$ respectively. Here we
used Proposition 2.1 (replacing the term $1_{\omega }^{\ast }U(t)$ in (\ref%
{31-0}) by $\mathcal{G}(\overline{y}(t))$ in (\ref{100})) and the estimate (%
\ref{32-1-0}) in which the last term on the right-hand side, $%
\int_{0}^{T}\left\Vert 1_{\omega }^{\ast }U(s)\right\Vert _{\mathcal{H}%
}^{2}ds,$ is replaced by the integral $\int_{0}^{T}\left\Vert G(\overline{y}%
(t))\right\Vert _{H}^{2}dt$ in (\ref{est-G}). Since $D(A^{3/2})\times
D(A^{3/4})$ is compactly embedded in $D(A^{3/2-\varepsilon })\times H$ it
follows by Lions-Aubin lemma (see \cite{lions}, p. 58) that the set $\Psi
_{T}(S_{T})$ is relatively compact in $L^{2}(0,T;D(A^{3/2-\varepsilon
})\times H).$

iii) Let $(\overline{y_{n}},\overline{z_{n}})\in S_{T}$, $(\overline{y_{n}},%
\overline{z_{n}})\rightarrow (\overline{y},\overline{z})$ strongly in $%
L^{2}(0,T;D(A^{3/2-\varepsilon })\times H),$ as $n\rightarrow \infty .$ We
have to prove that
\begin{equation*}
\Psi _{T}(\overline{y_{n}},\overline{z_{n}})\rightarrow \Psi _{T}(\overline{y%
},\overline{z})\text{ strongly in }L^{2}(0,T;D(A^{3/2-\varepsilon })\times
H).
\end{equation*}%
The solution $(y_{n},z_{n})$ to (\ref{100}) corresponding to $(\overline{%
y_{n}},\overline{z_{n}})$ belongs to the spaces specified in (\ref{sol-Lions}%
) and satisfies the estimate (\ref{131}). Also, $\left\{ \frac{d}{dt}%
(y_{n},z_{n})\right\} _{n}$ is bounded in $L^{2}(0,T;(D(A)\times
D(A^{1/2}))^{\prime }).$ Hence, on a subsequence $\{n\rightarrow \infty \}$
we have
\begin{equation*}
y_{n}\rightarrow y\text{ weakly in }L^{2}(0,T;D(A^{3/2}))\text{ and weak* in
}L^{\infty }(0,T;D(A^{1/2})),
\end{equation*}%
\begin{equation*}
z_{n}\rightarrow z\text{ weakly in }L^{2}(0,T;D(A^{3/4}))\text{ and weak* in
}L^{\infty }(0,T;D(A^{1/4})),
\end{equation*}%
\begin{equation*}
\frac{dy_{n}}{dt}\rightarrow \frac{dy}{dt}\text{ weakly in }%
L^{2}(0,T;(D(A))^{\prime }),
\end{equation*}%
\begin{equation*}
\frac{dz_{n}}{dt}\rightarrow \frac{dz}{dt}\text{ weakly in }%
L^{2}(0,T;(D(A^{1/2}))^{\prime }).
\end{equation*}
Then, by Aubin-Lions lemma%
\begin{equation*}
(y_{n},z_{n})\rightarrow (y,z)\text{ strongly in }L^{2}(0,T;D(A^{3/2-%
\varepsilon })\times H).
\end{equation*}%
Let us to show that
\begin{equation*}
G(\overline{y_{n}})\rightarrow G(\overline{y})\text{ weakly in }L^{2}(0,T;H),
\end{equation*}%
by treating the terms in (\ref{82+3})$.$ Taking into account that $%
D(A^{3/2-\varepsilon })\subset H^{3-2\varepsilon }(\Omega )$ and that $%
\overline{y_{n}}\rightarrow \overline{y}$ strongly in $L^{2}(0,T;D(A^{3/2-%
\varepsilon }))$ we infer that
\begin{equation*}
\overline{y_{n}}\rightarrow \overline{y},\text{ \ }\Delta \overline{y_{n}}%
\rightarrow \Delta \overline{y}\text{ strongly in }L^{2}(0,T;L^{2}(\Omega )),
\end{equation*}%
\begin{equation*}
\nabla \overline{y_{n}}\rightarrow \nabla \overline{y}\text{ strongly in }%
(L^{2}(0,T;L^{2}(\Omega )))^{d}
\end{equation*}%
and so, on a subsequence they tend a.e. in $Q.$ This implies that
\begin{equation*}
I_{j}(\overline{y_{n}})\rightarrow I_{j}(\overline{y}),\text{ a.e. on }Q,%
\text{ for }j=1,...,7.
\end{equation*}%
On the other hand, by estimates (\ref{I1-1})-(\ref{I7}) we conclude  that $%
\{I_{j}(\overline{y_{n}})\}_{n}$ is bounded in $L^{2}(0,T;H),$ and selecting
a subsequence we have
\begin{equation*}
I_{j}(\overline{y_{n}})\rightarrow \zeta _{j}\text{ weakly in }L^{2}(0,T;H),%
\text{ for }j=1,...,7.
\end{equation*}%
Thus, we deduce that $\zeta _{j}=I_{j}(\overline{y}),$ a.e. on $Q,$ for $%
j=1,...,7.$

Moreover, writing the weak form of (\ref{100}) corresponding to $(\overline{%
y_{n}},\overline{z_{n}})$ and passing to the limit we get that $(y,z)=\Psi
_{T}(\overline{y},\overline{z}).$ As the same holds for any subsequence this
ends the proof of the continuity of $\Psi _{T}.$

Then, by the Schauder fixed point theorem, applied to the mapping $\Psi _{T}$
on the space $L^{2}(0,T;D(A^{3/2-\varepsilon })\times H),$ it follows that
problem (\ref{100}) has at least a solution on the interval $[0,T],$ $%
(y,z)\in S_{T}.$

\subparagraph{Step 2.}

We prove here the uniqueness of the solution on $[0,T]$. Let us consider the
nonlinear system (\ref{13})-(\ref{15}) written in terms of $\varphi $ and $%
\sigma ,$ where by (\ref{control-fin})
\begin{equation}
(1_{\omega }^{\ast }v(t),1_{\omega }^{\ast }u(t))=-BB^{\ast }R(y(t),z(t)),%
\text{ }t\in (0,T).  \label{6p}
\end{equation}%
For the uniqueness proof we need to prove that
\begin{equation}
BB^{\ast }\text{ is linear continuous from }V^{\prime }\times V^{\prime
}\rightarrow V^{\prime }\times V^{\prime },  \label{7p}
\end{equation}%
while we already know that
\begin{equation}
R\text{ is linear continuous from }D(A^{1/2})\times D(A^{1/4})\rightarrow
(D(A^{1/2})\times D(A^{1/4}))^{\prime },  \label{8p}
\end{equation}%
the latter being known by the definition of $R$ (according to (\ref{64-0})).

In what concerns (\ref{7p}), using the definitions of $B$ and $B^{\ast },$
we have that if $q=(q_{1},q_{2})\in V^{\prime }\times V^{\prime }$, then%
\begin{equation}
BB^{\ast }q=\left[
\begin{tabular}{c}
$\sum\limits_{i=1}^{N}1_{\omega }^{\ast }\varphi _{i}(\left\langle
q_{1},1_{\omega }^{\ast }\varphi _{i}\right\rangle +\left\langle
q_{2},1_{\omega }^{\ast }\psi _{i}\right\rangle )$ \\
$\sum\limits_{i=1}^{N}1_{\omega }^{\ast }\psi _{i}(\left\langle
q_{1},1_{\omega }^{\ast }\varphi _{i}\right\rangle +\left\langle
q_{2},1_{\omega }^{\ast }\psi _{i}\right\rangle )$%
\end{tabular}%
\right]  \label{8'p}
\end{equation}%
is well defined because $1_{\omega }^{\ast }$ is a multiplicator in $V.$
Then, it is easily seen that%
\begin{equation}
\left\Vert BB^{\ast }q\right\Vert _{V^{\prime }\times V^{\prime }}\leq
C\left\Vert q\right\Vert _{V^{\prime }\times V^{\prime }}\text{ for }q\in
V^{\prime }\times V^{\prime }.  \label{8''p}
\end{equation}%
Now, we rewrite (\ref{13})-(\ref{14}) in terms of the operator $A$ and have
\begin{equation}
\varphi _{t}+\nu A^{2}\varphi +A(\varphi ^{3}+(l-1-2\nu )\varphi -\gamma
\sigma )-\varphi ^{3}-(l-1-\nu )\varphi +\gamma \sigma =1_{\omega }^{\ast }v,
\label{9p}
\end{equation}%
\begin{equation}
\sigma _{t}+A(\sigma -\gamma \varphi )-\sigma +\gamma \varphi =1_{\omega
}^{\ast }u,  \label{10p}
\end{equation}%
with the boundary and initial conditions (\ref{16})-(\ref{15}).

Assume that there are two solutions $(\varphi ^{i},\sigma ^{i}),$ $=1,2,$
corresponding to $U_{i}=(v_{i},u_{i})$ with $1_{\omega }^{\ast
}U_{i}=-BB^{\ast }R(y_{i},z_{i}),$ $i=1,2.$ We take the difference of the
equations (\ref{9p}) and test it by $A^{-1}(\varphi ^{1}-\varphi ^{2}).$
Then, test the difference of equations (\ref{10p}) by $\lambda (\sigma
^{1}-\sigma ^{2}),$ where $\lambda >0$ is a coefficient to be chosen later.
We use the simplified notation $\varphi =\varphi ^{1}-\varphi ^{2},$ $\sigma
=\sigma ^{1}-\sigma ^{2},$ $v=v_{1}-v_{2},$ $u=u_{1}-u_{2}.$

Moreover, we see the operator $A$ also from $V$ to $V^{\prime }$ defined by
\begin{equation*}
\left\langle Aw,\psi \right\rangle _{V^{\prime },V}=\int_{\Omega }(\nabla
w\cdot \nabla \psi +w\psi )dx\text{ \ for }\psi \in V.
\end{equation*}

For the first computation we obtain%
\begin{eqnarray}
&&\frac{1}{2}\frac{d}{dt}\left\Vert \varphi (t)\right\Vert _{V^{\prime
}}^{2}+\nu \left\Vert \varphi (t)\right\Vert _{V}^{2}+\int_{\Omega }\left(
(\varphi ^{1})^{3}-(\varphi ^{2})^{3}\right) (\varphi ^{1}-\varphi ^{2})dx
\label{11p} \\
&\leq &(2\nu +1-l)\int_{\Omega }\varphi ^{2}dx+\gamma \int_{\Omega }\sigma
(\varphi -A^{-1}\varphi )dx  \notag \\
&&+\int_{\Omega }((\varphi ^{1})^{3}-(\varphi ^{2})^{3})A^{-1}(\varphi
^{1}-\varphi ^{2})dx+\int_{\Omega }(l-1-\nu )\varphi A^{-1}\varphi dx  \notag
\\
&&+\left\langle 1_{\omega }^{\ast }v(t),A^{-1}\varphi (t)\right\rangle
_{V^{\prime },V},  \notag
\end{eqnarray}%
where we have used the property that $\left\langle A\varphi ,\varphi
\right\rangle _{V^{\prime },V}=\left\Vert \varphi \right\Vert _{V}^{2}$ and $%
\left\langle w,A^{-1}w\right\rangle _{V^{\prime },V}=\left\Vert w\right\Vert
_{V^{\prime }}^{2}$ for $w\in V^{\prime }.$

Let us estimate each term on the right-hand side:%
\begin{equation*}
(2\nu +1-l)\int_{\Omega }\varphi ^{2}dx\leq \left\vert 2\nu +1-l\right\vert
\left\Vert \varphi (t)\right\Vert _{V^{\prime }}\left\Vert \varphi
(t)\right\Vert _{V}\leq \frac{\nu }{8}\left\Vert \varphi (t)\right\Vert
_{V}^{2}+C\left\Vert \varphi (t)\right\Vert _{V^{\prime }}^{2}.
\end{equation*}%
By the Young inequality%
\begin{eqnarray*}
\gamma \int_{\Omega }\sigma (\varphi -A^{-1}\varphi )dx &\leq &\gamma
\left\Vert \sigma (t)\right\Vert _{V^{\prime }}(\left\Vert \varphi
(t)\right\Vert _{V}+\left\Vert A^{-1}\varphi (t)\right\Vert _{V}) \\
&\leq &\frac{\nu }{8}\left\Vert \varphi (t)\right\Vert _{V}^{2}+C(\left\Vert
\sigma (t)\right\Vert _{H}^{2}+\left\Vert \varphi (t)\right\Vert _{V^{\prime
}}^{2}),
\end{eqnarray*}%
where we have used the continuous embedding $H\subset V^{\prime }$ and the
property $\left\Vert A^{-1}\varphi \right\Vert _{V}\leq C\left\Vert \varphi
\right\Vert _{V^{\prime }}.$ Next, by using the assumptions $\varphi ^{i}\in
L^{\infty }(0,T;V)$ we have
\begin{eqnarray*}
&&\int_{\Omega }((\varphi ^{1})^{3}-(\varphi ^{2})^{3})A^{-1}(\varphi
^{1}-\varphi ^{2})dx=\int_{\Omega }\varphi \left( (\varphi ^{1})^{2}+\varphi
^{1}\varphi ^{2}+(\varphi ^{2})^{2}\right) A^{-1}\varphi dx \\
&\leq &C\left\Vert \varphi (t)\right\Vert _{L^{4}(\Omega )}(\left\Vert
(\varphi ^{1})^{2}(t)\right\Vert _{L^{2}(\Omega )}+\left\Vert (\varphi
^{2})^{2}(t)\right\Vert _{L^{2}(\Omega )})\left\Vert A^{-1}\varphi
(t)\right\Vert _{L^{4}(\Omega )} \\
&\leq &C\left\Vert \varphi (t)\right\Vert _{V}(\left\Vert \varphi
^{1}(t)\right\Vert _{L^{4}(\Omega )}^{2}+\left\Vert \varphi
^{2}(t)\right\Vert _{L^{4}(\Omega )}^{2})\left\Vert A^{-1}\varphi
(t)\right\Vert _{V} \\
&\leq &\frac{\nu }{8}\left\Vert \varphi (t)\right\Vert _{V}^{2}+C(\left\Vert
\varphi ^{1}(t)\right\Vert _{V}^{4}+\left\Vert \varphi ^{2}(t)\right\Vert
_{V}^{4})\left\Vert \varphi (t)\right\Vert _{V^{\prime }}^{2} \\
&\leq &\frac{\nu }{8}\left\Vert \varphi (t)\right\Vert _{V}^{2}+C\left\Vert
\varphi (t)\right\Vert _{V^{\prime }}^{2},
\end{eqnarray*}%
by H\"{o}lder's inequality and the continuous embedding $V\subset
L^{4}(\Omega )$. The last constant $C$ also accounts for the norms of $%
\varphi ^{i}$ in $L^{\infty }(0,T;V).$ Furthermore,
\begin{equation*}
\int_{\Omega }(l-1-\nu )\varphi A^{-1}\varphi dx\leq \left\vert l-1-\nu
\right\vert \left\Vert \varphi (t)\right\Vert _{V^{\prime }}\left\Vert
A^{-1}\varphi (t)\right\Vert _{V}\leq C\left\Vert \varphi (t)\right\Vert
_{V^{\prime }}^{2}
\end{equation*}%
and finally%
\begin{eqnarray*}
&&\left\langle 1_{\omega }^{\ast }v(t),A^{-1}\varphi (t)\right\rangle
_{V^{\prime },V}\leq C\left\Vert v(t)\right\Vert _{V^{\prime }}\left\Vert
A^{-1}\varphi (t)\right\Vert _{V} \\
&\leq &C(\left\Vert \varphi (t)\right\Vert _{D(A^{1/2})}+\left\Vert \sigma
(t)\right\Vert _{D(A^{/14})})\left\Vert \varphi (t)\right\Vert _{V^{\prime }}
\\
&\leq &C(\left\Vert \varphi (t)\right\Vert _{V}+\left\Vert \sigma
(t)\right\Vert _{V})\left\Vert \varphi (t)\right\Vert _{V^{\prime }}\leq
\frac{\nu }{8}\left\Vert \varphi (t)\right\Vert _{V}^{2}+\frac{\lambda }{8}%
\left\Vert \sigma (t)\right\Vert _{V}^{2}+C_{\lambda }\left\Vert \varphi
(t)\right\Vert _{V^{\prime }}^{2}
\end{eqnarray*}%
with the help of (\ref{6p})-(\ref{8p}). Note that in the last bound we have
already used the parameter $\lambda >0.$ Now, we detail the second
computation which makes use of the parameter $\lambda .$ We have
\begin{eqnarray}
\frac{\lambda }{2}\frac{d}{dt}\left\Vert \sigma (t)\right\Vert
_{H}^{2}+\lambda \left\Vert \sigma (t)\right\Vert _{V}^{2} &\leq &\gamma
\lambda \left\langle A\varphi (t),\sigma (t)\right\rangle _{V^{\prime
},V}+\lambda \left\Vert \sigma (t)\right\Vert _{H}^{2}  \label{12p} \\
&&+\gamma \lambda \int_{\Omega }\varphi \sigma dx+\lambda \left\langle
1_{\omega }^{\ast }u(t),\sigma (t)\right\rangle _{V^{\prime },V}  \notag
\end{eqnarray}%
and treat the terms on the right-hand side. Thus, we see that%
\begin{equation*}
\gamma \lambda \left\langle A\varphi (t),\sigma (t)\right\rangle _{V^{\prime
},V}\leq \gamma \lambda \left\Vert \varphi (t)\right\Vert _{V}\left\Vert
\sigma (t)\right\Vert _{V}\leq \frac{\lambda }{8}\left\Vert \sigma
(t)\right\Vert _{V}^{2}+2\gamma ^{2}\lambda \left\Vert \varphi
(t)\right\Vert _{V}^{2},
\end{equation*}%
then%
\begin{equation*}
\gamma \lambda \int_{\Omega }\varphi \sigma dx\leq \gamma \lambda \left\Vert
\varphi (t)\right\Vert _{V^{\prime }}\left\Vert \sigma (t)\right\Vert
_{V}\leq \frac{\lambda }{8}\left\Vert \sigma (t)\right\Vert
_{V}^{2}+C_{\lambda }\left\Vert \varphi (t)\right\Vert _{V^{\prime }}^{2}
\end{equation*}%
and using (\ref{54+0}) and $V=D(A^{1/2})$ we get finally%
\begin{align*}
&\quad\ \lambda \left\langle 1_{\omega }^{\ast }u(t),\sigma (t)\right\rangle
_{V^{\prime },V}\leq C\lambda \left\Vert u(t)\right\Vert _{V^{\prime
}}\left\Vert \sigma (t)\right\Vert _{V}\leq C_{1}\lambda (\left\Vert \varphi
(t)\right\Vert _{D(A^{1/2})}+\left\Vert \sigma (t)\right\Vert
_{D(A^{1/4})})\left\Vert \sigma (t)\right\Vert _{V} \\
&\leq C_{2}\lambda \left\Vert \varphi (t)\right\Vert _{V}\left\Vert \sigma
(t)\right\Vert _{V}+C_{3}\lambda \left\Vert \sigma (t)\right\Vert
_{H}^{1/2}\left\Vert \sigma (t)\right\Vert _{V}^{3/2} \\
&\leq \frac{\lambda }{8}\left\Vert \sigma (t)\right\Vert
_{V}^{2}+C_{4}\lambda \left\Vert \varphi (t)\right\Vert _{V}^{2}+C_{\lambda
}\left\Vert \sigma (t)\right\Vert _{H}^{2}.
\end{align*}%
Note that the Young inequality with two pairs of exponents $(2,2)$ and $(4,%
\frac{4}{3})$ has been used.

Next, we sum (\ref{11p}) and (\ref{12p}) by observing that the last integral
on the left-hand side of (\ref{11p}) is nonnegative and taking into account
the estimates in all terms. We infer that
\begin{eqnarray*}
&&\frac{1}{2}\frac{d}{dt}\left( \left\Vert \varphi (t)\right\Vert
_{V^{\prime }}^{2}+\lambda \left\Vert \sigma (t)\right\Vert _{H}^{2}\right)
+\left( \frac{\nu }{2}-(2\gamma ^{2}+C_{4})\lambda \right) \left\Vert
\varphi (t)\right\Vert _{V}^{2}+\frac{\lambda }{2}\left\Vert \sigma
(t)\right\Vert _{V}^{2} \\
&\leq &C\left\Vert \varphi (t)\right\Vert _{V^{\prime }}^{2}+C\left\Vert
\sigma (t)\right\Vert _{H}^{2}+C_{\lambda }\left\Vert \varphi (t)\right\Vert
_{V^{\prime }}^{2}+C_{\lambda }\left\Vert \sigma (t)\right\Vert _{H}^{2},%
\text{ a.e. }t\in (0,T).
\end{eqnarray*}%
At this point, we can choose
\begin{equation*}
\lambda =\frac{\nu }{4(2\gamma ^{2}+C_{4})},
\end{equation*}%
then integrate from $0$ to $t$ and apply the Gronwall lemma. Thus we deduce
that $\varphi =0,$ $\sigma =0,$ whence uniqueness follows.

\subparagraph{Continuation of the existence proof on $(0,\infty ).$}

At the end of these two steps, relying on the existence and uniqueness of
the solution on $[0,T],$ with $T$ arbitrary, we show that (\ref{82+1}) has a
unique solution. We recall that $r_{1}$ is independent of $T.$ For any $%
r\leq r_{1}$ (see (\ref{136})) we introduce
\begin{eqnarray}
S_{\infty } &=&\left\{ (y,z)\in L^{2}(0,\infty ;H\times H);\text{ }%
\sup_{t\in (0,\infty )}\left( \left\Vert y(t)\right\Vert
_{D(A^{1/2})}^{2}+\left\Vert z(t)\right\Vert _{D(A^{1/4})}^{2}\right) \right.
\label{S-inf} \\
&&\left. \text{ \ \ \ \ \ \ \ \ \ \ }+\int_{0}^{\infty }\left( \left\Vert
A^{3/2}y(t)\right\Vert _{H}^{2}+\left\Vert A^{3/4}z(t)\right\Vert
_{H}^{2}\right) dt\leq r^{2}\right\} .  \notag
\end{eqnarray}%
By assuming (\ref{130-2}) and (\ref{ro}) we show that there exists a unique
solution on $[0,\infty )$ which also belongs to $S_{\infty }.$

Consider the functions $(y,z):[0,\infty )\rightarrow D(A^{3/2-\varepsilon
})\times H,$ defined by
\begin{equation*}
(y(t),z(t))=(y_{T}(t),z_{T}(t)),\text{ for any }t\in \lbrack 0,T],
\end{equation*}%
where $(y_{T}(t),z_{T}(t))$ denotes here the solution on $[0,T]$ constructed
at Steps 1 and 2. By the uniqueness proof, $(y_{T}(t),z_{T}(t))=(y_{T^{%
\prime }}(t),z_{T^{\prime }}(t))$ on $[0,T]\subset \lbrack 0,T^{\prime }]$
and so $(y,z)$ is well defined. Moreover, by the first part of the proof,
under the assumption (\ref{ro}) it follows that $(y,z)\in S_{\infty }$ and
it is the solution to (\ref{82+1}) satisfying (\ref{y-z}).

\subparagraph{Step 3.}

To prove the stabilization result we multiply equation (\ref{82+1}) by $%
R(y(t),z(t))$ scalarly in $H\times H.$ Since $R$ is symmetric as an
unbounded operator in $H\times H$ and the Riccati equation (\ref{67}) and
the estimate (\ref{68}) hold, we get%
\begin{eqnarray}
&&\frac{1}{2}\frac{d}{dt}(R(y(t),z(t)),(y(t),z(t)))_{H\times H}  \label{91}
\\
&&+\frac{1}{2}\left( \left\Vert A^{3/2}y(t)\right\Vert _{H}^{2}+\left\Vert
A^{3/4}z(t)\right\Vert _{H}^{2}+\left\Vert B^{\ast }R(y(t),z(t))\right\Vert
_{\mathbb{R}^{N}}^{2}\right)  \notag \\
&\leq &\left\Vert \mathcal{G}(y(t))\right\Vert _{H\times H}\left\Vert
R(y(t),z(t))\right\Vert _{H\times H}\leq \left\Vert G(y(t))\right\Vert
_{H}\left\Vert R(y(t),z(t))\right\Vert _{H\times H}  \notag
\end{eqnarray}%
a.e. $t\in (0,T).$ This relation is used to compute an estimate for the norm
$\left\Vert (y(t),z(t))\right\Vert _{H\times H}$ using (\ref{68}), rewritten
for convenience in the following way
\begin{equation}
\left\Vert R(y(t),z(t))\right\Vert _{H\times H}\leq C_{R}\left\Vert
(y(t),z(t))\right\Vert _{D(A)\times D(A^{1/2})}\leq C_{R}\left( \left\Vert
Ay(t)\right\Vert _{H}+\left\Vert A^{1/2}z(t)\right\Vert _{H}\right) .
\label{91+1}
\end{equation}%
Now, we estimate the right-hand side $RHS$ of (\ref{91}) and for simplicity,
we shall write $(y,z)$ instead of $(y(t),z(t)).$ By recalling (\ref{82+3})
we see that
\begin{equation}
RHS\leq C_{R}\sum_{j=1}^{7}\left\Vert I_{j}(y(t))\right\Vert _{H}\left(
\left\Vert Ay(t)\right\Vert _{H}+\left\Vert A^{1/2}z(t)\right\Vert
_{H}\right)  \label{rhs-0}
\end{equation}%
and we treat each term $RHS_{j}$ of the sum, separately, as done in (\ref{I1}%
)-(\ref{I7}).

We have
\begin{eqnarray*}
&&RHS_{1}=\left\Vert I_{1}\right\Vert _{H}\left\Vert R(y,z)\right\Vert
_{H\times H}\leq C\left\Vert A^{3/2}y\right\Vert _{H}^{(2\alpha _{1}+\alpha
_{2}-1)/2}\left\Vert A^{1/2}y\right\Vert _{H}^{(7-(2\alpha _{1}+\alpha
_{2}))/2}\left\Vert Ay\right\Vert _{H} \\
&&+C\left\Vert A^{3/2}y\right\Vert _{H}^{(2\alpha _{1}+\alpha
_{2}-1)/2}\left\Vert A^{1/2}y\right\Vert _{H}^{(7-(2\alpha _{1}+\alpha
_{2}))/2}\left\Vert A^{1/2}z\right\Vert _{H} \\
&\leq &C\left\Vert A^{3/2}y\right\Vert _{H}^{(2\alpha _{1}+\alpha
_{2}+1)/2}\left\Vert A^{1/2}y\right\Vert _{H}^{(7-(2\alpha _{1}+\alpha
_{2}))/2} \\
&&+C\left\Vert A^{3/2}y\right\Vert _{H}^{(2\alpha _{1}+\alpha
_{2}-1)/2}\left\Vert A^{1/2}y\right\Vert _{H}^{(7-(2\alpha _{1}+\alpha
_{2}))/4}\left\Vert A^{1/2}y\right\Vert _{H}^{(7-(2\alpha _{1}+\alpha
_{2}))/4}\left\Vert A^{1/2}z\right\Vert _{H}.
\end{eqnarray*}%
We used the relation (\ref{54+1}). We should take into account that the
calculations before are true for $\alpha _{1}\geq \frac{9}{8},$ $\alpha
_{2}\geq \frac{3}{4}$ (see (\ref{I1}))$.$ By choosing $2\alpha _{1}+\alpha
_{2}=3$ (e.g., for $\alpha _{1}=\frac{9}{8}$ and $\alpha _{2}=\frac{3}{4}$)
we obtain%
\begin{eqnarray}
RHS_{1} &\leq &C\left( \left\Vert A^{3/2}y\right\Vert _{H}^{2}\left\Vert
A^{1/2}y\right\Vert _{H}^{2}+\left\Vert A^{3/2}y\right\Vert _{H}\left\Vert
A^{1/2}y\right\Vert _{H}\left\Vert A^{1/2}y\right\Vert _{H}\left\Vert
A^{1/2}z\right\Vert _{H}\right)  \notag \\
&\leq &C\left( \left\Vert A^{3/2}y\right\Vert _{H}^{2}\left\Vert
A^{1/2}y\right\Vert _{H}^{2}+\left\Vert A^{3/2}y\right\Vert
_{H}^{2}\left\Vert A^{1/2}y\right\Vert _{H}^{2}+\left\Vert
A^{1/2}y\right\Vert _{H}^{2}\left\Vert A^{1/2}z\right\Vert _{H}^{2}\right)
\notag \\
&\leq &C\left( \left\Vert A^{3/2}y\right\Vert _{H}^{2}+\left\Vert
A^{1/2}z\right\Vert _{H}^{2}\right) \left\Vert A^{1/2}y\right\Vert _{H}^{2}.
\label{rhs1}
\end{eqnarray}%
For $RHS_{2}$ we have
\begin{eqnarray}
&&RHS_{2}=\left\Vert I_{2}\right\Vert _{H}\left\Vert R(y,z)\right\Vert
_{H\times H}  \notag \\
&\leq &C\left\Vert A^{3/2}y\right\Vert _{H}^{(\alpha _{1}+2\alpha
_{2}-1)/2}\left\Vert A^{1/2}y\right\Vert _{H}^{(7-(\alpha _{1}+2\alpha
_{2}))/2}\left( \left\Vert Ay\right\Vert _{H}+\left\Vert A^{1/2}z\right\Vert
_{H}\right)  \notag \\
&\leq &C\left\Vert A^{3/2}y\right\Vert _{H}^{(\alpha _{1}+2\alpha
_{2}+1)/2}\left\Vert A^{1/2}y\right\Vert _{H}^{(7-(\alpha _{1}+2\alpha
_{2}))/2}  \notag \\
&&+C\left\Vert A^{3/2}y\right\Vert _{H}^{(\alpha _{1}+2\alpha
_{2}-1)/2}\left\Vert A^{1/2}y\right\Vert _{H}^{(7-(\alpha _{1}+2\alpha
_{2}))/4}\left\Vert A^{1/2}y\right\Vert _{H}^{(7-(\alpha _{1}+2\alpha
_{2}))/4}\left\Vert A^{1/2}z\right\Vert _{H}  \notag \\
&\leq &C\left( \left\Vert A^{3/2}y\right\Vert _{H}^{2}+\left\Vert
A^{1/2}z\right\Vert _{H}^{2}\right) \left\Vert A^{1/2}y\right\Vert _{H}^{2}
\label{rhs2}
\end{eqnarray}%
for $\alpha _{1}+2\alpha _{2}=3$ $\left( \alpha _{1}=\frac{3}{4},\text{ }%
\alpha _{2}=\frac{9}{8}\right) .$

In the same way we get the other necessary estimates. We have
\begin{equation}
RHS_{3}\leq C\left\Vert \nabla \varphi _{\infty }\right\Vert _{\infty
}\left( \left\Vert A^{3/2}y\right\Vert _{H}^{2}+\left\Vert
A^{1/2}z\right\Vert _{H}^{2}\right) \left\Vert A^{1/2}y\right\Vert _{H},
\label{rhs3}
\end{equation}%
\begin{equation}
RHS_{4}\leq C\left\Vert \Delta \varphi _{\infty }\right\Vert _{\infty
}\left( \left\Vert A^{3/2}y\right\Vert _{H}^{2}+\left\Vert
A^{1/2}z\right\Vert _{H}^{2}\right) \left\Vert A^{1/2}y\right\Vert _{H},
\label{rhs4}
\end{equation}%
\begin{equation}
RHS_{5}\leq C\left\Vert \varphi _{\infty }\right\Vert _{\infty }\left(
\left\Vert A^{3/2}y\right\Vert _{H}^{2}+\left\Vert A^{1/2}z\right\Vert
_{H}^{2}\right) \left\Vert A^{1/2}y\right\Vert _{H},  \label{rhs5}
\end{equation}%
\begin{equation}
RHS_{6}\leq C\left\Vert \varphi _{\infty }\right\Vert _{\infty }\left(
\left\Vert A^{3/2}y\right\Vert _{H}^{2}+\left\Vert A^{1/2}z\right\Vert
_{H}^{2}\right) \left\Vert A^{1/2}y\right\Vert _{H},  \label{rhs6}
\end{equation}%
\begin{equation}
RHS_{7}\leq C\left\Vert g_{\infty }\right\Vert _{2,\infty }\left( \left\Vert
A^{3/2}y\right\Vert _{H}^{2}+\left\Vert A^{1/2}z\right\Vert _{H}^{2}\right) .
\label{rhs7}
\end{equation}%
Therefore, owing to (\ref{54+1}), we see that (\ref{rhs-0}) yields
\begin{equation*}
RHS\leq \widehat{C}\left( \left\Vert A^{3/2}y\right\Vert _{H}^{2}+\left\Vert
A^{3/4}z\right\Vert _{H}^{2}\right) \left( \left\Vert A^{1/2}y\right\Vert
_{H}^{2}+\left\Vert \varphi _{\infty }\right\Vert _{2,\infty }\left\Vert
A^{1/2}y\right\Vert _{H}+\left\Vert g_{\infty }\right\Vert _{2,\infty
}\right) ,
\end{equation*}%
where we have marked $\widehat{C}$ by using a special symbol for future
convenience. We stress that $\widehat{C}$ depends on $\Omega .$

Coming back to (\ref{91}), by ignoring a nonnegative term on the left-hand
side and making use of the fact that $(y,z)\in S_{\infty }$ (see (\ref{S-inf}%
)) we conclude that for a.e. $t$%
\begin{eqnarray*}
&&\frac{d}{dt}(R(y(t),z(t)),(y(t),z(t)))_{H\times H}+\left\Vert
A^{3/2}y\right\Vert _{H}^{2}+\left\Vert A^{3/4}z\right\Vert _{H}^{2} \\
&\leq &\widehat{C}\left( \left\Vert A^{3/2}y\right\Vert _{H}^{2}+\left\Vert
A^{3/4}z\right\Vert _{H}^{2}\right) \left( r^{2}+\left\Vert \varphi _{\infty
}\right\Vert _{2,\infty }r+\left\Vert g_{\infty }\right\Vert _{2,\infty
}\right) .
\end{eqnarray*}%
We recall that the initial datum $(y_{0},z_{0}),$ the parameter $r$ and the
target $\varphi _{\infty }$ are already subject to some restrictions: see (%
\ref{130-2})-(\ref{hi-0}), where $C_{P}$ and $C$ are the precise constants
occurring in (\ref{130}) and (\ref{est-G}). Here, we require something more.
Namely, we impose that%
\begin{equation}
C_{\ast }:=1-\widehat{C}\left( r^{2}+\left\Vert \varphi _{\infty
}\right\Vert _{2,\infty }r+\left\Vert g_{\infty }\right\Vert _{2,\infty
}\right) >0.  \label{C*}
\end{equation}%
If we repeat the argument we have used to obtain the restriction $\chi
_{\infty }\leq \chi _{0}$ of (\ref{hi-0}) we see that (\ref{C*}) is
satisfied whenever%
\begin{equation*}
r\leq r_{2},\text{ \ }r_{2}:=\frac{-\left\Vert \varphi _{\infty }\right\Vert
_{2,\infty }+\sqrt{D_{\infty }^{\prime }}}{2},\text{ \ }D_{\infty }^{\prime
}:=\left\Vert \varphi _{\infty }\right\Vert _{2,\infty }^{2}-4(\left\Vert
g\right\Vert _{2,\infty }-(\widehat{C})^{-1})
\end{equation*}%
by assuming that $\left\Vert g\right\Vert _{2,\infty }\leq (\widehat{C}%
)^{-1}.$ This is true if%
\begin{equation}
\chi _{\infty }\leq \chi _{0}^{\prime \prime }:=\frac{-\left\Vert \varphi
_{\infty }\right\Vert _{\infty }+\sqrt{\left\Vert \varphi _{\infty
}\right\Vert _{\infty }^{2}+4(C_{\Omega }\widehat{C})^{-1}}}{2}
\label{hi-0-2}
\end{equation}%
with $C_{\Omega }$ (the same as before in Step 1) and $\widehat{C}$
depending only on $\Omega .$ Then we set (see (\ref{hi-0}))
\begin{equation*}
\chi _{0}:=\min \left\{ \chi _{0}^{\prime },\chi _{0}^{\prime \prime
}\right\} .
\end{equation*}%
At this point, by assuming $\chi _{\infty }\leq \chi _{0},$ we can fix $%
r=\min \{r_{1},r_{2}\}$ and $\rho $ satisfying (\ref{130-2}).

Under the assumption $\chi _{\infty }\leq \chi _{0}$, the previous
inequality implies that
\begin{equation}
\frac{d}{dt}(R(y(t),z(t)),(y(t),z(t)))_{H\times H}+C_{\ast }\left(
\left\Vert A^{3/2}y(t)\right\Vert _{H}^{2}+\left\Vert A^{3/4}z(t)\right\Vert
_{H}^{2}\right) \leq 0  \label{109}
\end{equation}%
a.e. $t\in (0,\infty ).$

Next, we owe to (\ref{54+1}) and (\ref{69}) and see that%
\begin{equation*}
\left\Vert A^{3/2}y(t)\right\Vert _{H}^{2}+\left\Vert A^{3/4}z(t)\right\Vert
_{H}^{2}\geq c_{0}(R(y(t),z(t)),(y(t),z(t)))_{H\times H}
\end{equation*}%
for some constant $c_{0}>0$ depending on the problem parameters and $%
\left\Vert \varphi _{\infty }\right\Vert _{L^{2}(\Omega )}.$ Therefore we
deduce that
\begin{equation}
\frac{d}{dt}(R(y(t),z(t)),(y(t),z(t)))_{H\times H}+C_{\ast
}c_{0}(R(y(t),z(t)),(y(t),z(t)))_{H\times H}\leq 0,\text{ a.e. }t\in
(0,\infty ),  \label{111}
\end{equation}%
and this immediately implies
\begin{equation}
(R(y(t),z(t)),(y(t),z(t)))_{H\times H}\leq
e^{-2kt}(R(y_{0},z_{0}),(y_{0},z_{0}))_{H\times H}  \label{113}
\end{equation}%
where $k:=\frac{C_{\ast }c_{0}}{2}.$

On account of (\ref{69}) we conclude that
\begin{equation*}
c_{1}\left\Vert (y(t),z(t))\right\Vert _{D(A^{1/2})\times
D(A^{1/4})}^{2}\leq c_{2}e^{-2kt}\left\Vert (y_{0},z_{0})\right\Vert
_{D(A^{1/2})\times D(A^{1/4})}^{2},\text{ a.e. }t>0.
\end{equation*}%
This is nothing but (\ref{83+5}) with a precise constant in front of the
exponential. Hence, we conclude the proof.\hfill $\square $

\medskip

\noindent \textbf{Proof of Theorem 1.1.}\textit{\ }The result given in
Theorem 1.1 follows immediately. Assume (\ref{c3}), that is
\begin{equation*}
\left\Vert \varphi _{0}-\varphi _{\infty }\right\Vert
_{D(A^{1/2})}+\left\Vert \theta _{0}-\theta _{\infty }\right\Vert
_{D(A^{1/4})}\leq \rho
\end{equation*}%
and go back to (\ref{1})-(\ref{4}) by transformation (\ref{8}). Then,
\begin{eqnarray*}
&&\left\Vert y_{0}\right\Vert _{D(A^{1/2})}+\left\Vert z_{0}\right\Vert
_{D(A^{1/4})}=\left\Vert \varphi _{0}-\varphi _{\infty }\right\Vert
_{D(A^{1/2})}+\left\Vert \sigma _{0}-\sigma _{\infty }\right\Vert
_{D(A^{1/4})} \\
&\leq &(\alpha _{0}l_{0}+1)\left\Vert \varphi _{0}-\varphi _{\infty
}\right\Vert _{D(A^{1/2})}+\alpha _{0}\left\Vert \theta _{0}-\theta _{\infty
}\right\Vert _{D(A^{1/4})}\leq \rho \max \{\alpha _{0}l_{0}+1,\alpha
_{0}\}=:\rho _{1}.
\end{eqnarray*}%
By Theorem 3.1, for $\rho _{1}$ small enough, we get (\ref{y-z}) and (\ref%
{83+5}). The latter becomes (\ref{c5}) by using the transformation (\ref{8}).

\section{Appendix}

\setcounter{equation}{0}

\noindent \textbf{Lemma A1 }\textit{The stationary system}\textbf{\ }(\ref%
{17'})\textbf{\ }\textit{has at least a solution }$(\theta _{\infty
},\varphi _{\infty })$\textit{, with }$\theta _{\infty }$\textit{\ constant
and} $\varphi _{\infty }\in H^{4}(\Omega )\subset C^{2}(\overline{\Omega }).$

\medskip

\noindent \textbf{Proof. }It\textbf{\ }is obvious that $\theta _{\infty }$
is constant and $\varphi _{\infty }$ satisfies equation
\begin{eqnarray*}
\Delta (\nu \Delta \varphi _{\infty }-\varphi _{\infty }^{3}+\varphi
_{\infty }) &=&0, \\
\frac{\partial \Delta \varphi _{\infty }}{\partial \nu } &=&\frac{\partial
\varphi _{\infty }}{\partial \nu }=0
\end{eqnarray*}%
whence
\begin{equation}
\nu \Delta \varphi _{\infty }-\varphi _{\infty }^{3}+\varphi _{\infty }=C.
\label{200}
\end{equation}%
This equation has a solution $\varphi _{\infty }\in H^{1}(\Omega )$. The
argument is that (\ref{200}) is a necessary condition for $\varphi _{\infty
} $ to be a minimizer of the functional
\begin{equation*}
\Upsilon (\varphi )=\int_{\Omega }\left( \frac{\nu }{2}\left\vert \nabla
\varphi \right\vert ^{2}+\frac{(\varphi ^{2}-1)^{2}}{4}+C\varphi \right) dx
\end{equation*}%
for $\varphi \in H^{1}(\Omega ).$ On the other hand $\Upsilon $ is l.s.c.
and coercive on $H^{1}(\Omega ),$ whence a minimizer exists.

It follows by (\ref{200}) that $\Delta \varphi _{\infty }\in L^{2}(\Omega )$
and together with $\frac{\partial \varphi _{\infty }}{\partial \nu }=0$ this
leads to $\varphi _{\infty }\in H^{2}(\Omega ).$ Then $\Delta \varphi
_{\infty }^{3}\in L^{2}(\Omega )$ by a simple computation and we deduce that
\begin{equation*}
\nu \Delta \Delta \varphi _{\infty }=\Delta (\varphi _{\infty }^{3}-\varphi
_{\infty })\in L^{2}(\Omega ),
\end{equation*}
whence $\varphi _{\infty }\in H^{4}(\Omega )$ since $\nu \Delta \varphi
_{\infty }=\varphi _{\infty }^{3}-\varphi _{\infty }+C$ satisfies the
homogeneous Neumann condition. Therefore, $\varphi _{\infty }\in C^{2}(%
\overline{\Omega }).$\hfill $\square $

\medskip

\noindent \textbf{Lemma A2 }(Kalman, see \cite{lee-markus})\textbf{. }%
\textit{Let us consider the finite dimensional system }%
\begin{eqnarray*}
X^{\prime }+MX &=&DW,\mbox{ }t\in \lbrack 0,T_{0}] \\
X(0) &=&X_{0}\in \mathbb{R}^{N},
\end{eqnarray*}%
\textit{where }$M$\textit{\ and }$D$\textit{\ are real }$N\times N$\textit{\
matrices}. \textit{Assume that }$D^{\ast }P(t)=0$\textit{\ implies }$P(t)=0$%
\textit{\ for all }$t\in \lbrack 0,T_{0}],$\textit{\ where }$P$\textit{\ is
any solution to the dual backward differential system }$P^{\prime
}(t)-M^{\ast }P(t)=0.$\textit{\ Then, there exists }$W\in L^{2}(0,T_{0};%
\mathbb{R}^{N})$\textit{\ such that }$X(T_{0})=0.$\textit{\ Moreover,}%
\begin{equation}
\int_{0}^{T_{0}}\left\Vert W(t)\right\Vert _{\mathbb{R}^{N}}^{2}dt\leq
C\left\Vert X_{0}\right\Vert _{\mathbb{R}^{N}}^{2}.  \label{36+6}
\end{equation}

\medskip

\noindent \textbf{Lemma A3. }\textit{Let }$(E\subset F\subset E^{\prime })$%
\textit{\ be a variational triplet and let }$L:E\rightarrow E^{\prime }$%
\textit{\ be a linear continuous operator such that }%
\begin{equation*}
\left\langle Lw,w\right\rangle _{E^{\prime },E}\geq C_{1}\left\Vert
w\right\Vert _{E}^{2}-C_{2}\left\Vert w\right\Vert _{F}^{2}.
\end{equation*}%
\textit{Let }$M\in \mathcal{L}(E,F),$\textit{\ that is }%
\begin{equation*}
\left\Vert Mw\right\Vert _{F}\leq C_{3}\left\Vert w\right\Vert _{E}.
\end{equation*}%
\textit{Then, }$L+M$\textit{\ is quasi }$m$\textit{-accretive in }$F\times
F. $

\medskip

\noindent \textbf{Proof. }The operator $\widetilde{L}=L+M$ is continuous
from $E$ to $E^{\prime }$ and
\begin{equation*}
\left\langle \widetilde{L}w,w\right\rangle _{E^{\prime },E}\geq
C_{1}\left\Vert w\right\Vert _{E}^{2}-C_{2}\left\Vert w\right\Vert
_{F}^{2}-C_{3}\left\Vert w\right\Vert _{F}^{2}\geq C_{1}\left\Vert
w\right\Vert _{E}^{2}-C_{4}\left\Vert w\right\Vert _{F}^{2}.
\end{equation*}%
Hence, $\lambda I+\widetilde{L}$ is $m$-accretive on $F\times F$ because
\begin{equation*}
\left( (\lambda I+\widetilde{L})w,w\right) _{F}\geq 0,\mbox{ for }\lambda
\geq \lambda _{0}
\end{equation*}%
and equation
\begin{equation*}
(\lambda I+\widetilde{L})w=f\in F
\end{equation*}%
has a solution for $\lambda $ sufficiently large.

\medskip

\noindent \textbf{Lemma A4 }(\cite{datko}, see \cite{pazy}, p. 116).\textit{%
\ Let }$-\mathcal{D}$\textit{\ be a }$C_{0}$\textit{-semigroup generator on
a Banach space }$\mathcal{X}$\textit{\ and consider the Cauchy problem }%
\begin{eqnarray*}
\frac{dZ}{dt}(t)+\mathcal{D}Z(t) &=&0,\mbox{ \textit{a.e.} }t>0, \\
Z(0) &=&Z_{0}.
\end{eqnarray*}%
\textit{\ If }%
\begin{equation*}
\int_{0}^{\infty }\left\Vert Z(t)\right\Vert _{\mathcal{X}}^{2}dt\leq
C\left\Vert Z_{0}\right\Vert _{\mathcal{X}}^{2},\mbox{ }\forall Z_{0}\in
\mathcal{X},
\end{equation*}%
\textit{\ then}
\begin{equation*}
\left\Vert Z(t)\right\Vert _{\mathcal{X}}^{2}\leq Ce^{-\kappa t}\left\Vert
Z_{0}\right\Vert _{\mathcal{X}}^{2},
\end{equation*}%
\textit{for some} $k>0.$

\mathstrut

\noindent \textbf{Acknowledgments}

\medskip

This research activity has been performed in the framework of the
Italian-Romanian project \textquotedblleft Nonlinear partial differential
equations (PDE) with applications in modeling cell growth, chemotaxis and
phase transition\textquotedblright\ financed by the Italian CNR and the
Romanian Academy. The present paper also benefits from the support of the
UEFISCDI project PNII-ID-PCE-2011-3-0027 for VB and GM. A partial support
from the GNAMPA (Gruppo Nazionale per l'Analisi Matematica, la Probabilit%
\`{a} e loro Applicazioni) of INDAM (Istituto Nazionale di Alta Matematica)
and the IMATI -- C.N.R. Pavia is gratefully acknowledged by PC and GG.

\end{document}